\documentclass[12pt]{amsart}
\usepackage{amsfonts,amssymb,amscd,amsmath}
\usepackage{latexsym}
\usepackage{tocvsec2}

\usepackage[
hyperindex=true,pagebackref=true,bookmarks=true,
colorlinks=true,linkcolor=blue,citecolor=red]
{hyperref}
\def\~{{\rm --}}

\settocdepth{subsubsection}

\begin{document}
\title
[Nonsymmetric difference Whittaker functions]
{Nonsymmetric difference Whittaker functions}
\author
[Ivan Cherednik and Daniel Orr]
{Ivan Cherednik $^\dag$ and Daniel Orr}

\begin{abstract}
Starting with nonsymmetric global difference spherical 
functions, we define and calculate spinor (nonsymmetric) 
global q-Whittaker functions for arbitrary reduced
root systems, which are reproducing kernels of 
the DAHA-Fourier transforms of Nil-DAHA and solutions
of the q-Toda-Dunkl eigenvalue problem. We introduce the
spinor q-Toda-Dunkl operators as limits of the difference
Dunkl operators in DAHA theory under the spinor variant of 
the Ruijsenaars procedure. Their general algebraic theory 
(any reduced root systems) is the key part of this paper, 
based on the new technique of W-spinors and corresponding 
developments in combinatorics of affine root systems.
\end{abstract}

\thanks{$^\dag$  \today\ \ \ Partially supported by NSF grant
DMS--1101535}

\address[I. Cherednik]{Department of Mathematics, UNC
Chapel Hill, North Carolina 27599, USA\\
chered@email.unc.edu}
\address[D. Orr]{Department of Mathematics, UNC
Chapel Hill, North Carolina 27599, USA\\
danorr@email.unc.edu}

\def\bysame{{\bf --- }}
\def\~{{\bf --}}
\newcommand{\comment}[1]{}
\renewcommand{\tilde}{\widetilde}
\renewcommand{\hat}{\widehat}
\newcommand{\dagx}{\hbox{\tiny\mathversion{bold}$\dag$}}
\newcommand{\ddagx}{\hbox{\tiny\mathversion{bold}$\ddag$}}

\renewcommand{\tilde}{\widetilde}
\renewcommand{\hat}{\widehat}

\newcommand{\BR}{{\mathbb R}}
\newcommand{\BQ}{{\mathbb Q}}
\newcommand{\BC}{{\mathbb C}}
\newcommand{\BP}{{\mathbb P}}
\newcommand{\BZ}{{\mathbb Z}}
\newcommand{\BN}{{\mathbb N}}
\newcommand{\BS}{{\mathbb S}}

\newcommand{\cH}{{\mathcal H}}
\newcommand{\cA}{{\mathcal A}}
\newcommand{\cB}{{\mathcal B}}
\newcommand{\ccF}{{\mathfrak F}}
\newcommand{\cD}{{\mathcal D}}
\newcommand{\cL}{{\mathcal L}}
\newcommand{\cF}{{\mathcal F}}
\newcommand{\cP}{{\mathcal P}}
\newcommand{\cX}{{\mathcal X}}
\newcommand{\cY}{{\mathcal Y}}
\newcommand{\cS}{{\mathcal S}}
\newcommand{\cSol}{\hbox{$\mathcal Sol$}}
\newcommand{\cT}{\hbox{$\mathcal T$}}

\newcommand{\Z}{{\mathbb Z}}
\newcommand{\Q}{{\mathbb Q}}
\newcommand{\N}{{\mathbb N}}
\newcommand{\C}{{\mathbb C}}
\newcommand{\R}{{\mathbb R}}

\newcommand{\CH}{{\mathcal H}}
\newcommand{\CA}{{\mathcal A}}

\def\HH{\mbox{${\mathcal H}$\kern-5.2pt${\mathcal H}$}}


\def\der{\partial}
\def\tensor{\otimes}
\def\gam{\gamma} \def\Gam{\Gamma}
\def\del{\delta} \def\Del{\Delta}
\def\kap{\kappa}
\def\lam{\lambda} \def\Lam{\Lambda}
\def\Comp{{\mathbb C}}
\def\sM{{\mathcal M}}

\newtheorem{theorem}{Theorem}[section]
\newtheorem{maintheorem}[theorem]{Main Theorem}
\newtheorem{proposition}[theorem]{Proposition}
\newtheorem{definition}[theorem]{Definition}
\newtheorem{lemma}[theorem]{Lemma}
\newtheorem{corollary}[theorem]{Corollary}
\newtheorem{notation}[theorem]{Notation}
\newtheorem{remark}[theorem]{Remark}
\newtheorem{example}[theorem]{Example}

\newtheorem{theorem }{Theorem}[section]
\newtheorem{maintheorem }[theorem]{Main Theorem}
\newtheorem{proposition }[theorem]{Proposition}
\newtheorem{definition }[theorem]{Definition}
\newtheorem{lemma }[theorem]{Lemma}
\newtheorem{corollary }[theorem]{Corollary}
\newtheorem{notation }[theorem]{Notation}
\newtheorem{remark }[theorem]{Remark}
\newtheorem{example }[theorem]{Example}

\newtheorem{ maintheorem }[theorem]{Main Theorem}
\newtheorem{ theorem}{Theorem}[section]
\newtheorem{ proposition}[theorem]{Proposition}
\newtheorem{ definition}[theorem]{Definition}
\newtheorem{ lemma}[theorem]{Lemma}
\newtheorem{ corollary}[theorem]{Corollary}
\newtheorem{ notation}[theorem]{Notation}
\newtheorem{ remark}[theorem]{Remark}
\newtheorem{ example}[theorem]{Example}

\newtheorem{thm}{Theorem}[section]
\newtheorem{prop}[thm]{Proposition}
\newtheorem{lem}[thm]{Lemma}
\newtheorem{cor}[thm]{Corollary}
\newtheorem{conj}[thm]{Conjecture}
\newtheorem{con}[thm]{Conjecture}
\newtheorem{dfn}[thm]{Definition}
\newtheorem{df}[thm]{Definition}
 \newcommand{\rem}{{\bf Comment.\ }}
 \newcommand{\rmk}{{\bf Comment.\ }}
 \newcommand{\exmp}{{\bf Example.\ }}
 \newcommand{\ex}{{\bf Example.\ }}
 \newcommand{\prob}{{\bf Problem.\ }}

\newtheorem{note}{Note} 
\renewcommand{\thenote}{}
\newtheorem*{acka}{Acknowledgments}
\newtheorem{ack}{Acknowledgments}
\renewcommand{\theack}{}
\renewcommand{\appendixname}{\bf Appendix}
\renewcommand{\proof}{{\em Proof.\ }}

\hyphenation{
ap-pen-dix as-ymp-tot-ic at-trib-uted at-trib-ut-able
Bry-li-n-sky com-mu-ta-tion de-ge-ne-rate
de-riv-a-tive dis-trib-ute equi-vari-ant ex-tra-or-di-nary  
geo-met-ric griev-ance griev-ous grad-ed ho-lo-no-my ho-mo-thetic
in-fin-ite-ly in-fin-i-tes-i-mal Ha-rish Cha-n-dra mul-ti-plic-able 
non-euclid-ean non-iso-mor-phic non-smooth par-a-digm 
par-a-bol-ic pa-rab-o-loid pa-ram-e-trize phe-nom-e-non 
post-script pseu-do-dif-fer-en-tial pseu-do-fi-nite 
qua-drat-ics quad-ra-ture Han-kel rec-tan-gle semi-def-i-nite 
set-up wide-spread Euler-ian Feb-ru-ary Gauss-ian Grothen-dieck 
Hamil-ton-ian Her-mi-t-ian her-mi-t-ian Jan-u-ary 
Japan-ese Ka-shi-wa-ra Kor-te-weg Le-gendre No-vem-ber Rie-mann-ian 
Sep-tem-ber Za-mo-lo-d-chi-kov Kni-zh-nik quan-tum Op-dam
Mac-do-nald Ca-lo-ge-ro Su-ther-land Mo-ser 
Ol-sha-net-sky  Pe-re-lo-mov in-de-pen-dent ope-ra-tors 
cy-clo-to-mic ra-tio-nal de-gen-er-a-tion 
in-ter-est-ing de-for-ma-tions de-for-ma-tion pro-ce-dure 
fol-lows ope-ra-tors  pre-serve suf-fices ap-proach 
for-mu-las con-sider its com-ple-tion cor-re-spond-ing 
au-to-mor-phism be-cause pro-por-tional fi-nal-ly let-ting 
equi-v-a-lence ge-n-er-al-ized Mac-do-nald iden-ti-ties 
cor-re-s-pond sub-dia-grams par-ti-tion na-t-u-ral-ly 
or-dered stan-dard de-for-ma-tion ar-gu-ment com-bined 
sphe-r-i-cal rep-re-sen-ta-tions tri-go-no-me-t-ric
ge-n-er-al-ly speak-ing pri-m-it-ive ir-re-du-cible 
sum-ma-tion  rep-re-sen-ta-tives pro-por-ti-o-na-li-ty
ultra-sphe-ri-cal Ro-gers}

\def\ffor{\quad\hbox{ for }\quad}
\def\wwhen{\quad\hbox{ when }\quad}
\def\wwhere{\quad\hbox{ where }\quad}
\def\aand{\quad\hbox{ and }\quad}
\def\for{\  \hbox{ for } \ }
\def\iif{ \ \hbox{ if } \ }
\def\when{ \ \hbox{ when } \ }
\def\where{\  \hbox{ where } \ }
\def\and{\  \hbox{ and } \ }
\def\and{\  \hbox{ and } \ }
\def\oor{\  \hbox{ or } \ }
\def\proof{{\em Proof. \  }}

\def\equal{\stackrel{\,\mathbf{def}}{= \kern-3pt =}}

\def\la{\lambda}
\def\La{\Lambda}
\def\om{\omega}
\def\Om{\Omega}
\def\Th{\Theta}
\def\th{\theta}
\def\al{\alpha}
\def\be{\beta}
\def\ga{\gamma}
\def\ep{\epsilon}
\def\up{\upsilon}
\def\Up{\Upsilon}
\def\de{\delta}
\def\De{\Delta}
\def\ka{\kappa}
\def\kapp{\hbox{\bf \ae}}
\def\si{\sigma}
\def\Si{\Sigma}
\def\Ga{\Gamma}
\def\ze{\zeta}
\def\io{\iota}
\def\bio{b^\iota}
\def\aio{a^\iota}
\def\twio{\tilde{w}^\iota}
\def\hwio{\hat{w}^\iota}
\def\gio{\g^\iota}
\def\Bio{B^\iota}

\def\del{\delta}
\def\pa{\partial}
\def\vp{\varphi}
\def\ve{\varepsilon}
\def\inf{\infty}

\def\vph{\varphi}
\def\vps{\varpsi}
\def\vPh{\varPhi}
\def\vep{\varepsilon}
\def\vpi{{\varpi}}
\def\vth{{\vartheta}}
\def\vsi{{\varsigma}}
\def\vrh{{\varrho}}

\def\bph{\bar{\phi}}
\def\bsi{\bar{\si}}
\def\bvp{\bar{\varphi}}

\newcommand{\bS}{{\mathbf S}}
\newcommand{\bH}{{\mathbf H}}
\newcommand{\bF}{{\mathbf F}}
\newcommand{\bE}{{\mathbf E}}

\def\tal{\tilde{\alpha}}
\def\tbe{\tilde{\beta}}
\def\tde{\tilde{\delta}}
\def\tpi{\tilde{\pi}}
\def\txi{\tilde{\xi}}
\def\tPi{\tilde{\Pi}}
\def\tPhi{\tilde{\Phi}}
\def\tV{\tilde{V}}
\def\tJ{\tilde{J}}
\def\tla{\tilde{\lambda}}
\def\tga{\tilde{\gamma}}
\def\tGa{\tilde{\Gamma}}
\def\tvs{\tilde{{\varsigma}}}
\def\tu{\tilde{u}}
\def\tU{\tilde{U}}
\def\tw{\widetilde w}
\def\tW{\widetilde W}
\def\tB{\tilde B}
\def\tv{\tilde v}
\def\tV{\tilde V}
\def\tz{\tilde z}
\def\tb{\tilde b}
\def\ta{\tilde a}
\def\tih{\tilde h}
\def\trh{\tilde {\rho}}
\def\tx{\tilde x}
\def\tf{\tilde f}
\def\tg{\tilde g}
\def\tG{\tilde G}
\def\tk{\tilde k}
\def\tl{\tilde l}
\def\tL{\tilde L}
\def\tD{\tilde D}
\def\tR{\tilde R}
\def\tP{\tilde P}
\def\tH{\tilde H}
\def\tp{\tilde p}

\def\hH{\hat{H}}
\def\hh{\hat{h}}
\def\hR{\hat{R}}
\def\hY{\hat{Y}}
\def\hX{\hat{X}}
\def\hP{\hat{P}}
\def\hT{\hat{T}}
\def\hV{\hat{V}}
\def\hG{\hat{G}}
\def\hF{\hat{F}}
\def\hw{\widehat{w}}
\def\hW{\widehat{W}}
\def\hu{\hat{u}}
\def\hs{\hat{s}}
\def\hv{\hat{v}}
\def\hb{\hat{b}}
\def\hB{\widehat{B}}
\def\hze{\hat{\zeta}}
\def\hsi{\hat{\sigma}}
\def\hrh{\hat{\rho}}
\def\hth{\hat{\theta}}
\def\hy{\hat{y}}
\def\hx{\hat{x}}
\def\hz{\hat{z}}
\def\hg{\hat{g}}
\def\he{\hat{e}}
\def\hE{\widehat{E}}

\def\B{\mathbf{B}}
\def\I{\mathbf{I}}
\def\P{\mathbf{P}}
\def\G{\mathbf{G}}
\def\S{\mathbf{S}}
\def\F{\mathbf{F}}
\def\one{\mathbf{1}}
\def\Sn{\mathbf{S}_n}
\def\0{\mathbf{0}}
\def\H{\mathbf{H}}
\def\V{\mathbf{V}}

\def\f{\mathcal{F}}
\def\çF{\mathcal{F}}
\def\o{\mathcal{O}}
\def\t{\mathcal{T}}
\def\r{\mathcal{R}}
\def\l{\mathcal{L}}
\def\m{\mathcal{M}}
\def\k{\mathcal{K}}
\def\n{\mathcal{N}}
\def\d{\mathcal{D}}
\def\p{\mathcal{P}}
\def\cP{\mathcal{P}}
\def\a{\mathcal{A}}
\def\h{\mathcal{H}}
\def\c{\mathcal{C}}
\def\y{\mathcal{Y}}
\def\e{\mathcal{E}}
\def\v{\mathcal{V}}
\def\z{\mathcal{Z}}
\def\x{\mathcal{X}}
\def\s{\mathcal{S}}
\def\g{\mathcal{G}}
\def\u{\mathcal{U}}
\def\w{\mathcal{W}}
\def\i{\mathcal{I}}
\def\j{\mathcal{J}}
\def\b{\mathcal{B}}

\def\lan{\langle}
\def\llb{(\!(}
\def\ran{\rangle}
\def\rrb{)\!)}
 \def\dim{{\hbox{\rm dim}}_{\mathbb C}\,}
\def\lng{\hbox{\rm{\tiny lng}}}
\def\sht{\hbox{\rm{\tiny sht}}}
\def\sph{\hbox{\rm{\tiny sph}}}
\def\inv{\hbox{\rm{\tiny inv}}}

\def\br#1{\langle #1 \rangle}

\def\rank{\hbox{rank}}
\def\gl{\mathfrak{gl}_N}

\newcommand{\Aut}{\operatorname{Aut}}
\newcommand{\Hom}{\operatorname{Hom}}
\newcommand{\End}{\operatorname{End}}
\newcommand{\Ind}{\operatorname{Ind}}
\newcommand{\ad}{\operatorname{ad}}
\newcommand{\pr}{\operatorname{pr}}
\newcommand{\aweyl}{\tilde{\mathbb S}_n}
\newcommand{\hec}{{\mathcal H}^t_n}
\newcommand{\Func}{{\mathcal F}({\mathbb C}^n,{\mathcal H}^t_n)}
\newcommand{\tr}{\operatorname{tr}}
\newcommand{\Out}{\operatorname{Out}}
\newcommand{\Rad}{\operatorname{Rad}}
\newcommand{\Spec}{\operatorname{Spec}}
\newcommand{\id}{\operatorname{id}}
\newcommand{\Int}{\operatorname{Int}}
\newcommand{\ct} {\operatorname{ct}}

\newcommand{\rat}{{\mathbb Q}}
\newcommand{\real}{{\mathbb R}}
\newcommand{\cplx}{{\mathbb C}}
\newcommand{\zint}{{\mathbb Z}}

\newcommand{\sq}{\phantom{1}\hfill$\qed$}
\newcommand{\Rea}{\Re}
\newcommand{\Ima}{\Im}

\newcommand{\st}{\bowtie}
\newcommand{\modd}{\mbox{\,mod\,}}
\newcommand{\lr}{\langle}
\newcommand{\rr}{\rangle}
\newcommand{\eps}{\varepsilon}
\newcommand{\phk}{\phi^{(k)}}
\newcommand{\psk}{\psi^{(k)}}
\newcommand{\Res}{\mbox{Res}\;}
\newcommand{\sgn}{\mbox{sgn}}
\newcommand{\mn} {\left\{ \begin{array}{c}m\\
n\end{array}\right\}}

\def\TT{\mathfrak{T}}
\def\JJ{\mathfrak{J}}
\def\HH{\mathfrak{H}}
\def\FF{\mathfrak{F}}
\def\GG{\mathfrak{G}}
\def\CC{\mathfrak{C}}
\def\LL{\mathfrak{L}}

\def\BB{\mathfrak{B}}
\def\AA{\mathfrak{A}}
\def\ZZ{\mathfrak{Z}}
\def\HH{\hbox{${\mathcal H}$\kern-5.2pt${\mathcal H}$}}
\def\tHH{\widetilde{\HH\ }}

\font\smm=msbm10 at 12pt 
\def\symbol#1{\hbox{\smm #1}}
\def\lsmash{{\symbol n}}
\def\rsmash{{\symbol o}}
\def\#{\sharp}

\font\tenbf=cmbx10
\font\tenrm=cmr10
\font\tenit=cmti10
\font\ninebf=cmbx9
\font\ninerm=cmr9
\font\nineit=cmti9
\font\eightbf=cmbx8
\font\eightrm=cmr8
\font\eightit=cmti8
\font\sevenrm=cmr7
\font\sevenbf=cmbx7


\def\dHH{\mbox{$\dot{\mathcal H}$\kern-5.2pt$\dot{\mathcal H}$}}
\newcommand{\bHH}{\overline{\HH}}
\newcommand{\ts}{\si}
\newcommand{\spin}{\mathrm{Spin}}
\newcommand{\fw}{\mathfrak{w}}
\newcommand{\hgt}{\mathrm{ht}}
\newcommand{\sfs}{\mathsf{s}}
\newcommand{\Red}{\hbox{Red}}
\newcommand{\ord}{\mathrm{ord}}
\newcommand{\dord}{{}^\dag\kern-1pt\ord}
\newcommand{\preceqq}{\preceq\!\!\!\preceq}
\newcommand{\succeqq}{\succeq\!\!\!\succeq}
\newcommand{\precq}{\prec\!\!\!\prec}
\newcommand{\succq}{\succ\!\!\!\succ}
\newcommand{\RE}{R\!E}
\def\fora{\ \hbox{ for all } \ }
\def\unl{\quad\hbox{ unless }\quad}
\newtheorem{conjecture}[theorem]{Conjecture}

\maketitle
\smallskip

{\small
\centerline{{\bf MSC} (2010): 20C08, 22E35, 33D80, 22E66, 20G44}
} 

\renewcommand{\baselinestretch}{1.2}
{\textmd
\tableofcontents
}
\renewcommand{\baselinestretch}{1.0}
\vfill\eject

\renewcommand{\natural}{\wr}

\setcounter{section}{-1}
\setcounter{equation}{0}
\section{\sc Introduction}
This paper is devoted to the theory of nonsymmetric (spinor) difference  
Whittaker functions and the corresponding Toda-Dunkl operators 
in the $q$\~Toda theory for arbitrary irreducible reduced
root systems, 
generalizing the $A_1$\~case 
considered in \cite{ChM} and \cite{ChO1,ChO2}.
Our approach is based on the new technique of $W$\~spinors,
multicomponent functions indexed by the elements of the
nonaffine Weyl group $W$ with the natural action of $W$
on the indices. 

The applications of this technique are deep (despite its
simple definition), with important links to the classical
harmonic analysis on symmetric spaces and the theory of
spherical, Whittaker and Bessel functions. For instance, 
spinors arise in the study of nonsymmetric or singular
symmetric solutions of symmetric systems such as the
Quantum Many-Body Problem; see \cite{C13,O2,ChM}.
However the main applications so far are for Dunkl-type
operators; the Toda-Dunkl operators simply cannot be defined
without $W$\~spinors.

The theory of {\em global} nonsymmetric spherical
functions from \cite{C5}--- the reproducing kernels of DAHA-Fourier
transforms --- is the starting point of this paper.
They
are eigenfunctions of the difference Dunkl operators,
but this is just one of their remarkable properties.
We continue their general theory (e.g.
Proposition \ref{PSIESTAR}) and then define 
global nonsymmetric $q$\~Whittaker functions
as the limits of global spherical functions
using a $W$\~spinor variant of the Ruijsenaars limiting procedure. 
These functions appear certain quadratic-type 
generating functions of the $q$\~Hermite polynomials, denoted by
$\overline{E}_b$ in the paper, for $b$ from
the weight lattice $P$. See Theorem~\ref{GLOBNSWH} and 
Proposition~\ref{OMLIMIT}.  

One can expect such a generating function to be
a series in terms of $q^{(x,b)}\overline{E}_b(\La)$ for 
generic $x\in \C^n$ and all weights $b\in P$.
However, this is not the case. The asymptotic behavior
of such series is inconsistent with the analytic Whittaker 
theory, where the presence of $q^{(x,b)}$ can be expected 
only for (anti-)dominant $b$. 
The ``right" generating function requires $W$\~spinors.
\smallskip

{\bf The symmetric theory.}
Symmetric global $q$\~Whittaker functions from \cite{C10}
solve the $q$\~Toda eigenvalue problem; the $q$\~Toda operators are
due to Ruijsenaars for $GL_n$ \cite{Ru}
and Etingof, Sevostyanov \cite{Et1,Sev} 
for any root systems (via Quantum Groups).
The corresponding {\em global\,} eigenfunctions are given in 
terms of $\overline{E}_{b}(\La)$ for anti-dominant $b\in P_-$ only;
then $\overline{E}_{b}(\La)$ become $W$\~invariant.
They are {\em not\,} $W$\~invariant (symmetric) in terms of $X=q^x$,
just like the classical Whittaker functions. Calling
them ``symmetric" can be confusing, though they are indeed
$W$\~invariant with respect to $\La$.

These functions generalize
the Whittaker functions from the classical harmonic
analysis on symmetric spaces \cite{GW,Wa} and their
$p$\~adic counterparts from \cite{CS} given by 
Shintani-type formulas. $W$\~spinors are not needed in their 
definition (they are functions, not spinors). 
\smallskip

It is necessary to mention a connection with the
$q$\~Whittaker functions obtained in \cite{GiL} from the 
quantum $K$\~theory of the flag varieties; 
see also \cite{GLO1}, \cite{ChO1,ChO2} and \cite{BeF,BF}.
Establishing the relation to our global ``symmetric"
$q$\~Whittaker functions is essentially equivalent to the 
theory of Harish-Chandra decompositions of 
$q,t$\~spherical and
$q$\~Whittaker functions as weighted sums (over $W$) 
of their asymptotic expansions; see \cite{HC,C10,Sto2,ChO1}. 

It is essential here (and in other geometric applications)
that the $q$\~Hermite polynomials coincide with 
the level-one Demazure characters for affine Kac-Moody algebras
(for all weights, not only dominant), due to \cite{San,Ion1}.
Thus, the present paper solves the algebraic-geometric problem 
of finding the generating function for {\em all} level-one
Demazure characters; the answer is the spinor $q$\~Whittaker function.
\smallskip

We note that the theory of global 
$q$\~functions is actually very algebraic (in contrast
to the classical differential theory), including 
special algebraic techniques in the difference 
Harish-Chandra theory. The Harish-Chandra-type theory of
asymptotic decompositions for 
nonsymmetric (spinor) global spherical and Whittaker functions,
including the $p$\~adic limit ($q\to 0$), will be a subject of our 
further paper(s). The theory becomes even more 
algebraic in the nonsymmetric setting due to the use of 
DAHA intertwining operators,
the main tool in the theory of nonsymmetric Macdonald polynomials
and its variants/applications. 
\smallskip

{\bf The main results.}
The construction of Toda-Dunkl operators is the key result
of this paper. The global spinor Whittaker functions are
their eigenfunctions, though we can avoid this fundamental
connection in the theory of these operators. We define
the operator spinor Ruijsenaars-type procedure,
which results in a completely algebraic (though involved) theory
of the Toda-Dunkl operators. 

We give two related approaches to calculating these operators.
They are based on Parts A and B of 
Proposition~\ref{prop-ord-bound} and Lemma~\ref{lem-ord-est}.
Using Part B enables one to obtain arbitrary Toda-Dunkl
operators, but not in a very explicit way.
Part A provides the formulas for certain basic
Toda-Dunkl operators, including those for minuscule weights 
(which are involved even for $A_n$ --- see Section~\ref{examples} 
for some examples), which are then used to calculate {\em all\,} 
operators; see Proposition \ref{SpinDunkl}.
\medskip

The formula for the global spinor Whittaker function 
for this function and its interpretation via
the DAHA-Fourier transform is the second key result of this paper.
This provides the most direct approach to
justification of the existence of Toda-Dunkl operators,
however inconvenient for clarifying their structure.
This is in sharp 
contrast with previously known families of Dunkl operators, 
where such operators directly resulted from 
the existence of the polynomial representation
of the corresponding Double Affine Hecke Algebras, DAHA.

The construction of Toda-Dunkl operators 
in \cite{ChM} was a 
surprising development not expected by specialists. The  
families of Dunkl operators known at that time served only 
$W$\~invariant families of operators; QMBP, 
the Quantum Many-Body Problem (also called the Heckman-Opdam system in
the differential setting),
is a major and the most universal example.
In contrast to QMBP, the Toda operators are not symmetric, which makes 
their theory very different. Only the case of $A_1$ was considered in 
\cite{ChM} and its continuation
\cite{ChO1,ChO2}; it was not clear after these papers how to generalize 
the formulas for the one-dimensional spinor Toda-Dunkl operators 
obtained there. 

The general theory of nonsymmetric Whittaker functions 
follows essentially that for $A_1$, though a 
significant development of the theory of nonsymmetric Hermite 
polynomials was necessary (as was expected). 
However, the intrinsic theory of spinor Toda-Dunkl operators
(without using the global $q$\~Whittaker 
functions) required new tools, 
which deserve thorough analysis and may
have applications beyond our paper. 
\vskip 0.2cm

{\bf Perspectives.} 
We consider this paper a major
step toward the general theory (for arbitrary root systems)
in the following directions:

a) {\em The theory of pseudo-polynomial representations
of Nil-DAHA\,}, which will explain the algebraic origins of the
Toda-Dunkl operators. They are {\em not\,} induced representations
of nil-DAHA in contrast to the polynomial DAHA modules in all other 
theories. However, as in \cite{ChO2} (the case of $A_1$), they are 
induced modules of the {\em core subalgebras\,} of Nil-DAHA; 
the Toda-Dunkl operators naturally act there. 
This construction is part of {\em the theory of canonical-crystal 
bases of Nil-DAHA\,} (in process now), which is expected to have
significant applications in the representation theory of 
DAHA and beyond. 
The theory of canonical-crystal bases was a major development in the 
theory of {\em quantum groups}, closely related to 
{\em Kazhdan-Lusztig polynomials} (missing in the DAHA theory so far)
and {\em cluster algebras}.
\vskip 0.2cm

b) {\em The analytic theory of global nonsymmetric $q,t$\~spherical 
and $q$\~Whittaker functions}, 
including the Harish-Chandra theory of their decompositions in
terms of the asymptotic expansions. This theory requires
the development of analytic techniques for dealing with $W$\~spinors.
The first application (already reached) is a  
nonsymmetric generalization of the existence of {\em the asymptotic 
expansions of global spherical functions\,} from \cite{Sto2} 
(see \cite{ChO1} for the case of $A_1$), which includes the global 
$q$\~Whittaker functions
as well. The nonsymmetric methods significantly simplify 
considerations here, similar to their role at all other levels
of the theory of DAHA and Macdonald polynomials.
\vskip 0.2cm

c) {\em Applications of ``symmetric" $q$\~Whittaker functions},
including the {\em $K$\~theory of flag varieties} \cite{GiL},
the theory of {\em affine flag varieties\,} \cite{BF},\ the theory of 
Demazure characters and {\em local and global Weyl modules\,}
in the Kac-Moody theory, {\em Rogers-Ramanujan identities\,}
(see \cite{ChF} and the references there), 
{\em $q$\~Whittaker processes\,} (random discrete polymers, see
\cite{BC}) and, presumably, the quantum Langlands program. 
We expect these directions to be 
enriched by our theory of spinor $q$\~Whittaker functions and 
Toda-Dunkl operators.
One of the applications (last but not least), which was outlined in
\cite{ChM} for $A_1$, will be the theory of ``nonsymmetric"
{\em Matsumoto-type  $p$\~adic Whittaker functions\,} in the classical
$p$\~adic harmonic analysis.
\medskip

{\bf Acknowledgments.}
The first author thanks Tetsuji Miwa and the Mathematics
Department of Kyoto University for the invitation and hospitality. 
We thank Eric Opdam for useful remarks and Evgeny Feigin for
his help with Conjecture \ref{CONJDOMINB}. The second author 
thanks the organizers of the 5th Southeast Lie Theory 
Workshop, where the results of this paper were reported.
\medskip

\setcounter{equation}{0}
\section{\sc Double Hecke algebra}
Let $R=\{\al\} \subset \R^n$ be a root system of type
$A,B,...,F,G$
with respect to a Euclidean form $(z,z')$ on $\R^n
\ni z,z'$,
$W$ the {\em Weyl group}
generated by the reflections $s_\al$,
$R_{+}$ the set of positive  roots ($R_-=-R_+$)
corresponding to fixed simple
roots $\al_1,...,\al_n,$
$\Ga$ the Dynkin diagram
with $\{\al_i , 1 \le i \le n\}$ as the vertices.
Accordingly,
$$R^\vee=\{\al^\vee =2\al/(\al,\al)\}.$$

The root lattice and the weight lattice are:
\begin{align}
&Q=\oplus^n_{i=1}\Z \al_i \subset P=\oplus^n_{i=1}\Z \om_i,
\notag
\end{align}
where $\{\om_i\}$ are fundamental weights:
$(\om_i,\al_j^\vee)=\de_{ij}$ for the
simple coroots $\al_i^\vee.$
Replacing $\Z$ by $\Z_{\pm}=\{m\in\Z, \pm m\ge 0\}$ we obtain
$Q_\pm, P_\pm.$
Here and further, see \cite{Bo}.

The form will be normalized
by the condition  $(\al,\al)=2$ for the
{\em short} roots in this paper.
Thus,
$$
\nu_R=\{\nu_\al\equal (\al,\al)/2,\al\in R\}
\text{ can be either $\{1\}, \{1,2\},$ or $\{1,3\}.$ }
$$
We will use the notation $\nu_{\lng}$ for long roots 
and $\nu_{\sht}=1$ for short roots.

The normalization leads to the inclusions
$Q\subset Q^\vee,  P\subset P^\vee,$ where $P^\vee$ is
generated by the fundamental coweights $\{\om_i^\vee\}$
dual to $\{\al_i\}$.

We set
$\nu_i\ =\ \nu_{\al_i}$
and
\begin{align}\label{partialrho}
&\rho_\nu\equal \frac{1}{2}\sum_{\nu_{\al}=\nu} \al \ =
\ \sum_{\nu_i=\nu}  \om_i, \hbox{\ where\ } \al\in R_+,
\ \nu\in\nu_R.
\end{align}
Accordingly, $(\rho_\nu,\al_i^\vee)=1$ for $\nu_i=\nu.$
Together with $\rho_\nu$, we will also use the notation 
$\rho_\diamond$ for $\diamond=\lng,\sht$.
\medskip

\subsection{\bf Affine Weyl group}
The vectors $\ \tal=[\al,\nu_\al j] \in
\R^n\times \R \subset \R^{n+1}$
for $\al \in R, j \in \Z $ form the
{\em affine root system}
$\tR \supset R$ ($z\in \R^n$ are identified with $ [z,0]$).
We add $\al_0 \equal [-\vth,1]$ to the simple
roots for the {\em maximal short root} $\vth\in R_+$.
It is also the {\em maximal positive
coroot} because of the choice of normalization.

The corresponding set
$\tR_+$ of positive roots equals 
$R_+\cup \{[\al,\nu_\al j] ,\ \al\in R, \ j > 0\}$.
Indeed, any positive affine root $[\al,\nu_\al j]$
is a linear combination of
$\{\al_i , 0\le i\le n\}$
with coefficients from $\Z_+$.

We complete the Dynkin diagram $\Ga$ of $R$
by adding $\al_0$ ($-\vth$, to be more
exact); it is called {\em affine Dynkin diagram}
$\tGa$. One can obtain it from the
completed Dynkin diagram from \cite{Bo} for
the {\em dual system}
$R^\vee$ by reversing all arrows.
The number of laces between $\al_i$ and
$\al_j$ in $\tGa$ will be denoted by $m_{ij}.$

The set of indices of the images of $\al_0$ by all
the automorphisms of $\tGa$ will be denoted by $O$
($O=\{0\} \for E_8,F_4,G_2$). Let $O'=\{r\in O , r\neq 0\}$.
The elements $\om_r$ for $r\in O'$ are the minuscule
weights: $(\om_r,\al^\vee)\le 1$ for all $\al \in R_+$.

Given $\tal=[\al,\nu_\al j]\in \tR$ and $b \in P$, let
\begin{align}
&s_{\tal}(\tz)\ =\  \tz-(z,\al^\vee)\tal,\
\ b'(\tz)\ =\ [z,\ze-(z,b)]
\label{ondon}
\end{align}
for $\tz=[z,\ze] \in \R^{n+1}$.

The
{\em affine Weyl group} $\tW$ is generated by all $s_{\tal}$
for $\tal\in\tR_+$; 
we write $\tW = \lan s_{\tal} (\tal\in \tR_+\ran)$. One can take
the simple reflections $s_i=s_{\al_i}\ (0 \le i \le n)$
as its
generators and introduce the corresponding notion of the
length. 
The group $\tW$ is the semidirect product $W\lsmash Q'$ of
its subgroups $W=\lan s_\al, \al \in R_+\ran$ and $Q'=\{a', a\in Q\}$,
where
\begin{align}
& \al'=\ s_{\al}s_{[\al,\,\nu_{\al}]}=\
s_{[-\al,\,\nu_\al]}s_{\al}\for
\al\in R.
\label{ondtwo}
\end{align}

The {\em extended affine Weyl group} $\hW$ generated by $W$ and $P'$
(instead of $Q'$) is isomorphic to $W\lsmash P'$:
\begin{align}
&(wb')([z,\ze])\ =\ [w(z),\ze-(z,b)] \for w\in W,\ b\in B.
\label{ondthr}
\end{align}
From now on,  $b$ and $b',$ $P$ and $P'$ will be identified.

Given $b\in P_+$, let $w^b_0$ be the longest element
in the subgroup $W^{b}\subset W$ of the elements
preserving $b$. This subgroup is generated by simple
reflections. We set
\begin{align}
&u_{b} = w_0w^b_0  \in  W,\ \pi_{b} =
b( u_{b})^{-1} \in \hW, \  u_i= u_{\om_i},\ \pi_i=\pi_{\om_i},
\label{xwo}
\end{align}
where $w_0$ is the longest element in $W,$
$1\le i\le n.$

The elements $\pi_r\equal\pi_{\om_r}\ (r \in O')$ and
$\pi_0=\id$ leave $\tGa$ invariant
and form a group denoted by $\Pi$,
which is isomorphic to $P/Q$ by the natural
projection $\{\om_r \mapsto \pi_r\}$. As to $u_r \ (r\in O')$,
they preserve the set $\{-\vth\}\cup\{\al_i, i>0\}$.
The relations $\pi_r(\al_0)= \al_r= ( u_r)^{-1}(-\vth)$
distinguish the indices $r \in O'$. Moreover, one has
\begin{align}
& \hW  = \Pi \lsmash \tW, \where
  \pi_rs_i\pi_r^{-1}  =  s_j \iif \pi_r(\al_i)=\al_j,\ i,j\ge 0.
\end{align}
We note that $\pi_r^{-1}=\pi_{r^*}$ and $u_r^{-1}=u_{r^*}$ 
for $r\in O'$,
where $r^*\in O'$ is determined by the action of $-w_0$ on
the nonaffine Dynkin diagram $\Gamma$:
$-w_0(\al_r)=\al_{r^*}$.

We will need
the following
{\em affine action} of $\hW$ on
$z \in \R^n$:
\begin{align}
& (wb)\llb z \rrb \ =\ w(b+z),\ w\in W,\ b\in P,\notag\\
& s_{\tal}\llb z\rrb\ =\ z - ((z,\al^\vee)+j)\al,
\ \tal=[\al,\nu_\al j]\in \tR.
\label{afaction}
\end{align}
For instance, $(b w)\llb 0\rrb=b$ for any $w\in W.$
The relation to the above action is given in terms of
the {\em affine pairing}
$([z,l], z'+d)\equal (z,z')+l:$
\begin{align}
& (\hw([z,l]),\hw\llb z' \rrb+d) \ =\
([z,l], z'+d) \for \hw\in \hW,
\label{dform}
\end{align}
where we treat $d$ formally.
\smallskip

\subsection{\bf The length 
\texorpdfstring{on $\hat{W}$}{}}
Setting
$\hw = \pi_r\tw \in \hW\ (\pi_r\in \Pi, \tw\in \tW),$
the length $l(\hw)$
is by definition the length of the reduced decomposition
$\tw= $ $s_{i_l}\cdots s_{i_2} s_{i_1} $
in terms of the simple reflections
$s_i, 0\le i\le n.$ The number of  $s_{i}$
in this decomposition
such that $\nu_i=\nu$ is denoted by   $ l_\nu(\hw).$
We will also use the notation $l_\diamond(\hw)$
for $\diamond=\lng,\sht$.

The {\em length} can be
also defined as the
cardinality $|\la(\hw)|$
of the {\em $\la$\~set} of $\hw$\,:
\begin{align}\label{lasetdef}
&\la(\hw)\equal\tR_+\cap \hw^{-1}(\tR_-)=\{\tal\in \tR_+,\
\hw(\tal)\in \tR_-\},\
\hw\in \hW.
\end{align}
One has:
\begin{align}
&\la(\hw)=\cup_\nu\la_\nu(\hw),\
\la_\nu(\hw)\
\equal\ \{\tal\in \la(\hw),\nu({\tal})=\nu \}.
\label{xlambda}
\end{align}

The coincidence with the previous definition
is based on the equivalence of the {\em length equality}
\begin{align}\label{ltutwa}
&(a)\ \  l_\nu(\hw\hu)=
 l_\nu(\hw)+ l_\nu(\hu)
\for \hw,\hu\in\hW
\end{align}
and the {\em cocycle relation}
\begin{align}
&  (b)\ \ \la_\nu(\hw\hu) = \la_\nu(\hu) \cup
\hu^{-1}(\la_\nu(\hw)),
\label{ltutw}
\end{align}
which, in turn, is equivalent to
the {\em positivity condition}
\begin{align}\label{ltutwc}
& (c)\ \  \hu^{-1}(\la_\nu(\hw))
\subset \tR_+
\end{align}
and is also equivalent to the {\em embedding condition}
\begin{align}\label{ltutwd}
& (d)\ \  \la_\nu(\hu)\subset \la_\nu(\hw).
\end{align}

See, e.g., \cite{C4,C101} and also \cite{Bo,Hu}.
Applying (\ref{ltutw}) to the reduced decomposition
$\hw=\pi_rs_{i_l}\cdots s_{i_2}s_{i_1},$ we obtain an
ordering of the $\la$\~set:
\begin{align}
\la(\hw) = \{\ \tal^1=\al_{i_1},\ &\tal^2=s_{i_1}(\al_{i_2}),\
\tal^3=s_{i_1}s_{i_2}(\al_{i_3}),\notag\\ 
\ldots,\ &\tal^l=\tw^{-1}s_{i_l}(\al_{i_l})\ \}.
\label{tal}
\end{align}
We will call (\ref{tal}) the {\em $\la$\~sequence associated
with the given decomposition of $\hw$}. Such sequences are exactly
those in $\tR_+$ satisfying properties $(i,ii)$ from
the following lemma. We consider $\la(\hw)$ as sets (not sequences)
in quite a few statements below, which then depend only on elements
$\hw$ (not on their reduced decompositions).
\begin{lemma}
\label{la-lem}
Given a reduced decomposition
$\hat{w}=\pi_r s_{j_l}\cdots s_{j_1}\in\hW$, form the
$\la$\~sequence $\la(\hat{w})=\{\tal^1,\ldots,\tal^l\}$
using (\ref{tal}).

(i) If $\tal=\tal^q+\tal^r\in\tilde{R}_+$,
then $\tal=\tal^p$ for some $p$ between $q$ and $r$.
The same holds if \,$\tal=c_1\tal^q+c_2\tal^r\in \tR_+$ 
for positive rational $c_1,c_2$.

(ii) If \,$\la(\hw)\ni
\tal=\tbe+\tga$ for $\tbe,\tga\in\tilde{R}_+\cup [0,\Z_+]$,
then at least one of $\tbe,\tga$ belongs to $\la(\hat{w})$
and exactly one of $\tbe,\tga$ comes before 
$\tal$ in $\la(\hat{w})$.
\end{lemma}

See Main Theorem 2.1 of \cite{C103}. Note that the corresponding
reduced decomposition of $\hw$ can be uniquely recovered from the
$\la$\~sequence $\la(\hw)$; considered as a $\la$\~set, the 
latter is sufficient to recover $\hw$.
\medskip

{\em Reduction modulo $W$.}
The following proposition generalizes the construction of
the elements $\pi_{b}$ for $b\in P_+;$ see \cite{C4} or \cite{C101}.

\begin{proposition} \label{PIOM}
Given $ b\in P$, there exists a unique decomposition
$b= \pi_b  u_b,$
$u_b \in W$ satisfying one of the following equivalent conditions:

{(i) \  } $l(\pi_b)+l( u_b)\ =\ l(b)$ and
$l( u_b)$ is the greatest possible,

{(ii)\  }
$ \la(\pi_b)\cap R\ =\ \emptyset$.

The latter condition implies that
$l(\pi_b)+l(w)\ =\ l(\pi_b w)$
for any $w\in W.$ Besides, the relation $ u_b(b)
\equal b_-\in P_-=-P_+$
holds, which, in turn,
determines $ u_b$ uniquely if one of the following equivalent
conditions is imposed:

{(iii) }
$l( u_b)$ is the smallest possible,

{(iv)\ }
if\, $\al\in \la( u_b)$  then $(\al,b)\neq 0$.
\end{proposition}
\qed
\smallskip

Condition ($ii$) readily
gives a complete description of the set
$\pi_P=\{\pi_b, b\in P\}$, namely,
only
$\,[\,\al<0,\,\nu_\al j>0\,]\,$
can appear in $\la(\pi_b)$.
\smallskip

Explicitly,
\begin{align}
\la(b) = \{ \tal>0,\  &( b, \al^\vee )>j\ge 0 \iif \al\in R_+,
\label{xlambi}\\
&( b, \al^\vee )\ge j> 0 \iif \al\in R_-\},
\notag \\
\la(\pi_b) = \{ \tal>0,\ \al\in R_-,\
&( b_-, \al^\vee )>j> 0
\iif  u_b^{-1}(\al)\in R_+,
\label{lambpi} \\
&( b_-, \al^\vee )\ge j > 0 \iif
u_b^{-1}(\al)\in R_- \}, \notag
\end{align}
For instance,
$l(b)=l(b_-)=-2(\rho^\vee,b_-)$ for $2\rho^\vee=
\sum_{\al>0}\al^\vee.$ 

Switching here to $l_\diamond$ for
$\diamond=\lng,\sht$, one has
$l_\diamond(b)=-2(\rho^\vee_\diamond,b_-)$, where 
$$
\rho^\vee_\diamond\,=\,
\frac{1}{2}\sum_{\al>0,\nu_\al\,=\,
\nu_\diamond}\al^\vee=\rho_\diamond/\nu_\diamond.
$$
\medskip

The element $b_{-}= u_b(b)$ is the unique element
from $P_{-}$ that belongs to the orbit $W(b)$.
Thus the equality $c_-=b_- $ means that $b,c$
belong to the same orbit. We will also use
$b_{+} \equal w_0(b_-),$ the unique element in $W(b)\cap P_{+}.$
In terms of $\pi_b,$
$$u_b\pi_b = b_-,\ \ \pi_b u_b = b_+.$$

Note that $l_\diamond(\pi_b w)=l_\diamond(\pi_b)+l_\diamond(w)$ 
for all $b\in P,\ w\in W.$
For instance,
\begin{align}
&l_\diamond(b_- w)=l_\diamond(b_-)+l_\diamond(w),\ 
l_\diamond(wb_+)=l_\diamond(b_+)+l_\diamond(w),
\label{lupiw}
\\
&l_\diamond(u_b\pi_b w)=l_\diamond(u_b)+l_\diamond(\pi_b)+
l_\diamond(w) \for b\in P,\,
 w\in W.\notag
\end{align}
\smallskip

\comment{
For further references, let us mention the relation 
\begin{align}\label{wlbplus}
&\{[\al,\nu_\al j]\in \la_\diamond(b_+)\,\mid\,
1\le j<(\al,b_+)\}\\
=\,
w\bigl(&\{[\al,\nu_\al j]\,\in\, \la_\diamond(b)\ \mid\ 
1\le j<(\al,b)\}\,\bigr) \hbox{\ if\ } w(b)=b_+,\notag
\end{align}
which readily results from (\ref{xlambi}).
}

{\em Partial orderings on $P$.}
The following two partial orderings on $P$ are commonly
used in the theory of Dunkl operators and
nonsymmetric Macdonald polynomials.
See \cite{C2,O2,M4}. 

We mainly need the partial ordering on $P$ defined by:
\begin{align}
b \preceq c,\, c\succeq b \iif b_-< c_- \hbox{\ or\  }
\{b_-=c_- \hbox{\ and\ } b\le c\}&,\label{succ}\\
\where b \le c,\, c\ge b \for b, c\in P \iif c-b \in Q_+&.\notag
\end{align}
Recall that $b_-=c_- $ means that $b,c$
belong to the same $W$\~orbit.
We write  $<,>,\prec, \succ$ respectively if $b \neq c$.
This ordering was also used in \cite{C2} in the
process of calculating the coefficients of the $Y$\~operators.

For any $b\in P$, we define the sets
\begin{align}
&\si(b)\equal \{c\in P, c\succeq b\},\
\si_*(b)\equal \{c\in P, c\succ b\},\label{cones}\\
&\si_-(b)\equal \si(b_-),\
\si_+(b)\equal \si_*(b_+)= \{c\in P, c_->b_-\}.\notag
\end{align}

\comment{
The following sets
\begin{align}
&\si(b)\equal \{c\in P, c\succeq b\},\
\si_*(b)\equal \{c\in P, c\succ b\}, \notag\\
&\si_-(b)\equal \si(b_-),\
\si_+(b)\equal \si_*(b_+)= \{c\in P, c_->b_-\}.
\label{cones}
\end{align}
are convex.
By {\em convex}, we mean that if
$ c, d= c+r\al\in \si$
for $\al\in R_+, r\in \Z_+$, then
\begin{align}
&\{c,\ c+\al,...,c+(r-1)\al,\ d\}\subset \si.
\label{convex}
\end{align}
}
The second partial ordering is defined by
\begin{align}
&b \preceqq c,\, c\succeqq b \iif b_-< c_- \hbox{\ or\  }
\{b_-=c_- \hbox{\ and\ } u_b \leq u_c\},
\label{succpr}
\end{align}
where $u_b$ is from Proposition \ref{PIOM}
and $\leq$ applied to elements of $W$ is the Bruhat ordering.

It is not hard to show that $c\succeqq b\Rightarrow c\succeq b$
and that the converse is false (see \cite[(2.7.7)]{M1}).
We remark that, if $w_b$ is the unique shortest element of $W$
satisfying $w_b(b_+)=b$, then $u_b\leq u_c$ if and only if
$w_c\leq w_b$.
\medskip

\subsection{\bf On  
\texorpdfstring{$\la$}{Lambda}-sequences of reflections}
\label{SEC:PROPW}
The construction of the Toda-Dunkl operators
will heavily use the sequences $\la(s_{\tal})$
for reflections $s_{\tal}$, where
$\tal=[\al,\nu_\al j]\in \tR_+$. 
As above, $\la$\~sequences will be 
frequently considered as sets, i.e. without the orderings 
determined by reduced decompositions; these sequences are
described intrinsically by Lemma~\ref{la-lem}. 
Formula (1.19) from \cite{C103} states 
that all $\la$\~sequences for $s_{\tal}$ are as follows:
\begin{align}\label{reflambda}
&\la(s_{\tal})=\{s_{\tal}(-\la(\tw))\}_{op}\,\cup\,\tal\,\cup\,
\la(\tw)\ \, \hbox{(as sequences)},
\end{align}
where $\tw\in\tW$ is of minimal possible length such that
$\tw(\tal)=\al_m$ among all $m\ge 0$;
by $\{\ \cdot\ \}_{op}\ $, we mean the inversion of
the ordering of a given sequence  $\{\ \cdot\ \}$. Then
$s_{\tal}=\tw^{-1}s_m\tw$ is reduced for any reduced decomposition
of $\tw$ and an arbitrary reduced decomposition of $s_{\tal}$
can be presented in this
form for proper $\tw,s_m$ satisfying the above minimality condition:
\begin{align}
\label{sal-red}
&s_{\tal}\ =\ s_{j_1}\cdots s_{j_p}s_m s_{j_p}\cdots s_{j_1}, 
\where\\
&\tw=s_{j_p}\cdots s_{j_1},\, p=l(\tw),\ j_1,\ldots,j_p\ge 0.
\end{align}
See e.g. Proposition 1.1 from \cite{C103}
for a proof of these (standard) facts. 
We observe that if such a reduced decomposition is used
to construct the $\la$\~sequence 
$\la(s_{\tal})=\{\tbe^1,\ldots,\tbe^l\}$, where
$l=l(s_{\tal})=2p+1$, then
\begin{align}
\tbe^{l-i+1}=-s_{\tal}(\tbe^i), \for 1\leq i\leq l.
\label{tbe-opp}
\end{align}

Furthermore, for any $\hw\in \hW$ and any sequence $\la(\hw)$, 
one has
\begin{align} \label{talinla}
\tal\in \la(\hw)\Leftrightarrow 
\la(s_{\tal})\setminus\{\tal\} =  
\bigcup_{\tbe}\,\{\tbe,\,\tbe'=-s_{\tal}(\tbe)\}
\hbox{\,(as sets),\ \, where}&\notag\\ 
\tbe\,\in\,\la(\hw)\cap \la(s_{\tal}) \hbox{\ \,such that\ \,} \tbe 
\hbox{\ \,appears in\ \,} \la(\hw) \hbox{\,\ before\ \,} \tal.&
\end{align}
See formula (1.20) in \cite{C103}.
Here $\tbe$ and $\tbe'=-s_{\tal}(\tbe)$ do not coincide 
unless $\tbe=\tal$. Indeed, $\tbe'
=2\frac{(\al,\be)}{(\al,\al)}\tal-\tbe$ and
$\tbe'=\tbe$ if and only if $\tbe$ is proportional to
$\tal$, which occurs exactly for $\tbe=\tal$.
All pairs $\{\tbe,\tbe'\}$ are pairwise distinct. Indeed, if 
$\tga=\tbe'$ for any $\tbe,\tga\in \la(s_{\tal})$,  
then $\tal$, which is $c(\tbe+\tbe')$ for $c>0$,
occurs between $\tbe$ and $\tga$ in this sequence, 
which is impossible by construction. Here and below,
see Lemma~\ref{la-lem}.
\smallskip

Let us list some other properties of the sets $\la(s_{\tal})$.
First of all, always $(\al,\be)>0$ for any 
$\tbe=[\be,\nu_\be k]\in \la(s_{\tal})$ and we have the
following equivalent inequalities :
\begin{align}\label{stalineq}
2\frac{(\al,\be)}{(\al,\al)}\nu_\al j\ge \nu_\be k \,
\Longleftrightarrow\, (\al,\be)j\ge \nu_\be k \,
\Longleftrightarrow\, 
2\frac{(\al,\be)}{(\be,\be)}j\ge k.
\end{align}

The strict inequality $(\al,\be)j>\nu_\be k$
for $\tbe>0$ gives that $\tbe\in\la(s_{\tal})$. If 
$(\al,\be)j=\nu_\be k$, then the conditions
$(\al,\be)>0$ and $s_\al(\be)<0$ are necessary and sufficient.

If $\be\neq \al$, then the inequality in (\ref{stalineq})
becomes $j\ge k$
unless $\be$ is short and $\al$ is long; in the latter
case, it becomes $\nu_\al j\ge k$. 
If $\be=\al$, then
$\tbe=[\al,\nu_\al k]$ are in $\la(s_{\tal})$ if and only
if $0\le k\le 2j$ when $\al>0$, and $0<k<2j$ when $\al<0$.
\smallskip

Following  the calculation of
the sets $\la(s_{\tal})$ for $\al\in \tR_-,j>0$ 
in formula (1.28) from \cite{C103} and the
action of $-s_{\tal}$ in these sets described there, 
let us calculate such sets for arbitrary positive 
affine roots. 

\begin{lemma}\label{LASTAL}
For $\tal=[\al,\nu_\al j]\in \tR_+$ and \,$\diamond=\sht,\lng$,
we set 
\begin{align*}
&\de_{\al,\diamond}\!=\!\de_{\nu_{\al},\nu_\diamond},\ \,
\eta_{\al\diamond}\!=\!1 \hbox{\ \,unless \ }
\eta_{\al\diamond}\!=\!\nu_{\lng} \hbox{ when } \nu_\al\!=\!1\!
=\!\nu_{\sht}
\hbox{ and } \diamond\!=\!\lng.
\end{align*}

(i) For $\al\in R_+$, there exists a set 
\,$\{\be^i\}\subset R_+
\setminus\{\al\}$\,
such that 
\begin{align}\label{lashtlng1}
&\la_\diamond(s_{\al})\setminus\{\al\}=
\bigl\{\,\{\be^i,-s_\al(\be^i)\}\ \mid \
\nu_{\be^i}=\nu_\diamond,\ s_\al(\be^i)<0\,\bigr\},\\
&\hbox{where\ \,} 1\le i\le \frac{l_\diamond(s_{\al})-
\de_{\al,\diamond}}{2},\ \,
l_{\diamond}(s_{\al})\, =\, 
2\frac{(\al,\rho_\diamond)}{\eta_{\al\diamond}\,\nu_\al}-
\de_{\al,\diamond}.\notag
\end{align}
More explicitly, $l_\diamond(s_{\al})=2(\al^\vee,\rho_\diamond)-
\de_{\al,\diamond}$
for long roots $\al\in R_+$ and 
$l_\diamond(s_{\al})=2(\al,\rho^\vee_\diamond)-
\de_{\al,\diamond}$ for short $\al$.
\smallskip

(ii) Let $\tal=[\al,\nu_\al j]$ for $\al\in R_+,j>0$. 
Then provided that $\tbe=[\be,\nu_\be k]>0$ such that\,
$\be\neq \al$,\,  $\nu_\be=\nu_\diamond$\, and\, $(\al,\be)>0$, 
\begin{align}
\la_\diamond(s_{\al})= 
\{\,[\al,\nu_\al k]\,&,\, 0\le k\le 2j\,\}\,\cup\,
\{\,\tbe>0\,,\, 0\le \nu_\be k<(\al,\be)j\,\}
\notag\\
\cup\,&
\{\,[\be,(\al,\be)j]\,;\, \be<0 \hbox{\ or\ }
\be\in \la(s_{\al})\,\}.\label{lashtlnga}
\end{align}
Accordingly, $l_\diamond(s_{\tal})=l_\diamond(j\al)
+l_\diamond(s_\al)$, 
where $j\al=s_\al s_{\tal}=s_{[-\al,\nu_\al j]}s_{\al}$.
\smallskip

(iii) Let $\tal=[-\al,\nu_\al j]$ for $\al\in R_+,j>0$.
Then assuming that $\tbe=[-\be,\nu_\be k]>0$ and that
$\nu_\be=\nu_\diamond$,
\begin{align}\label{lamtbe}
\la_\diamond(s_{\tal})=\{[-\be,\nu_\al k]\in \tR_+\ &,\
\ 0\le \nu_\be k <(\al,\be)j\}\notag\\
\cup\, \{[-\be,\,(\al,\be)j]\ &;\ 
\be>0<s_{\al}(\be),\, (\al,\be)>0 \}.
\end{align}
One has\  $\la_\diamond(s_{\tal})=\la_\diamond (-j\al)\setminus
\{\,[-\be,\, (\al,\be)j]\,\mid\, \be\in \la_\diamond(s_\al)\,\}$.
\end{lemma}
\smallskip
{\em Proof.} The presentation of $\la(s_\al)$
from $(i)$ is a particular case
of (\ref{talinla}). Then we will use that
$\rho_\diamond-w(\rho_\diamond)=
\sum_{\be\in\la_\diamond(w)} \be,
$  
combining it with (\ref{lashtlng1}) and the
formula $\be'=-s_\al(\be)=\eta_{\al\be}\al-\be$ for 
$\be\in \la(s_\al)\setminus\{\al\}$, 
where $\eta_{\al\be}=1$ unless
$\eta_{\al\be}=\nu_{\lng}$ for short $\al$ and long $\be$. 
One has
\begin{align}\label{rhodiamal}
&(\rho_\diamond-s_\al(\rho_\diamond),\al)/\nu_\al\ =\ 
2(\rho_\diamond,\al)/\nu_\al\\
= \de_{\al,\diamond}+
\!\!\sum_{\be\in\la_\diamond(s_\al)}& 
(\frac{\eta_{\al\be}\al}{2}\,,\al)/\nu_\al= 
\de_{\al,\diamond}+
\!\!\sum_{\be\in\la_\diamond(s_\al)} \eta_{\al\be}= 
\de_{\al,\diamond}+
\eta_{\al\diamond}\,l_{\diamond}(s_\al).\notag
\end{align}

Claim $(ii)$ follows from the previous considerations
and formula (\ref{xlambi}). One can also use that
$\la(s_\al)\cap \la(-j\al)=\emptyset$. Claim $(iii)$
is formula (1.28) from \cite{C103}. It follows
from $(i)$ or $(ii)$ using  the decomposition
$s_{[-\al,\nu_\al j]}= s_\al\cdot(-j\al)$; 
here $-j\al= s_{[-\al,\nu_\al j]}s_\al$ is reduced.
\sq
\smallskip

{\em Nonaffine reflections.}
Let $s_\al=w^{-1}s_iw=s_{j_1}\cdots s_{j_p}s_m s_{j_p}\cdots s_{j_1}$ 
be a reduced decomposition from (\ref{sal-red})
for $\al\in R_+$ and proper $m>0$; here 
$w=s_{j_p}\cdots s_{j_1}$, $l(w)=p$.
Then one has the inequalities 
\begin{align}\label{negbemin}
&(\al,\be)>0 \for \be \in \la(w), \hbox{\ \, equivalently,\,\ }\\
&(\al_m,\be')<0 \for \be'=-w(\be)\in \la(w^{-1})=
-w(\la(w)).\notag
\end{align}

\begin{lemma}
\label{lem-sht-pos}
For any $\al\in R_+$, let $w=s_{j_k}\cdots s_{j_1}$ be a
reduced decomposition of an element $w\in W$ such 
that $w(\al)=\al_m$ and $k=l(w)$ is
minimal possible among all $\al_m$ (which is then
$k=(l(s_\al)-1)/2$). Equivalently, 
$s_\al=s_{j_1}\cdots s_{j_k}s_m s_{j_k}\cdots s_{j_1}$
is reduced.

(i) For such $w$, $\al_m$ and any reduced decomposition of $w$,
there exists its extension $\breve{w}=s_{j_p}\cdots s_{j_1}$ 
of length $p=(l(s_{\th'})-1)/2$ and the corresponding 
reduced decomposition
$s_{\th'}=s_{j_1}\cdots s_{j_p}s_m s_{j_p}\cdots s_{j_1}$, where 
$\th'$ is the maximal root $\th$ for long $\al$ and $\vth$ for
short $\al$;\, $m,j_1,\ldots,j_p>0$.

(ii) When $\al=\th'$, any $m=1,\ldots,n$ can be taken here 
provided that $|\al_m|=|\th'|$;
the corresponding element $w=w^{(m)}=s_{j_p}\cdots s_{j_1}$
(but not its reduced decomposition) is
uniquely determined by the choice of $m$, equivalently,
by the condition $w(\th')=\al_m$ together with the 
inequalities $(\be,\th')>0$ for all $\be\in \la(w)$.
\end{lemma}
\proof
Let us demonstrate that the inequalities 
$(\al_m,\be')<0$ from (\ref{negbemin}) 
for all $\be'\in \la(w^{-1})$
are actually sufficient to ensure that $l(w)$ is minimal
possible among all $w$ such that $w(\al)=\al_m$
for a given $\al_m$. 

We argue by ``descending" induction on $l(s_\al)$. 
Unless $\al$ satisfying
inequalities (\ref{negbemin}) is
maximal long or short root, we can find a simple root
$\al_j$ such that $(\al,\al_j)<0$. Then $\al_j\not\in \la(s_{\al})$.
Therefore $\al_j\not\in \la(w)$, $w s_j$ is reduced,
$(\al_m,w(\al_j))=(w^{-1}(\al_m),\al_j)<0$  and, finally,
$\breve{w}\equal ws_j$ is of length $l(w)+1$ satisfying the
inequalities from (\ref{negbemin}) for any 
$\breve{\be}'\in \la(\breve{w}^{-1})=w(\al_j)\cup \la(w^{-1})$. 
Continuing this way by induction, we eventually construct $w$  
satisfying $w(\th')=\al_m$ and 
the inequalities from (\ref{negbemin}), where $|\th'|=|\al_m|$.

Let $\th'=\vth$ here for the sake of definiteness. 
We claim that the resulting $w$ is a minimal element 
satisfying $w(\vth)=\al_m$; moreover, it is unique (depends
only on $\al_m$).  Indeed, the subgroup
$W^{\vth}=\{u\in W\,\mid\,u(\vth)=\vth\}$ is parabolic 
generated by simple $s_r$ such that 
$r \neq k$ for $\al_k$
connected with $\al_0$ in the affine Dynkin diagram
$\tilde{\Ga}$ for the (twisted) root system $\tR$.
We use here that $\vth=\om_k$ unless for $A_n$,
where $\vth=\om_1+\om_n$. Due to the inequalities from
(\ref{negbemin}), all products $wu$ are reduced for such $w$
and any elements $u\in W^{\vth}$, so $w$ is really 
minimal and unique such.

Here one can begin with any short simple $\al=\al_m$, 
which proves $(ii)$. Moreover, 
the induction process above automatically guarantees
that $\breve{w}=w s_j$ has to be minimal (though maybe 
not unique such) for $\breve{\al}=s_j(\al)=\breve{w}^{-1}(\al_i)$, 
as well as for all consecutive $w$ serving
$\{\al,\breve{\al},\ldots,\vth\}$, since
the last $w$ in this chain has been proven to be minimal.
This justifies that the inequalities in (\ref{negbemin}) 
are sufficient for the minimality of $w$ and gives $(ii)$.
\sq
\smallskip

As a by-product, we obtain that for any reduced decomposition
$w=s_{j_p}\cdots s_{j_1}$ of minimal $w$ such that
$w(\vth)=\al_m$, where $\al_m$ is any given short root,
\begin{align}
\label{inn-prod-min1}
(s_{j_{k+1}}\cdots s_{j_p}(\al_m),\al_{j_k}^\vee)
=(s_{j_k}\cdots s_{j_1}(\vth),\al_{j_k}^\vee)=-1,
\ \ 1\leq k\leq p.
\end{align}
Indeed, using (\ref{negbemin})
\begin{align*}
&(s_{j_{k+1}}\cdots s_{j_p}(\al_m),\al_{j_k})=
(\al_m,s_{j_p}\cdots s_{j_{k+1}}(\al_{j_k})), \\
&\hbox{where \ \, } s_{j_p}\cdots s_{j_{k+1}}(\al_{j_k})\,\in\, 
\la(w^{-1}),
\end{align*}
which gives that the right-hand side in (\ref{inn-prod-min1})
is negative; so it must be $-1$ since $\al_m$ is short.

It is of interest to calculate explicitly the elements 
$w^{(m)}$ from $(ii)$
and their $\la$\~sets. Let us do the latter for $m=k$ 
for short $\al_k$ 
connected with $\al_0$ in $\tilde{\Ga}$.  Then $w'(\vth)=\al_k$ for
$w'=s_{\vth}s_k$ and $\la(w')=s_k(\la(s_{\vth})\setminus \al_k)$.
Dividing $w'$ by the maximal possible $u\in W{\vth}$ on the right, 
we obtain that
\begin{align*}
\la(w^{(k)})=\{\,\be\in R_+\,|\, (\be,\vth)>0, (\be,\al_k)=0\,\}.
\end{align*}
Indeed, we need to remove $s_k(\be)$ from $\la(w')$ for
$\be\in \la(s_{\vth})\setminus \al_k$ such that 
$(s_k(\be),\vth)=0=(\be,\vth-\al_k)$. Recall that 
$\be+\be'=\nu_\be \vth$ for $\be'=-s_{\vth}(\be)$.
Therefore, either $(\be,\al_k)=\nu_\be$ and 
$(\be',\al_k)=0$ or  $(\be,\al_k)=0$ and 
$(\be',\al_k)=\nu_\be$. Thus $(\be,\vth-\al_k)\neq 0$
is equivalent to $(\be,\al_k)=0$. This also 
gives that the number of such 
$\be$ in $\la(s_{\vth})$ is $(l(s_{\vth})-1)/2$ because exactly 
one root from each pair $\{\be,\be'\}$ is orthogonal to $\al_k$.
 
\comment{
\begin{lemma}
\label{lem-sht-pos}
For any short $\al\in R_+\setminus\{\vth\}$,
there is a reduced decomposition
$s_\vth=s_{j_1}\cdots s_{j_p}s_m s_{j_p}\cdots s_{j_1}$ (where 
$m,j_1,\ldots,j_p>0$) such that 
$\al=s_{j_r}\cdots s_{j_1}(\vth)$ for some $1\leq r\leq p$.
\end{lemma}
\proof
Write $\vth=s_{j_1}\cdots s_{j_r}(\al)$ and
$\al=s_{j_{r+1}}\cdots s_{j_p}(\al_m)$ so that
\begin{align}
(s_{j_{i+1}}\cdots s_{j_p}(\al_m),\al_{j_i}^\vee)<0, \ 
\ 1\leq i\leq p.
\end{align}
Then $\vth=s_{j_1}\cdots s_{j_p}(\al_m)$ and one has 
$p\leq(\vth,\rho^\vee)-1$. The decomposition
$s_\vth=s_{j_1}\cdots s_{j_p}s_m s_{j_p}\cdots s_{j_1}$
must be reduced, because
\begin{align}
l(s_\vth)=l(\vth)-l(s_0)=2(\vth,\rho^\vee)-1,
\end{align}
and hence $p=(\vth,\rho^\vee)-1$. We note that, consequently,
one must have
\begin{align}
\label{inn-prod-min1}
(s_{j_{i+1}}\cdots s_{j_p}(\al_m),\al_{j_i}^\vee)
=(s_{j_i}\cdots s_{j_1}(\vth),\al_{j_i}^\vee)=-1,
\ \ 1\leq i\leq p.
\end{align}
\sq
}
\smallskip

\subsection{\bf Main definition}
Let $m$ denote the least natural number
such that $(P,P)=(1/m)\Z.$  Thus
$m=2$ for $D_{2k}$, $m=1$ for $B_{2k}$ and $C_{k}$,
and $m=|\Pi|$ otherwise.

The double affine Hecke algebra depends
on the parameters
$q, t_\nu \ (\nu\in\nu_R).$ 
It will be defined
over the ring
$\Q_{q,t}\equal \Q[q^{\pm 1/(2m)},t_\nu^{\pm 1/2}]$.
Later we will need the field of fractions
$\Q_{q,t}'\equal\Q(q^{\pm 1/(2m)},t_\nu^{1/2})$
and its subrings
\begin{align}
\label{dQprime}
&\ddot{\Q}_{q,t}'\equal\{c\in \Q_{q,t}' \, ;\,
c \hbox{\,\, is well defined when all\ \,} t_\nu^{1/2}=0\},\\
\label{dQdag}
&\ddot{\Q}_{q,t}^\dag\equal\{c\in \Q_{q,t}' \, ;\,
c \hbox{\,\, is well defined when all\ \,} t_\nu^{-1/2}=0\}.
\end{align}
We set
\begin{align}
&t_{\tal}=t_{\al}=t_{\nu_\al},\ \ t_i = t_{\al_i},\ \
q_{\tal}=q^{\nu_\al},\ \ q_i=q^{\nu_{\al_i}},\notag\\
&\where \tal=[\al,\nu_\al j] \in \tR,\ \ 0\le i\le n.
\label{taljx}
\end{align}

It will be convenient to use parameters
$\{k_\nu\}$ together with $\{t_\nu \},$ setting
here and further:
$$
t_\al=t_\nu=q_\al^{k_\nu} \for \nu=\nu_\al, \and
\rho_k=(1/2)\sum_{\al>0} k_\al \al.
$$
Note that $(\rho_k,\al_i^\vee)=k_i=k_{\al_i}=
(\rho_k^\vee,\al_i)$ for $i>0$, where
$$
\rho_k^\vee\ \equal\ \sum k_{\nu}\rho_\nu^\vee\, \for\,
\rho_\nu^\vee\ \equal\ \rho_\nu/\nu.
$$

Using that $w_0(\rho_k)=-\rho_k$, we obtain that
 $(\rho_k,-w_0(b))=(\rho_k,b)$. For instance,
$(\rho_k,b_+)=-(\rho_k,b_-)$, where
$b_{+} \equal w_0(b_-)$ (see above).

By $q^{(\rho_k,\al)}$, we mean
$\prod_{\nu\in\nu_R}t_\nu^{(\rho_\nu^\vee,\al)}$;
here $\al\in R$, $(\rho_\nu^\vee,\al)\in \Z$ and
this product contains
only {\em integral} powers of $t_{\sht}$ and $t_{\lng}$
(non-negative if $\al>0$).

For pairwise commutative $X_1,\ldots,X_n,$ let
\begin{align}
& X_{\tb}\ =\ \prod_{i=1}^nX_i^{l_i} q^{ j}
\iif \tb=[b,j],\ \hw(X_{\tb})\ =\ X_{\hw(\tb)}.
\label{Xdex}\\
&\hbox{where\ } b=\sum_{i=1}^n l_i \om_i\in P,\ j \in
\frac{1}{ m}\Z,\ \hw\in \hW.
\notag \end{align}
For instance, $X_0\equal X_{\al_0}=qX_\vth^{-1}$.

We set $(\tilde{b},\tilde{c})=(b,c)$,
ignoring the affine extensions
in this pairing.

Recall that $m_{ij}$ denotes the order of $s_is_j$ in $\tW$
$(0\le i,j\le n)$ and that $r,r^*\in O'$ are related by
$-w_0(\om_r)=\om_{r^*}$.

\begin{definition}
The double affine Hecke algebra $\HH\ $
is generated over
$\Q_{q,t}$ by the elements $\{T_i , 0\le i\le n\}$,
pairwise commutative $\{X_b , b\in P\}$ satisfying
(\ref{Xdex}),
and the group $\Pi$,
where the following relations are imposed:

(o)\ \  $ (T_i-t_i^{1/2})(T_i+t_i^{-1/2})\ =\
0,\ 0\ \le\ i\ \le\ n$;

(i)\ \ \ $ T_iT_jT_i\cdots = T_jT_iT_j\cdots, m_{ij}$
factors on each side, $0\le i\ne j\le n$;

(ii)\ \   $ \pi_rT_i\pi_r^{-1}\ =\ T_j \iif
\pi_r(\al_i)=\al_j$;

(iii)\  $T_iX_b \ =\ X_b X_{\al_i}^{-1} T_i^{-1} \iif
(b,\al^\vee_i)=1,\
0 \le i\le  n$;

(iv)\ $T_iX_b\ =\ X_b T_i \iif (b,\al^\vee_i)=0,\
 0 \le i\le  n$;

(v)\ \ $\pi_rX_b \pi_r^{-1}\ =\ X_{\pi_r(b)}\ =\
X_{ u^{-1}_r(b)}
 q^{(\om_{r^*},b)},\  r\in O'$.
\label{double}
\end{definition}

%
%

One can rewrite ($iii, iv$) as in \cite{L}:
\begin{align}
&T_iX_b -X_{s_i(b)}T_i\ =\
(t_i^{1/2}-t_i^{-1/2})\frac{X_{s_i(b)}-X_b}
{X_{\al_i}-1},\ 0 \le i\le  n.
\label{tixi}
\end{align}

Given $\tw \in \tW, r\in O,\ $ the product
\begin{align}
&T_{\pi_r\tw}\equal \pi_r\prod_{k=1}^l T_{i_k},\where
\tw=\prod_{k=1}^l s_{i_k},
l=l(\tw),
\label{Twx}
\end{align}
does not depend on the choice of the reduced decomposition
(because $T_i$ satisfy the same ``braid'' relations
as $s_i$ do).
Moreover,
\begin{align}
&T_{\hv}T_{\hw}\ =\ T_{\hv\hw}\  \hbox{ whenever}\
 l(\hv\hw)=l(\hv)+l(\hw) \for
\hv,\hw \in \hW. \label{TTx}
\end{align}
In particular, we arrive at the pairwise
commutative elements:
\begin{align}
& Y_{b} = \prod_{i=1}^nY_i^{l_i} \iif
b=\sum_{i=1}^n l_i\om_i\in P,\
Y_i\equal T_{\om_i},\ b\in P\label{Ybx}.
\end{align}

For any $b\in P$, the element $Y_b$ can be presented
as the product $\pi_r T_{j_l}^{\pm 1}\cdots T_{j_1}^{\pm 1}$ for
any reduced decomposition $b=\pi_r s_{j_l}\cdots s_{j_1}$
and a proper choice of signs $\pm$ (see (\ref{def-ep-p}) below).
Note that $l=l(b)=2(\rho^\vee,\,b_+)$ depends only on $b_+$.
The total number of factors $T_j^{\pm 1}$ in this product
with $\nu_j=\nu$ equals $2(\rho_\nu^\vee,\,b_+)$. 

The signs $\pm$ can be described as follows
(see \cite[(3.2.10)]{M1}).
Given the reduced decomposition above, form $\la(b)$ using (\ref{tal})
and write $\tal^p=[\al^p,\nu_{\al^p} j]$.
Then one has $Y_b=\pi_r T_{j_l}^{\ep_l}\cdots T_{j_1}^{\ep_1}$, where
\begin{align}
\label{def-ep-p}
&\ep_p=
\begin{cases}
+1&\text{if $\al^p>0$,}\\
-1&\text{if $\al^p<0$.}
\end{cases}
\end{align}

{\em Duality anti-involution.}
There exists a unique anti-involution $\vph$ of $\HH$ satisfying
(see \cite{C15}):
\begin{align}
\vph: X_b \leftrightarrow Y_{-b},\
T_i\mapsto T_i \,(1\leq i\leq n),\
q^{\frac{1}{2m}}\mapsto q^{\frac{1}{2m}},
\ t_\nu^{1/2}\mapsto t_\nu^{1/2}.
\label{phianti}
\end{align}
Using $Y_\vartheta=T_0T_{s_\vartheta}$ and $Y_{\om_r}=\pi_r T_{u_r}$,
one finds that
\begin{align}
\label{vp0r}
\vph(T_0)=T_{s_\vth}^{-1}X_\vth^{-1},\ \ 
\vph(\pi_r)=T_{u_r^{-1}}^{-1}X_{\om_r}^{-1}=
X_{\om_{r^\ast}}T_{u_r}=\vph(\pi_{r^\ast}^{-1}).
\end{align}

Applying $\vph$ to ($iii,iv$) in the definition of $\HH$, 
we obtain the dual relations (for $i>0$ only):
\begin{align}\label{TYTL}
&T^{-1}_i Y_b\ =\ Y_{s_i(b)} T_i\iif
(b,\al^\vee_i)=1,\\
& T_iY_b\ =\ Y_b T_i \iif (b,\al^\vee_i)=0,
 \ 1 \le i\le  n.\notag
\end{align}
\smallskip
The counterpart of (\ref{tixi}) reads as follows:
\begin{align}
&T_iY_b -Y_{s_i(b)}T_i\ =\
(t_i^{1/2}-t_i^{-1/2})\frac{Y_b-Y_{s_i(b)}}
{1-Y_{-\al_i}},\ 1 \le i\le  n.\label{tiyi}
\end{align}

{\em Automorphisms.}
We will need the following automorphisms of $\HH$.
We refer to \cite{C101} and the references therein for 
proofs and for a discussion of how these automorphisms
can be described in terms of an action of $PSL(2,\Z)$.

We say that an automorphism (or anti-automorphism) of $\HH$
{\em preserves} $q,t_\nu$ if it fixes all fractional powers of
these parameters (i.e., it is $\Q_{q,t}$\~linear). We say
that an automorphism {\em conjugates} $q,t_\nu$ to mean that
it sends all fractional powers of these parameters to their
inverses (so such a map is only $\Q$\~linear).

The following map can be uniquely extended to
an automorphism of $\HH$ fixing
$T_i\ (i\ge 1)$ and preserving $q,t_\nu$:
\begin{align}
\label{taux}
\tau_+ :\ &X_b \mapsto X_b,\
\pi_r \mapsto q^{-\frac{(\om_r,\om_r)}{2}}X_{\om_r}\pi_r,
\ Y_{\om_r} \mapsto X_{\om_r}Y_{\om_r} 
q^{-\frac{(\om_r,\om_r)}{2}},\\
&T_0\mapsto X_{\al_0}^{-1} T_0^{-1},\
Y_\vth \mapsto  X_{\al_0}^{-1}T_0^{-1}T_{s_\vth}.\notag
\end{align}

Define the automorphism $\tau_-\equal\vph\tau_+\vph$.
Explicitly, $\tau_-$ fixes $T_i \ (i\ge 1)$, as well as $\tau_+$,
preserves $q, t_\nu$, and satisfies
\begin{align}
\label{tauy}
\tau_- : X_{\om_r}\mapsto q^{(\om_r,\om_r)/2}Y_{\om_r}X_{\om_r},\
\pi_r\mapsto \pi_r,\ T_0\mapsto T_0,\ Y_b\mapsto Y_b.
\end{align}

We also need the following automorphism of $\HH$:
\begin{align}
\si\equal\tau_+\tau_-^{-1}\tau_+=\tau_-^{-1}\tau_+\tau_-^{-1};
\end{align}
it preserves $q,t_\nu$ and satisfies
\begin{align}
\label{si}
\si:\ &T_i \mapsto T_i \ (i>0),\
X_b \mapsto Y_b^{-1},\
Y_{\om_r}\mapsto q^{-(\om_r,\om_r)}Y_{\om_r}^{-1}
X_{\om_r}Y_{\om_r},\\
&\pi_r \mapsto X_{\om_r}T_{u_r^{-1}}=T_{u_r}^{-1}
X_{\om_{r^\ast}}^{-1},\
T_0\mapsto T_{s_\vth}^{-1}X_{\vth}^{-1}.\notag
\end{align}
The equality of the two expressions for $\si(\pi_r)$ follows
from (\ref{vp0r}).

The following map can be uniquely extended to an involution
of $\HH$ conjugating $q,t_\nu$:
\begin{align}
\label{ve}
\vep : T_i\mapsto T_i^{-1}\ (i>0),\ X_b \leftrightarrow Y_b.
\end{align}
Using $Y_\vartheta=T_0T_{s_\vartheta}$ and $Y_{\om_r}=\pi_r T_{u_r}$,
one finds that
\begin{align}
\label{ve0r}
\vep(T_0)=X_\vartheta T_{s_\vartheta}, \ \
\vep(\pi_r)=X_{\om_r}T_{u_r^{-1}}.
\end{align}

We will also need the involution $\eta\equal\vep\si=\si^{-1}\vep$.
Explicitly, $\eta$ conjugates $q,t_\nu$ and satisfies
\begin{align}
\label{eta}
\eta: T_i\mapsto T_i^{-1} \ (i\ge 0),\
X_b \mapsto X_b^{-1},\ \pi_r\mapsto \pi_r.
\end{align}

\subsection{\bf Double-dot normalization}
The case when all $t_\nu=0$ will play an important role in this paper.
Definition \ref{double}, which matches that from \cite{C101},
\cite{C4}, and other first author's papers, is not suited to this
specialization. We introduce the following normalization to handle
the specialization $t_\nu=0$.

We set
\begin{align}
&\ddot{T}_i\equal t_i^{1/2}T_i,\ \ \ \
\ddot{T}'_i\equal t_i^{1/2}T_i^{-1}=\ddot{T}_i-(t_i-1).
\end{align}
Note that the same normalization is used for both
$T_i$ and $T_i^{-1}$, so that $\ddot{T}_i \ddot{T}'_i=t_i$.
Thus $(\ddot{T}_i)^{-1}$ and $\ddot T_i^{-1}$ do not coincide.

We observe that $\{\ddot{T}_i\}$ satisfy the same braid relations
as $\{T_i\}$. The quadratic relations read:
$(\ddot{T}_i-t_i)(\ddot{T}_i+1)=0.$

The dot-normalization $\ddot{\{\,\} }$ will be applied
term-wise to the products of
$\pi_r\in \Pi, T_i, T_i^{-1}$
provided that the corresponding word in $\hW$ is reduced.
We set $\ddot{\pi}_r=\pi_r$.

Thus, using the description of $Y_b$ from the previous
section, one has $\ddot{Y}_{b}=q^{(b_+,\rho_k)}Y_b$.
Equivalently, one can set $\ddot{Y}_i\equal \ddot{T}_{\om_i}$ and
define
\begin{align}
& \ddot{Y}_b = q^{(b_+-b,\,\rho_k)}
\prod_{i=1}^n\ddot{Y}_i^{l_i} \iif
b=\sum_{i=1}^n l_i\om_i.
\label{Ybxx}
\end{align}
Note that $\ddot{Y}_b \ddot{Y}_{-b}= q^{2(b_+,\,\rho_k)}$.

The first line in (\ref{TYTL}) can be rewritten as
\begin{align}
\ddot T_i'\, \ddot Y_b\ &=\ \ddot Y_{s_i(b)} \ddot T_i
\iif (b,\al^\vee_i)=1,
\end{align}
and (\ref{tiyi}) becomes
\begin{align}
&\ddot{T}_i\ddot{Y}_b -\ddot{Y}_{s_i(b)}\ddot{T}_i\ =\
(t_i-1)\frac{\ddot{Y}_b-\ddot{Y}_{s_i(b)}}
{1-q^{-(\th_i,\,\rho_k)}\ddot{Y}_{-\al_i}},\ 1 \le i\le  n,
\label{tiyidot}
\end{align}
where $\th_i=\th,\vth$ for long, short $\al_i$ respectively;
to see this, use that $\th_i$ is the only root in the intersection
$W(\al_i)\cap P_+$.

We come to the following definition.
Let us extend the scalars of $\HH$ to $\Q_{q,t}'$ and define
$\dHH\subset\HH$ to be the
subalgebra generated over $\ddot{\Q}_{q,t}'$ by the elements
$$
X_a \ (a\in P),\ \ \ddot{T}_{\hw} \ (\hw\in \hW),\ \ 
\ddot{Y}_b \ (b\in P).
$$
It suffices to take here only
$$
X_a \ (a\in P),\ \ \ddot{T}_{i} \ (i\ge 0),\ \ \Pi.
$$

\begin{definition}\label{DEFHHDOT}
The defining relations of 
$\dHH$ as an abstract algebra are as follows:

($o$)\ \  $ (\ddot{T}_i-t_i)(\ddot{T}_i+1)\ =\
0,\ 0\ \le\ i\ \le\ n$;

($i$)\ \ \ $ \ddot{T}_i\ddot{T}_j\ddot{T}_i\cdots =
\ddot{T}_j\ddot{T}_i\ddot{T}_j\cdots, m_{ij}$
factors on each side, $0\le i\ne j\le n$;

($ii$)\ \   $ \pi_r\ddot{T}_i\pi_r^{-1}\ =\ \ddot{T}_j \iif
\pi_r(\al_i)=\al_j$,\, $\pi_r\in \Pi$;

($iii$)\  $\ddot{T}_iX_b \ =\ X_b X_{\al_i}^{-1} \ddot{T}_i' \iif
(b,\al^\vee_i)=1,\ 0 \le i\le  n$;

($iv$)\ $\ddot{T}_iX_b\ =\ X_b \ddot{T}_i \iif (b,\al^\vee_i)=0
\for 0 \le i\le  n$;

($v$)\ \ $\pi_rX_b \pi_r^{-1}\ =\ X_{\pi_r(b)}\ =\
X_{ u^{-1}_r(b)} q^{(\om_{r^*},b)}$.
\sq
\end{definition}


The limit (reduction) of
$\dHH\,$ when $t_\nu=0$ for all or some $\nu$
will be denoted by $\bHH\,$ later and will be
called the {\em nil-DAHA}.

\setcounter{equation}{0}
\section{\sc Polynomial representation}
\label{SEC:POLREP}
From now on, we will switch from $\HH$ to its
{\em intermediate subalgebra} $\HH^\flat\subset\HH$
generated by $T_i \ (i\ge 0)$, $X_a \ (a\in B)$,
and $Y_a \ (a\in B)$, where $B$ is any lattice
between $Q$ and $P$ (see \cite{C12}).
Accordingly, we replace $\Pi$ by the preimage $\Pi^\flat$ of $B/Q$ in
$\Pi.$ Generally, there can be two different
lattices $B_X$ and $B_Y$ for $X$ and $Y.$
We consider only $B_X=B=B_Y$ in the paper.
It is straightforward to check that all (anti-)automorphisms
introduced in the previous section preserve $\HH^\flat$.
The $\ddot{\Q}_{q,t}'$\~subalgebra
$\dHH^\flat\subset \HH^\flat$ is defined
accordingly.

We also set $\hW^{\flat}=B\cdot W\subset \hW$
and replace $m$ by the least $\tilde{m}\in\N$ such that
$\tilde{m}(B,B)\subset \Z$ in the definition of
$\Q_{q,t}, \Q_{q,t}', \ddot{\Q}_{q,t}'$, and $\ddot{\Q}_{q,t}^\dag$.

We point out that $\HH^\flat$ and the polynomial representation
can be defined over the ring $\Z[q^{\pm 1/(2\tilde{m})}, 
t_\nu^{\pm 1/2}]$.
However, the ring $\Q_{q,t}$ and its localizations above
will be sufficient in this paper.
\smallskip

The {\em Demazure-Lusztig operators} are
\begin{align}
&T_i\  = \  t_i^{1/2} s_i\ +\
(t_i^{1/2}-t_i^{-1/2})(X_{\al_i}-1)^{-1}(s_i-1),
\ 0\le i\le n;
\label{Demazx}
\end{align}
they obviously preserve $\Q_{q,t}[X]\equal\Q_{q,t}[X_b, b\in B]$.
We note that only the formula for $T_0$ involves $q$:
\begin{align}
&T_0\  = \ t_0^{1/2}s_0\ +\ (t_0^{1/2}-t_0^{-1/2})
(X_0 -1)^{-1}(s_0-1),\hbox{\ where\ }\notag\\
&X_0=qX_\vth^{-1},\
s_0(X_b)\ =\ X_bX_{\vth}^{-(b,\vth)}
 q^{(b,\vth)},\
\al_0=[-\vth,1].
\end{align}

The map sending $T_i \ (i\ge 0)$ to the corresponding operator from
(\ref{Demazx}), $X_b$ to $X_b$
(see (\ref{Xdex})) and
$\pi_r\mapsto \pi_r$ induces a
$ \Q_{ q,t}$\~linear
homomorphism from $\HH^\flat$ to the algebra of
linear endomorphisms
of $\Q_{ q,t}[X]$. We will extend the ring of constants here
to the field of fractions $\Q_{q,t}'$ or its subring
$\ddot{\Q}_{q,t}'$ of the rationals well defined when
$t_\nu=0$ for all $\nu$.

This $\HH^\flat$\~module is faithful
and remains faithful when $q,t$ take
any complex values assuming that
$q\neq 0$ is not a root of unity.
It will be called the
{\em polynomial representation};
the notation is
$
\v\equal \Q_{q,t}'[X].
$
We also set
$\ddot{\v}\equal \ddot{\Q}_{q,t}'[X]$
and
$\ddot{\v}^\dag\equal\ddot{\Q}_{q,t}^\dag[X].$

Given any $H\in \HH^\flat$, we continue to denote by $H$ the
corresponding operator in $\v$.

The polynomial representation can be described as
the $\HH^\flat$\~module induced from the one-dimensional
representation $T_i\mapsto t_i^{1/2}\ (i\ge 0),$
$Y_b\mapsto q^{(\rho_k,b)}$
of the affine Hecke subalgebra
$\h_Y^\flat=\lan T_i,Y_b\ran$.
\smallskip

Elements of $\HH^\flat$ act in $\v$
by {\em difference-reflection operators}, which are
operators of the form
\begin{align}
\sum_{w\in W,\,b\in B} g_{b,w}\,\Ga_b\,w, \ \
g_{b,w}\in\Q_{q,t}'(X),
\label{def-drop}
\end{align}
where $\Q_{q,t}'(X)$ is the field of rational functions
in the $X_b \ (b\in B)$. We denote the algebra of 
all such operators by $\cA$;
its defining relations are as follows:
\begin{align}
&q^{(a,b)}X_a\Ga_b=\Ga_b X_a,\ wX_a=X_{w(a)}w,\
w\,\Ga_b=\Ga_{w(b)}\,w.
\end{align}
The algebra of {\em difference operators} is the subalgebra of
$\cA$ generated by $\Q_{q,t}'(X)$ and $\Ga_b \ (b\in B)$.
There is a natural linear map
\begin{align}
\Red:\ \sum_{w\in W,\,b\in B} g_{b,w}(X)\,\Ga_b\,w
\mapsto\sum_{w\in W,\,b\in B}g_{b,w}\,\Ga_b,
\label{def-red}
\end{align}
sending difference-reflection operators to difference operators.
Clearly, $\Red$ is {\em not} a homomorphism of algebras.
\smallskip

The images of the $Y_b$ in the polynomial representation
are called the {\em difference Dunkl operators}.
Later we will make use of the following explicit
description of these operators. 
Let $b=\pi_r s_{j_l}\cdots s_{j_1}$
be a reduced decomposition of any $b\in B$,
and recall the definition
of $\ep_p$ from (\ref{def-ep-p}). Then
\begin{align}\label{ybsgn}
Y_b=\pi_r T_{j_l}^{\ep_l}\cdots T_{j_1}^{\ep_1}=
b\,G_{\tal^l}^{\hbox{\tiny\sgn}(\ep_l)} \cdots 
G_{\tal^1}^{\hbox{\tiny\sgn}(\ep_1)},
\end{align}
where $\sgn(\pm 1)=\pm$ and 
\begin{align}
\label{Gtal}
&G_{\tal}^+\equal t_{\al}^{1/2}+\frac{t_{\al}^{1/2}-t_{\al}^{-1/2}}
{X_{\tal}^{-1}-1}(1-s_{\tal})=
t_{\al}^{-1/2}(f_{\tal}+g_{\tal}\,s_{\tal}),\\
&\notag \ \ \ \ \, 
f_{\tal}=\frac{t_{\al} X_{\tal}^{-1}-1}{X_{\tal}^{-1}-1}, \ \
g_{\tal}=\frac{t_{\al}-1}{1-X_{\tal}^{-1}};\\
&G_{\tal}^-\equal t_{\al}^{-1/2}+
\frac{t_{\al}^{1/2}-t_{\al}^{-1/2}}{1-X_{\tal}}(1-s_{\tal})
=t_{\al}^{-1/2}(f_{\tal}-s_{\tal}\,g_{\tal}).
\comment{
&\notag \ \ \ \ \, f_{\tal}'=
\frac{t-X_{\tal}^{-1}}{1-X_{\tal}^{-1}}, \ \
g_{\tal}'=\frac{t-1}{X_{\tal}-1}.
}
\end{align}
Note that
$$
G_{\al_i}^+=s_i\,T_i,\ \
G_{-\al_i}^+=T_i\,s_i,\ \
G_{\al_i}^-=s_i\,T_i^{-1},\and
G_{-\al_i}^-=T_i^{-1}\,s_i.
$$
Let $\ddot{G}_{\tal}^\pm\equal t_{\tal}^{1/2}G_{\tal}^\pm$, so that
$\ddot{Y}_b=b\,\ddot{G}_{\tal^l}^{\hbox{\tiny\sgn}(\ep_l)}\cdots
\ddot{G}_{\tal^1}^{\hbox{\tiny\sgn}(\ep_1)}.$
\comment{
\begin{align}
\ddot{Y}_b=b\,\ddot{G}_{\tal^l}^{\hbox{\tiny\sgn}(\ep_l)}\cdots
\ddot{G}_{\tal^1}^{\hbox{\tiny\sgn}(\ep_1)}.
\end{align}
}
\medskip

\subsection{\bf Macdonald polynomials}
There are two equivalent definitions of the
{\em nonsymmetric Macdonald polynomials}, denoted
$E_b(X) = E_b^{(k)}$ for $b\in B$.

The first definition
is based on the {\em truncated theta function}:
\begin{align}
&\mu(X;t) = \mu^{(k)}(X)=\prod_{\al \in R_+}
\prod_{j=0}^\infty \frac{(1-X_\al q_\al^{j})
(1-X_\al^{-1}q_\al^{j+1})
}{
(1-X_\al t_\al q_\al^{j})
(1-X_\al^{-1}t_\al^{}q_\al^{j+1})}.\
\label{mu}
\end{align}
We will mainly consider $\mu$ as a Laurent series with
coefficients in  the ring $\Q[t_\nu][[q_\nu]]$ for
$\nu\in \nu_R=\{\nu_{\sht},\nu_{\lng}\}$.
The constant term of a Laurent series $f(X)$ will be denoted
by $\langle  f \rangle.$ One has
\begin{align}
&\langle\mu\rangle\ =\ \prod_{\al \in R_+}
\prod_{j=1}^{\infty} \frac{ (1- q^{(\rho_k,\al)+j\,\nu_\al})^2
}{
(1-t_\al q^{(\rho_k,\al)+j\,\nu_\al})
(1-t_\al^{-1}q^{(\rho_k,\al)+j\,\nu_\al})
}.
\label{consterm}
\end{align}
This equality is equivalent
to the Macdonald constant term conjecture
proved in complete generality in \cite{C2}.

Let $\mu_\circ\equal \mu/\langle \mu \rangle$.
The coefficients of the Laurent series $\mu_\circ$
belong to $\Q(q,t_\nu)$ for $\nu\in\nu_R$ and are
well defined at $t_\nu=0$.
Define an inner product on $\v$ by
\begin{align}
\label{innp}
\lan f, g\ran_\circ \equal \lan f g^\star \mu_\circ \ran,
\end{align}
where $\star$ is the $\Q$\~linear involution on $\v$ defined by
\begin{align}
\label{stv}
X_b^\star = X_b^{-1}, \ (q^{1/(2m)})^\star=q^{-1/(2m)},\
(t_\nu^{1/2})^\star=t_\nu^{-1/2}.
\end{align}
One has $\mu_\circ^\star=\mu_\circ$ and consequently 
$\lan g,f\ran=\lan f,g \ran^\star$.

The polynomials $E_b$ are uniquely determined from
the relations
\begin{align}
&E_b-X_b\ \in\ \oplus_{c\succ b}\Q_{q,t}' X_c,\
\lan E_b, X_c \ran = 0 \for B \ni c\succ b
\label{macd}
\end{align}
and for generic $q,t$; they form a basis of $\v$.
Their coefficients actually belong to $\Q(q,t_\nu)$
and are well defined at $t_\nu=0$; for the latter,
see (\ref{epolexists}) below.

This definition is due to Macdonald (for
$k_{\sht}=k_{\lng}\in \Z_+ $), who extended
the construction from \cite{O2}.
The general (reduced) case was considered in \cite{C4}.

We note that the $E_b$ satisfy the stronger condition
(see \cite[(2.7.5)]{M1}):
\begin{align}
\label{macsucc}
&E_b-X_b\ \in\ \oplus_{c\succq b}\Q_{q,t}' X_c,
\end{align}
in terms of the partial ordering $\succq$ defined in
(\ref{succpr}).
\smallskip

The second definition of the $E$\~polynomials
is based on the Dunkl operators:

\begin{proposition}
The polynomials $\{E_b , b\in B\}$ are the unique
(up to proportionality) eigenfunctions of the operators
$Y_a$ ($a\in B$) acting in $\v$:
\begin{align}
&Y_a(E_b)\ =\ q^{-(a,b_\#)}\,E_b\, \hbox{\ for\ }\,
b_\#\equal b- u_b^{-1}(\rho_k),
\label{Yone}\\
&\hbox{equivalently,\ \ }
\ddot{Y}_a(E_b)\ =\ q^{(a_+,\,\rho_k)-(a,b_\#)}\,E_b,
\notag
\end{align}
where $u_b=\pi_b^{-1}b$ is from Proposition \ref{PIOM},\
$b_\#=\pi_b(\!(-\rho_k)\!)$.\sq
\label{YONE}
\end{proposition}

The second definition readily leads to the orthogonality
of the $E$\~polynomials. 
This is due to the fact that
\begin{align}
\label{astinnp}
\lan H(f),g \ran = \lan f, H^\star(g) \ran,
\for f,g\in\v,\ H\in\HH^\flat,
\end{align}
where $\star$ is the anti-involution of $\HH^\flat$ extending
(\ref{stv}) and defined by
\begin{align}
\label{star}
\star:\ &T_i\mapsto T_i^{-1} \ (i>0),\ X_b \mapsto X_b^{-1},\
Y_b\mapsto Y_b^{-1},\ \pi_r\mapsto\pi_r^{-1}.
\end{align}

The norms of the $E$\~polynomials are given explicitly by
\begin{align}
\lan E_b,E_c \ran_\circ
=\de_{bc}\prod_{[\al,j]\in\la'(\pi_b)}
\frac{(1-q_\al^j t_\al^{-1} X_\al(q^{\rho_k}))
(1-q_\al^j t_\al X_\al(q^{\rho_k}))}
{(1-q_\al^j X_\al(q^{\rho_k}))
(1-q_\al^j X_\al(q^{\rho_k}))},
\label{epolnorms}
\end{align}
where we set $X_a(q^z)=q^{(a,z)}$ for $a\in B$; see \cite{C4}.
\smallskip

We note that
\begin{align}
\label{etastar}
(H(f))^\star=\eta(H)(f^\star), \for f\in\v,\ H\in\HH^\flat,
\end{align}
where $\eta$ is the involution defined in (\ref{eta}).
\medskip

\subsection{\bf Symmetric polynomials}
For $f\in \Q_{q,t}'[X]^W$, let
\begin{align}
\label{def-Lf}
\cL_f\equal\, &f(Y_{\om_1},\ldots,Y_{\om_n})=\notag
\sum_{w\in W,\,b\in P}g_{b,w}\,
\Ga_b\,w, \ \ g_{b,w}\in\Q_{q,t}'(X),\\
&L_f\equal\Red(\cL_f)=\sum_{w\in W,\,b\in P}g_{b,w}\,\Ga_b.
\end{align}
The operators $\cL_f$ and $L_f$ preserve $\v^W$
and they coincide upon the restriction to this space.
Moreover, the $L_f$ are 
$W$\~invariant difference operators, i.e.,
$w L_f w^{-1}=L_f$ for any $w\in W$.
If $f$ has coefficients in $\ddot{\Q}_{q,t}'$, we set
$\ddot{\cL}_f\equal f(\ddot{Y}_{\om_1},\ldots,\ddot{Y}_{\om_n})$
and $\ddot{L}_f=\Red(\ddot{\cL}_f)$.

Following Proposition \ref{YONE}, the
{\em symmetric Macdonald polynomials} $P_b=P_b^{(k)}$
can be introduced as eigenfunctions of these operators.
Explicitly,
\begin{align}
&L_f(P_{b_-})=f(q^{-b_-+\rho_k})\,P_{b_-},\ \ b_-\in B_-,\notag\\
&P_{b_-}=\sum_{b\in W(b_-)}X_b \mod 
\oplus_{c_-\succ b_-}\Q_{q,t}'X_c.\
\label{Lf}
\end{align}
These polynomials were introduced in \cite{M2};
For classical root systems, they were first used
in an unpublished work of Kadell.
In the case of $A_1$, they are due to Rogers.

For $a\in B$, let $\cL_a\equal\cL_f$ and $L_a\equal L_f$
where $f=\sum_{w\in W/W^a}Y_{w(a)}$; recall that $W^a$ is the
stabilizer of $a$ in $W$.
For these operators, (\ref{Lf}) reads
\begin{align}
&L_a(P_{b_-})=\bigl(\sum_{a'\in W(a)}
q^{-(a',b_--\rho_k)}\bigr)P_{b_-},\notag\\
&\ddot{L}_a(P_{b_-})=q^{(a_+,\rho_k)}
\bigl(\sum_{a'\in W(a)}q^{-(a',b_--\rho_k)}\bigr)P_{b_-}.
\label{dLa}
\end{align}

The connection between $E$ and $P$ is as follows:
\begin{align}
&P_{b_-}\ =\ \P_{b_+} E_{b_+}, \ b_-\in B_-,\
b_+=w_0(b_-),\label{symmetr}\\
&\P_{b_+}\equal\sum_{c\in W(b_+)} \ddot{T}_{w_c}
=\sum_{c\in W(b_+)} t^{l(w_c)/2}T_{w_c},\notag\\
&t^{l(\hw)}\equal \prod_{\nu}t_\nu^{l_\nu(\hw)}
\for  l_\nu(\hw)=|\{\tal\in\la(\hw) , \nu_\al=\nu\}|,\notag
\end{align}
where $w_c\in W$ is the element of least length such that
$c=w_c(b_+)$. Hence $P_{b_-}$ belongs to $\ddot{\Q}_{q,t}'[X]$. 
Taking the complete $t$\~symmetrization
$\P$ here (with the summation over all $w$), one obtains
$P_{b_-}$ up to proportionality. See \cite{O2,M4,C4}.
\medskip

There are two different kinds of inner products in
$\v$ from \cite{C101} and other works, with and without
using $\star$\,; we will mainly need the former in this work,
which is $\lan\,,\,\ran_\circ$ from (\ref{innp}).
In the symmetric
setting, they essentially coincide.
The inner products
of the symmetric polynomials
$P_b$ for $b=b_-$ read as follows (see \cite{C101}):
\begin{align}\label{normppols}
\lan P_b, P_c \ran_\circ =&\\
\de_{bc}\prod_{\al>0}\prod_{j=0}^{-(\al^\vee,b)-1}&
\frac{
(1-q_\al^{j+1}t_\al^{-1} X_\al(q^{\rho_k}))
(1-q_\al^j t_\al  X_\al(q^{\rho_k}))
}{
(1-q_\al^j X_\al(q^{\rho_k})) (1-q_\al^{j+1} X_\al(q^{\rho_k}))
}.\notag
\end{align}
\medskip

\subsection{\bf Using intertwiners}
\label{SEC:INT}
The $Y$\~intertwiners serve as creation operators
in the theory of nonsymmetric Macdonald polynomials.
Let $0\le i\le n$, $Y_0=Y_{\al_0}\equal q^{-1}Y_\vth^{-1}$.
Following 
\cite{C101} here and below, we set 
\begin{align}
&\Psi_i=\tau_+(T_i)+
\frac{t_i^{1/2}-t_i^{-1/2}}{Y_{\al_i}^{-1}-1},\ 
\Psi_i^b=
\tau_+(T_i) + \frac{t_i^{1/2}-t_i^{-1/2}}{
X_{\al_i}(q^{b_{\#}})-1}.
\label{Phijb}
\end{align}

Recall the definition of the automorphism $\tau_+$ of
$\HH^\flat$ from (\ref{taux}).
In the following theorem, 
which explicitly describes the action of the
intertwiners on the $E$\~polynomials,
we also need the pairing from (\ref{dform})
and the affine action $\,\hw\llb b\rrb\,$
from (\ref{afaction}).

\begin{theorem}\label{PHIEB}
Given $b\in B,\ 0\le i\le n$ such that
$(\al_i, b+d)> 0,$ 
\begin{align}\label{Phieb}
&q^{-\frac{(c,c)}{2}} E_{c} = 
t_i^{\frac{1}{2}}\Psi_i^b(q^{-\frac{(b,b)}{2}}E_b),\ \,
\frac{q^{-\frac{(b,b)}{2}} 
E_{b}}{\lan E_b,E_b\ran_\circ} =
t_i^{-\frac{1}{2}}\Psi_i^c\bigl(\frac{q^{-\frac{(c,c)}{2}} 
E_{c}}{\lan E_c,E_c\ran_\circ}\bigr),
\end{align}
where $c= s_i\llb b\rrb$.
If  $(\al_i, b+d)=0,$ then
\begin{align}
&\tau_+(T_i) (E_b) = t_i^{1/2} E_b,\ \tau_+(\ddot T_i) (E_b)
= t_i E_b,\ \ 0\le i\le n,
\label{Tjeco}
\end{align}
which results in the relations $s_i(E_b)=E_b$ as
$i>0$. For $c=\pi_r\llb b\rrb,$ where
the indices $\,r\,$ are from $O',$
\begin{align}
&q^{(b,b)/2-(c,c)/2}E_c\ =\ \tau_+(\pi_r)(E_b)\ =\
X_{\om_r}q^{-(\om_r,\om_r)/2}\pi_r(E_b).
\label{pireb}
\end{align}
Also $\tau_+(\pi_r)(E_b)\neq E_b$ for $\pi_r\neq\id$,
since $\pi_r\llb b\rrb\neq b$ for any $b\in B$.\sq
\end{theorem}
\smallskip

If $(\al_i, b)>0$ and $i>0$, then
the set $\la(\pi_c)$ is obtained from
$\la(\pi_b)$ by adding $[\al,(b_-,\al)]$
for $\al=u_b(\al_i)\in R_-$
and $(b_-,\al^\vee)=(b,\al_i^\vee)>0.$
When $i=0$ and $(\al_0, b+d)=-(b,\vth)+1>0$,
then the root $[\al,(b_-,\al)+1]$ is added to
$\la(\pi_b)$ for $\al= u_b(-\vth)=\al^\vee\in R_-$
and $(b_-, \al)=-(b,\vth)\ge 0.$

In each of these two cases,
$(\al_i,u_b^{-1}(\rho))=(\al,\rho)<0$ and
the powers of $t_\nu$ in
\begin{align}\label{xaliq}
&X_{\al_i}(q^{b_{\#}})=q^{(\al_i,b-u_b^{-1}(\rho_k)+d)}=
q^{(\al_i,b+d)}\prod_\nu t_\nu^{-(\al,\rho_\nu^\vee)}
\end{align}
are from $\Z_+$, with that of $t_i$ strictly
positive.
\medskip

Due to Theorem \ref{PHIEB}
(see also \cite{C1}, Corollary 5.3),
the coefficients of polynomial $E_b$ belong
to $\Q_{q,t}$ divided by 
\begin{align}\label{epolexists}
&\prod_{[-\al,\,\nu_\al j]\in \la(\pi_b)}
\bigl(
1- q_\al^{j}X_\al(q^{\rho_k})\bigr).
\end{align}

\comment{If $b\in B_-$ and the latter inequality holds for
$b_+=w_0(b)\in B_+,$ then the symmetric polynomial $P_b$ is 
well defined.}
More exactly, the ring $\Q_{q,t}$ can be replaced here by
$\Q[q,t_\nu]$,
which readily results in the existence
of the limits of the $E$\~polynomials when  $t_\nu=0$ for
all $\nu\in \nu_R$. Thus their coefficients become
polynomials in terms of $q$ in this limit.

Here the key is that powers of $q$ appearing in (\ref{epolexists}) 
are always multiplied by nonzero powers of $t_\nu$.
The same argument and a relatively straightforward analysis
of the leading $t$\~powers of the coefficients of $E_b$
can be used to see that the limits of $E_b$ exist
when $t_\nu\to\infty$. Moreover, their coefficients become 
polynomials
in $q^{-1}$ in this limit; see Corollary~\ref{EDAGCOEFF} below.


\comment{
\subsection{\bf Spherical polynomials}
The following renormalization
of the $E$-polynomials is of major importance in
the Fourier analysis (see \cite{C4}):
\begin{align}\label{ebebs}
\e_b\ \equal&\ E_b(X)(E_b(q^{-\rho_k}))^{-1},\where b\in B,
\\
E_{b}(q^{-\rho_k}) \ =&\ q^{(\rho_k,b_-)}
\prod_{[\al,j]\in \la\,'\,(\pi_b)}
\Bigl(
\frac{
1- q_\al^{j}t_\al X_\al(q^{\rho_k})
 }{
1- q_\al^{j}X_\al(q^{\rho_k})
}
\Bigr).\notag
\end{align}
This definition requires the $t$\~localization.
We call them nonsymmetric
{\em spherical polynomials\,}.
Formula (\ref{ebebs}) is the Macdonald
{\em evaluation conjecture}
in the nonsymmetric variant from \cite{C4}.
See \cite{C3} for the symmetric evaluation conjecture.
The following {\em duality formula} holds
for $b,c\in B\, :$
\begin{align}
&\e_b(q^{c_{\#}})\ =\ \e_c(q^{b_{\#}}),\
b_\# = b- u_b^{-1}(\rho_k),
\label{ebdual}
\end{align}
which is the main justification of the definition of $\e_b$.
Given $b\in B,$
the polynomial
$\e_b$
is well defined for $q,t\in \C^*$ if
\begin{align}
&\prod_{[\al,j]\in \la\,'\,(\pi_b)}
\bigl(
1- q_\al^{j}t_\al X_\al(q^{\rho_k})\bigr)\ \neq\ 0.
\label{esphexists}
\end{align}
\smallskip
In the symmetric setting,
\begin{align}\label{pebebs}
\p_b\ \equal&\ P_b(X)(P_b(q^{-\rho_k}))^{-1} \where b\in B_-\,,
\\
P_{b}(q^{-\rho_k})=P_{b}(q^{\rho_k}) \ =&\ q^{(\rho_k,b_-)}
\prod_{\al>0}\prod_{j=0}^{-(\al^\vee,b)-1}
\Bigl(
\frac{
1- q_\al^{j}t_\al X_\al(q^{\rho_k})
 }{
1- q_\al^{j}X_\al(q^{\rho_k})
}
\Bigr).\notag
\end{align}
The symmetric duality reads as follows:
\begin{align}
&\p_b(q^{c-\rho_k})\ =\ \p_c(q^{b-\rho_k}),\
\for b,c\in B_-\,.
\label{pebdual}
\end{align}
\medskip
The norm formula becomes entirely conceptual:
\begin{align}\label{normpipols}
&(\lan \p_{b}(X)\p_{b}(X^{-1})\mu_\circ\ran)^{-1}\
=\ \sum_{a\in W(b)}\,\mu(\pi_{a})\mu(\hbox{id})^{-1},\\
& \where \mu(\hw)\equal\mu(\hw(\!(q^{-\rho_k})\!))
\for \hw\in \hW.\notag
\end{align}
It is a direct corollary of the fact that the Fourier
transform sends the $\p$\~polynomials to the delta-functions;
see \cite{C101}.
\medskip
}
\medskip

\subsection{\bf The limit
\texorpdfstring{$t\to$ 0}{\em t=0}}
The limit (reduction) of $\dHH^\flat$ introduced
in Definition \ref{DEFHHDOT} when $t_\nu=0$ for all
or some $\nu$ will be denoted by $\bHH^\flat$ and 
called the {\em nil-DAHA}. For the sake of definiteness,
{\em all\,} $t_\nu=0$ {\em in this section\,}.

The polynomials $\overline{E}_b$, the images of $E_b$
as $t_\nu=0$ for $\nu\in \nu_R$ 
linearly generate the bar-polynomial representation:
$$
\overline{\v}\equal\Q_q'[X_b , b\in B],\where
\Q_q'\equal\Q(q^{1/(2m)}).
$$
Thus, the $\bHH^\flat$ action in $\overline{\v}$ is given by the 
operators
\begin{align}
&\overline{T}_i\equal\ddot{T}_i(t_i=0)=(X_{\al_i}-1)^{-1}(1-s_i),\
0\le i\le n,\\
&\overline{Y}_a\equal\ddot{Y}_a(t_\nu=0),\ a\in B,
\end{align}
and the action of $\Pi$ and multiplication by $X_b \ (b\in B)$.

The $\overline{E}$\~polynomials are eigenfunctions of the
$\overline{Y}$\~operators.
Using (\ref{Yone}), one has explicitly:
\begin{align}
\label{YbarE}
\overline{Y}_a(\overline{E}_b)=
\begin{cases}
q^{-(a,b)}\overline{E}_b, &\text{if $u_b(a)=a_-$,}\\
0, &\text{otherwise.}
\end{cases}
\end{align}
Note that using only $\overline{Y}_a\, (a\in B_+)\,$ is obviously
not sufficient to split $\{\overline{E_b}\}$ for generic
$q$; all $a\in B$ must be involved. 

Theorem \ref{PHIEB} holds under this limit
and gives quite a constructive approach to the
$\overline{E}$\~polynomials. The
reductions $\overline{\Psi}_i^c$ of the intertwiners
$\ddot{\Psi}_i^c\equal t_i^{1/2}\Psi_i^c$ from
(\ref{Phijb}) can be used
to generate $\overline{E}_b$; these intertwiners become
$\tau_+(\overline{T}_i)+1$ in this limit.

This simplification is directly connected with the fact that
$$\overline{T}_i'\equal
\ddot{T}_i'(t_i=0)\ =\ \overline{T}_i+1
$$
satisfy the same
homogeneous Coxeter relations as $\{T_i,\, 0\le i\le n\}$ do,
a special feature of the nil-DAHA. It readily results from the
theory of intertwiners, and, of course, can be checked directly
as well.

Let us provide some details of the construction
of bar-polynomials. The action
of $\pi_r$ on $\{\overline{T}_i'\}$ by conjugation
obviously remains unchanged. Thus relations
(i,ii) from Definition \ref{double} (and above)
hold for  $\{\overline{T}_i'\}$. Therefore,
given $\hw\in \hW$,
the element $\overline{T}'_{\hw}=\pi_r
\overline{T}'_{i_l}\cdots \overline{T}'_{i_1}$
does not depend on the choice of the reduced
decomposition $\hw=\pi_r s_{i_l}\cdots s_{i_1}$.

For instance, the operators
$\overline{\Pi}'_i\equal\tau_+(\overline{T}'_{-\om_i})$ for
$i=1,\ldots,n$ are pairwise commutative and, importantly,
are $W$\~invariant.

Indeed, one has\,:\,
$\overline{\Pi}'_{b}=\prod_{i=1}^n\,
(\overline{\Pi}'_i)^{n_i}$ for
$B_-\ni b=-\sum n_i\,\om_i$. Provided that all $n_i>0$,
the {\em reduced} decomposition $b=b_-=w_0\pi_{b_+}$ holds
for the longest element $w_0\in W$ and $b_+=w_0(b)\in B_+$;
see (\ref{lupiw}).
Thus $\overline{\Pi}'_b$ is divisible on the left 
by $(\overline{T}_i+1)$
for any $i>0$ and therefore divisible by the $W$\~symmetrizer on the
left. It results in the $W$\~invariance of $\overline{\Pi}'_b$
provided that $b\in B_-$.

The $W$\~invariance of $\{\overline{\Pi}'_b,\,b\in B_-\}$ simplifies
significantly the relation of the $\overline{E}$\~polynomials to
the $\overline{P}$\~polynomials. It becomes
\begin{align}\label{eplimbar}
&\overline{P}_{b}\ =\ \overline{E}_{b} \mbox{\ \ for\ \ }
b=b_-\in B_-\,.
\end{align}

We come to the following explicit proposition.

\begin{proposition} \label{BAREP}
In the representation $\overline{\v}$ of $\bHH^{\flat}$, 
\begin{align}\label{barpform}
&\tau_+(\overline{T}'_{\hw})(1)= 
q^{-(b,b)/2}\,\overline{E}_b \for \hw\ =\ \pi_b,\, b\in B,\notag\\
&\overline{\Pi}'_b(1)\ =\ q^{(b,b)/2}\overline{P}_b
\ =\ q^{(b,b)/2}\overline{E}_b
\for b\in B_-,
\end{align}
where $\overline{\Pi}'_i$ can be replaced by their restrictions
$\hbox{\rm Red}_{\,W}(\overline{\Pi}'_i)$ to $\overline{\v}^{\,W}$,
which are pairwise commutative $W$\~invariant
difference operators.\sq
\end{proposition}
\smallskip

Define $\Si(b)$ (resp. $\Si_\ast(b), \Si_+(b)$) to be the span of
those monomials $X_c$ with $c\in\si(b)$ 
(resp. $\si_\ast(b), \si_+(b)$); cf. (\ref{cones}).
Recall the partial ordering $\succeqq$ from (\ref{succpr}).
\comment{
Consider the spaces
\begin{align}
\v(b_-)\equal\Si(b_-)/\Si_+(b_-),\
\overline{\v}(b_-)\equal\ddot{\v}(b_-)\mid_{t_\nu\to 0}.
\end{align}
Thus $\v(b_-)$ and $\overline{\v}(b_-)$
have a basis consisting of the images of monomials
labeled by the weights $w(b_-) \ (w\in W)$; we continue to denote
these by $X_{w(b_-)}$.
}

\begin{proposition}
For any $b\in B$, one has
\begin{align}
\overline{E}_b &= \sum_{c\succeqq b,\, c\in W(b)} X_c\mod \Si_+(b).
\label{ebarmodsi}
\end{align}
\end{proposition}

\proof
We argue by induction on $l(w_b)$; recall that $w_b$ is, 
by definition,
the unique element of $W$ of shortest length such that $w_b(b_+)=b$.
For $b=b_+$, (\ref{ebarmodsi}) clearly holds.
Before the inductive step, we remark that the intertwiners
$\overline{T}_i' \ (i>0)$ preserve $\Si_+(b)$.
Modulo $\Si_+(b)$, one has for $i>0$ that
\begin{align}
\overline{T}_i'(X_b)=
\begin{cases}
X_{s_i(b)}+X_b, & \text{if }(b,\al_i)>0,\\
X_b, & \text{if }(b,\al_i)=0,\\
0, & \text{if }(b,\al_i)<0.
\label{intmodsi}
\end{cases}
\end{align}

Now suppose $l(w_b)>1$ and choose any $i>0$ such that
$l(s_i w_b)<l(w_b)$. Then $w_b=s_i w_{c}$ is reduced for
$c=s_i(b)=s_i w_b(b_+)\succeqq b$. One has $(c,\al_i)>0$
and hence $\overline{E}_b=\overline{T}_i'(\overline{E}_c)$.

We must show that for any $w\leq w_b$, the monomial $X_{w(b_+)}$
appears in $\overline{E}_b$ with coefficient $1$.
Now $w \leq w_b$ implies

\centerline{
($i$)\, $w\leq w_c$ \ \ or\ \ 
($ii$)\, $s_i w\leq w_c$ \ (or both).}

Assuming ($i$), $\overline{E}_c$ contains $X_{w(b_+)}$ with
coefficient $1$, by induction. Then either $(w(b_+),\al_i)\geq 0$,
in which case (\ref{intmodsi}) applies directly,
or $(w(b_+),\al_i)<0$. In the latter case,
$s_i w < w\leq w_c$ and one applies (\ref{intmodsi})
to $X_{s_i w(b_+)}$.

Case ($ii$) can be handled by a similar argument.
\sq
\smallskip

\rmk
For the affine root systems considered in this paper
(with $\al_0$ defined in terms of the maximal {\em short} 
root $\vth$),
a connection was established between the polynomials
$\overline{E}_b$ and the level-one Demazure characters
of the corresponding irreducible affine Lie algebras; see
\cite{San} and, especially, Theorem 1 from \cite{Ion1}.
Paper \cite{Ion1} is based on the technique of
intertwiners (from \cite{KnS} in the $GL_n$\~case
and \cite{C1} for arbitrary reduced root systems).

It is important that only positive powers of $q$ appear in the
coefficients of $\overline{E}_b$ (see the discussion following
(\ref{epolexists}) above). In fact, the coefficients of these
$q$\~polynomials are non-negative. One can obtain
it from the interpretation via Demazure characters or
using the intertwiners (we are going to discuss the latter in further
papers). As $q\to 0$, the polynomials $\overline{P}_{b_-}$ become
the classical finite dimensional Lie characters,
which can be seen, for instance, from (\ref{normppolsbar})
below.
\sq
\smallskip

Concerning the orthogonality of the
$\overline{E}$\~polynomials,
the $\mu$\~function from (\ref{mubar}) becomes
\begin{align}
&\overline{\mu}\equal\mu(t_\nu=0)=\prod_{\al \in R_+}
\prod_{j=0}^\infty (1-X_\al q_\al^{j})
(1-X_\al^{-1}q_\al^{j+1}).
\label{mubar}
\end{align}
The constant term formula becomes a well known identity:
\begin{align}
&\langle\overline{\mu}\rangle\ =\ \prod_{i=1}^{n}
\prod_{j=1}^{\infty} \frac{1}
{1-q_i^j}.
\label{constermbar}
\end{align}

The polynomials $\overline{E}_b$ can be uniquely determined
from the relations (\ref{macd}) as $t_\nu=0$:
\begin{align}
&\overline{E}_b-X_b\ \in\ \oplus_{c\succ b}\Q_q' X_c,\
\lan \overline{E}_b X_{c}^{-1}
\overline{\mu}\ran = 0 \for B \ni c\succ b.
\label{macdbar}
\end{align}

To state the counterpart of the norm formula (\ref{epolnorms})
as $t_\nu=0$ we will need the limits
$\overline{E}_b^\dag$ of the $E$\~polynomials as 
$t_\nu\to\infty$. More generally, we set
$\overline{f}^\dag\equal 
\lim_{t_\nu\to \infty} f$ for any Laurent
polynomial or series depending on $q,t_\nu$,
provided the existence of this limit.
Using the conjugation $\star$ from (\ref{stv}) 
(sending $t_\nu^{1/2}$ to 
$t_\nu^{-1/2}$),
\begin{align}
\label{ebastbar}
\overline{(f^\star)}=(\overline{f}^{\,\dag}\,)^\ast,
\ \ \ \ \overline{(f^\star)}^{\,\dag}=(\,\overline{f}\,)^\ast,
\end{align}
where $X_b^\ast\equal X_b^{-1},\ (q^{1/(2m)})^\ast\equal q^{-1/(2m)}$.

Using this notation, the limit of the norm formula
from (\ref{epolnorms}) as $t\to 0$ is as follows:
\begin{align}\label{ebarnorms}
&\lan \overline{E}_b, \overline{E}_c \ran_\circ\equal
\lim_{t\to 0}\,\lan E_b,E_c \ran_\circ\\
&=
\lan \overline{E}_b\, (\overline{E}_c^\dag)^\ast\, 
\overline{\mu}_\circ\ran_\circ
=\de_{bc}\prod_{[\al,j]} (1-q_\al^j),
\notag
\end{align}
where the last product runs over all $[\al,j]\in\la'(\pi_b)$ with 
simple $\al=\al_i\in R_+$. 
Use that $(\rho_k,\al)=(\rho_k,\nu_\al \al^\vee)$
and $(\rho_k,\al_i^\vee)=k_i\,$;  recall that 
$q_\al=q^{\nu_\al}$ and $t_\al=q_\al^{k_\al}$.

Assuming that $q$ is not a root of unity,
these formulas readily provide the existence of the
polynomials $\overline{E}_b^\dag$ for any $b\in B$
and that they form a basis of $\overline{\v}$. This holds in fact
for any nonzero $q$, but the justification requires a different
approach. For instance, we use (\ref{epolexists}) in the
limit $t_\nu\to \infty$ to establish Corollary~\ref{EDAGCOEFF} below.
\smallskip

The relation between the limits $t_\nu\to 0$ and
$t_\nu\to \infty$ goes through the general formula
\begin{align}
&E_b^\star =\
\prod_{\nu\in \nu_R} t_\nu ^{l_\nu(u_b)- l_\nu(w_0)/2}\,
T_{w_0}(E_{\varsigma(b)}), \where\label{ebast}
\varsigma(b)=-w_0(b),
\end{align}
from \cite{C101} and other first author's works.
This connection becomes especially
simple for the symmetric polynomials:
$$
P_{b}(X)^\star=P_{b}(X^{-1}) \for b=b_-,\
\overline{P}_{b}=P_{b}(t_\nu\to 0)=P_{\varsigma(b)}(t_\nu\to \infty).
$$
We use that $P_b(X^{-1})=P_{\varsigma(b)}(X)$. 
This formula readily follows from 
the relations $\hw(\mu)/\mu=(\hw(\mu)/\mu)^\star$
for any $\hw\in \hW$ if one uses the definition of $\{P_b\,\}$ 
via $\mu$. The ratios $\hw(\mu)/\mu$ are rational functions in
terms of $X_{\tal},q,t_\al$, so the conjugation (applying $\star$)
is well defined.
\smallskip

For $b,c\in B_-\,$,
the norm formula from (\ref{normppols}) reads as:
\begin{align}\label{normppolsbar}
&\lan \overline{P}_{b}(X)\overline{P}_{c}(X^{-1})
\overline{\mu}_\circ\ran\
=\ \de_{bc}\prod_{i=1}^n\prod_{j=1}^{-(\al_i^\vee,b)}
(1-q_i^{j})\,.
\end{align}
We use that $P_b^\star=P_{\varsigma (b)}=P_b(X^{-1})$
for $b\in B_-$, which makes defining $\overline{P}_b^\dag$
unnecessary in this case.
\medskip

\subsection{\bf The limit
\texorpdfstring{$t\to \infty$}{at infinity}}
\label{sec:infty}
Let us discuss the limits $\overline{E}_b^\dag$ of the
$E$\~poly\-nomials more systematically. As orthogonal polynomials,
they can be introduced using (\ref{ebarnorms}); let us 
outline an approach based on the real integration instead
of taking the constant term.

We will use that 
$$
\mu^\dag\equal \mu(X;q,t^{-1})^{-1}=\prod_{\tal\in \tR_+}
\frac{1-t_\al^{-1}X_{\tal}}{1-X_{\tal}}
$$
satisfies the relations $\hw^{-1}(\mu)/\mu=
\hw^{-1}(\mu^\dag)/\mu^\dag$ for $\hw\in \hW$.
Thus $\mu/\mu^\dag$ is (formally) $\hW$\~invariant
and either function can be used for the corresponding
orthogonal polynomials and operators depending on the
setting of the theory. 
See formula (2.7) from Part I of \cite{ChM}. Recall that
$X_{\tal}=q_{\al}^jX_\al$ for $\tal=[\al,\nu_\al j]\in \tR$.

In the limit $t_\nu\to\infty$, 
\begin{align}
&\overline{\mu}^\dag\ = \ \prod_{\al \in R_+}
\prod_{j=0}^\infty \frac{1}{(1-X_\al q_\al^{j})
(1-X_\al^{-1}q_\al^{j+1})}.
\label{mubard}
\end{align}
We need to replace the constant term functional
$\lan f \ran$ by 
$$
\int_\ep\,f(X)\equal \sum_{w\in W}
\int_{w(\imath\ep+\R^n)}f(q^x)\,dx 
\for \ep\in \R^n,
$$
provided that $(\ep,\al)>0$ for $\al\in R_+$.
We set $X=q^x$, $X_b=q^{(x,b)}$, and
$$
\int_{\imath\ep+\R^n}(\ .\ )dx\ =\
\int_{\imath\ep_n-\infty}^{\imath\ep_n+\infty}
\cdots\int_{\imath\ep_1-\infty}^{\imath\ep_1+\infty}\,(\ .\ )\,
dx_1\cdots dx_n, 
$$ 
where $x_i=(x,\om_i)$; the integration contours $w(\imath\ep+\R^n)$ 
are the images of $\imath\ep+\R^n$.
The integral $\int_\ep \overline{\mu}^\dag$ is connected 
with the Appell-Lerch sums and will not be discussed here;
it does not depend on the choice of $\ep$. Cf. Section 2.3.5,
``Etingof's theorem", from \cite{C101}.

Now $\{\overline{E}_b^\dag\}$ can be introduced by means of the
relations
\begin{align}
&\overline{E}_b^\dag-X_b\ \in\ \oplus_{c\succ b}\Q_q' X_c,\
\int_\ep \overline{E}^\dag_b 
X_{c}^{-1}\overline{\mu}^\dag=0 \for B \ni c\succ b.
\label{macdbardag}
\end{align}

Setting $\overline{\mu}^\dag_\ep =\overline{\mu}^\dag/
\int_{\ep}\overline{\mu}^\dag$,
the norm formula (\ref{epolnorms}) becomes
\begin{align}\label{ebarnormd}
&\lan \overline{E}_b^\dag, \overline{E}_c^\dag \ran_\ep
\equal
\int_\ep \overline{E}_b^\dag\, (\overline{E}_c)^\ast\, 
\overline{\mu}_\ep^\dag
=\de_{bc}\prod_{[\al,j]} (1-q_\al^{-j}),
\end{align}
where $[\al,j]\in\la'(\pi_b)$ for simple $\al$.

Calculating the integrals here can be reduced to taking the
constant term. Namely,
$$
\int_\ep \overline{E}_b^\dag\, (\overline{E}_c)^\ast\, 
\overline{\mu}_\ep^\dag\ =\ 
\lan\,\overline{E}_b^\dag\, (\overline{E}_c)^\ast\, 
\overline{\mu}_\circ^\ast\,\ran \for b,c\in B.
$$
This follows from the relation 
$\overline{\mu}\,\overline{\mu}^\dag=1$, which connects the
action of $\hW$ on 
$\overline{\mu}_\circ$ and $\overline{\mu}^\dag_\ep$. 
Formally,  
$\overline{\mu}_\circ^\ast/\overline{\mu}^\dag$
is a $\hW$\~invariant function. We conclude that
(\ref{ebarnormd}) and (\ref{ebarnorms}) 
result in coinciding families of polynomials.
\smallskip

Proposition \ref{EDAG} below uses the intertwining operators
to establish the existence of $\{\overline{E}_b^\dag\}$ 
in a more direct way; it also relates them to the
$\overline{E}$\~polynomials when $b\in B_-$.
In Corollary \ref{EDAGCOEFF}, we show that the coefficients
of $\{\overline{E}_b^\dag\}$ belong to $\Z[q^{-1}]$.

We set
\begin{align}
\ddot{T}_i^\dag &\equal t_i^{-1/2}T_i,\notag
&&(\ddot{T}_i^\dag)' \equal t_i^{-1/2}T_i^{-1},\notag \\
\overline{T}_i^\dag &\equal\ddot{T}_i^\dag(t_i=\infty),\notag
&&(\overline{T}_i^\dag)' \equal(\ddot{T}_i^\dag)'(t_i=\infty)=
\overline{T}_i^\dag-1,\\
\overline{T}_i^\dag &=(s_i+1)\frac{1}{1-X_{\al_i}}, 
&&(\overline{T}_i^\dag)' =\frac{X_{\al_i}}{X_{\al_i}-1}(s_i-1)
\hbox{\ \ \, in\ } \overline{\v}. \notag
\end{align}
Recall that $\overline{\v}=\Q_q'[X_b , b\in B]$ as a space; 
by $\overline{\v}^\dag$, we will mean this space with the
action of $\overline{T}_i^\dag$ and other $\dag$\~operators.

Correspondingly, we set $\ddot{Y}_a^\dag=q^{-(a_+,\rho_k)}Y_a$.
Then it is straightforward to see that the limit
\begin{align}
\overline{Y}_a^\dag\equal\lim_{t\to 0}\ddot{Y}_a^\dag
\end{align}
exists.
Using (\ref{Yone}), we arrive at the analog of
(\ref{YbarE}):
\begin{align}
\label{YbarEdag}
\overline{Y}_a^\dag(\overline{E}_b^\dag)=
\begin{cases}
q^{-(a,b)}\overline{E}_b^\dag &\text{if $u_b(a)=a_+$,}\\
0 &\text{otherwise.}
\end{cases}
\end{align}

\begin{proposition}
\label{EDAG}
(i) For $b\in B_-$,
\begin{align}
\overline{E}_b^\dag\ &=\ q^{(\rho,b)}
(\overline{T}_{\pi_\rho}^{\,\prime}(\overline{E}_{-w_0(b)}))^\ast,
\label{ebminusdag}
\end{align}
where $\overline{T}_{\pi_\rho}^{\,\prime}=
\overline{T}_{j_1}^{\,\prime}\cdots
\overline{T}_{j_l}^{\,\prime}\,\pi_r^{-1}$
is defined for a reduced decomposition
$\pi_\rho=\pi_r s_{j_l}\cdots s_{j_1}\, (r\in O)$ 
and does not depend on the choice of such a decomposition.

(ii) In the opposite direction, the polynomials 
$\{\overline{E}_b^{\,\ast}\}$
for $b \in B_-$ can be obtained from $\overline{E}_{-w_0(b)}^\dag$
as follows:
\begin{align}\label{Edagbminus}
\overline{E}_b^{\,\ast}\ =\ 
\overline{T}_{w_0}^\dag(\overline{E}_{-w_0(b)}^\dag).
\end{align}
In particular, the polynomials on the right-hand side are 
$W$\~invariant and nonzero for any $b\in B_-$
with coefficients in $\Z_+[q^{-1}]$.

(iii) If $(b,\al_i)<0 \ (1\leq i\leq n)$, then
\begin{align}
\overline{E}_{s_i(b)}^\dag\ &=\
\begin{cases}
(1-q^{(b,\al_i)})^{-1}
(\overline{T}_i^\dag)'(\overline{E}_{b}^\dag) &
\text{if }(u_b(\al_i),\rho^\vee)=1,\\
(\overline{T}_i^\dag)'(\overline{E}_{b}^\dag) &
\text{if }(u_b(\al_i),\rho^\vee)>1.
\end{cases}
\label{daginteb}
\end{align}
Combining this with (\ref{ebminusdag}),
we obtain another proof of the existence of
$\overline{E}_b^\dag$ for any $b\in B$, which holds 
for any $q\neq 0$ due to (\ref{epolexists}) 
(see Corollary \ref{EDAGCOEFF} below).

(iv) Assuming that $(b,\al_i)>0 \ (1\leq i\leq n$), we set
$\overline{T}_{i,b}^\dag=\lim_{t_\nu\to\infty}\ddot{T}_{i,b}^\dag$
for 
$\ddot{T}_{i,b}^\dag=q^{-(b,b_\sharp)}t_i^{1/2}T_iY_b^{-1}$. Then
\begin{align}
&\overline{E}_{s_i(b)}^\dag=\label{daginteb-rev}
\begin{cases}
(\overline{T}_{i,b}^\dag+q^{-(b,\al_i)})(\overline{E}_b^\dag) &
\text{if }(u_b(\al_i),\rho^\vee)=-1,\\
\overline{T}_{i,b}^\dag(\overline{E}_b^\dag) &
\text{if }(u_b(\al_i),\rho^\vee)<-1.
\end{cases}
\end{align}

\comment{
(v) Suppose $b$ is antidominant. Let
$$ \tilde{T}_{0,b}^\dag=\lim_{t_\nu\to\infty}
t_0^{1/2}q^{(\vth,b_\sharp)}T_0^{-1}Y_\vth. $$
Then
\begin{align}
&\overline{E}_{s_0\llb b\rrb}=
X_0^{-1}\tilde{T}_{0,b}^\dag(E_b).
\end{align}
}

%
\end{proposition}

\proof
Claim $(i)$. Let $b=b_-$.
We will normalize (\ref{ebast}) to prove (\ref{ebminusdag}).
Note that $q^{(c,\,w_0(b)+\rho_k)}Y_c^{-1}$
acts as the identity on $E_{-w_0(b)}$
for any $c\in B$.
Taking $c=c_+$, so that $l_\nu(c)=2(c,\rho_\nu^\vee)$,
one therefore has
\begin{align}
E_b^\star &=q^{(c,w_0(b))}\prod_\nu t_\nu^{-l_\nu(w_0)/2+l_\nu(c)/2}
T_{w_0}Y_c^{-1}(E_{-w_0(b)}).
\label{E-ast-renorm}
\end{align}
Specializing further to $c=\rho$, we have
$Y_\rho=T_{\pi_\rho}T_{w_0}$
and (\ref{E-ast-renorm}) becomes
\begin{align}
E_b^\star &=q^{-(\rho,b)}\prod_\nu t_\nu^{l_\nu(\pi_\rho)/2}
T_{\pi_\rho}^{-1}(E_{-w_0(b)})=q^{-(\rho,b)}
\ddot{T}_{\pi_\rho}^{\,\prime}(E_{-w_0(b)}),
\end{align}
where by definition
$\ddot{T}_{\pi_\rho}^{\prime}=\ddot{T}_{j_1}^{\prime}\cdots
\ddot{T}_{j_l}^{\prime}\pi_r^{-1}$
for any reduced decomposition $\pi_\rho=\pi_r s_{j_l}\cdots s_{j_1}$.
Moving $\star$ to the right-hand side and taking
$t_\nu\to\infty$, we obtain (\ref{ebminusdag}).

Claim $(ii)$. This is immediate from (\ref{ebast}), which is
(\ref{E-ast-renorm}) for $c=0$. Note that (\ref{Edagbminus})
modulo $\Si_+(b)$ results in (\ref{ebarmodsi}) for $b=b_-$.

\comment{
We again use (\ref{ebast}), which gives that
$\overline{E}_{-w_0(b)}^{\,\ast}=
\overline{T}_{w_0}^\dag(\overline{E}_b^\dag)$
when $b=b_-$ and $t_\nu\to\infty$. Now since
$\overline{E}_{-w_0(b)}=\overline{P}_{-w_0(b)}$,
the claim follows.
The fact that
$\overline{T}_{w_0}^\dag(\overline{E}_b^\dag)$
is symmetric is immediate from
$(\overline{T}_i^\dag)'\overline{T}_{w_0}^\dag=0$.
}

Claim $(iii)$. This follows from a modification of 
(\ref{Phijb}) and (\ref{Phieb}).
It is convenient to use the normalized intertwiners 
$\mathfrak{G}_i\equal\psi_i^{-1}\Psi_i$ for
\begin{align}\label{normPsiinter}
\Psi_i &=\tau_+(T_i)+
\frac{t_i^{1/2}-t_i^{-1/2}}{Y_{\al_i}^{-1}-1},\ \
\psi_i =t_i^{1/2}+
\frac{t_i^{1/2}-t_i^{-1/2}}{Y_{\al_i}^{-1}-1}.
\end{align}
For simplicity, let us take here $1\leq i\leq n$.
In addition to the braid relations, the
normalized intertwiners satisfy $\mathfrak{G}_i^2=1$.
Hence $\Psi_i^{-1}=\psi_i^{-1}\Psi_i\psi_i^{-1}$.
Now, when $(b,\al_i)<0$, (\ref{Phieb}) gives
$E_{s_i(b)}=t_i^{-1/2}\Psi_i^{-1}(E_{b})$. Applying 
$\Psi_i^{-1}=\psi_i^{-1}\Psi_i\psi_i^{-1}$
to $E_b$, the first $\psi_i^{-1}$ produces
\begin{align*}
&\Bigl(
t_i^{1/2}+\frac{t_i^{1/2}-t_i^{-1/2}}{q^{(\al_i,b_\sharp)}-1}
\Bigr)^{-1}
=t_i^{1/2}\frac{q^{(\al_i,b_\sharp)}-1}
{q^{(\al_i,b_\sharp)}t_i-1}.
\end{align*}
Moving $\psi_i^{-1}$ through $\Psi_i$ changes $Y_{\al_i}^{-1}$
to $Y_{\al_i}$. Thus, the second (left) $\psi_i^{-1}$ produces 
the factor
\begin{align*}
&\Bigl(
t_i^{1/2}+\frac{t_i^{1/2}-t_i^{-1/2}}{q^{-(\al_i,b_\sharp)}-1}
\Bigr)^{-1}
=t_i^{-1/2}\frac{q^{-(\al_i,b_\sharp)}-1}
{q^{-(\al_i,b_\sharp)}-t_i^{-1}}
=t_i^{-1/2}\frac{q^{(\al_i,b_\sharp)}-1}
{q^{(\al_i,b_\sharp)}t_i^{-1}-1}.
\end{align*}
Multiplying these two factors and taking $t_\nu\to\infty$,
one arrives at (\ref{daginteb}).
Note that $u_b(\al_i)>0$ and hence $q^{(\al_i,b_\sharp)}$ contains
non-positive powers of $t_\nu$ and at least one $t_i^{-1}$.

Formula (\ref{daginteb}) can be directly deduced from
(\ref{Phieb}) and the norm formula (\ref{epolnorms})
in the limit $t\to\infty$.
Indeed, for $i>0$ and due to the inequality $(b,\al_i)<0$,
\begin{align}\label{Phiebx}
&\frac{E_{c}}{\lan E_c,E_c\ran_\circ} =
t_i^{-\frac{1}{2}}\Psi_i^b\bigl(
\frac{E_{b}}{\lan E_b,E_b\ran_\circ}\bigr).
\end{align}
Then $\,\lim_{t\to\infty} t_i^{-\frac{1}{2}}\Psi_i^b=
\overline{T_i}^\dag-1=(\overline{T_i}^\dag)'\,$ and 
\begin{align*}
&\lim_{t\to\infty}\lan E_b,E_b\ran_{\circ}
=\prod_{[\al,j]} (1-q_\al^{-j})=
(1-q^{(b,\al_i)})\lim_{t\to\infty}\lan E_c,E_c\ran_{\circ}.
\end{align*}
The product here is over 
$[\al,j]\in\la'(\pi_b)$ for simple $\al$.
Cf. (\ref{ebarnormd}) and recall that $q_{\al}=q^{\nu_\al}=
q^{2/(\al,\al)}$. 
\smallskip

Claim $(iv)$. Using (\ref{Phieb}) and (\ref{Yone}),
we can write
\begin{align}
&E_{s_i(b)}=t_i^{1/2}
\Bigl(T_i+\frac{t_i^{1/2}-t_i^{-1/2}}{q^{(\al_i,b_\sharp)}-1}\Bigr)
q^{-(b,b_\sharp)}Y_b^{-1}(E_b).
\label{Phieb-adj}
\end{align}

We claim that
$t_i^{1/2}q^{-(b,b_\sharp)}T_iY_b^{-1}$
is well defined as $t_\nu=\infty$.
Indeed, $$ l_\nu(b)=2(b_+,\rho_\nu^\vee)=-2(b_-,\rho_\nu^\vee) $$
and since $(b,\al_i)>0$, $u_b$ has a reduced decomposition
ending in $s_i$. 
Hence (\ref{def-ep-p}) gives the claim.

Now taking $t_\nu\to\infty$ in (\ref{Phieb-adj}),
while noting that
\begin{align}
\lim_{t_\nu\rightarrow\infty}
\frac{t_i-1}{q^{(\al_i,b_\sharp)}-1}&=
\begin{cases}
q^{-(\al_i,b)},&\text{if }(u_b(\al_i),\rho^\vee)=-1,\\
1,&\text{if }(u_b(\al_i),\rho^\vee)<-1.
\end{cases}
\end{align}
we arrive at (\ref{daginteb-rev}).
\sq
\smallskip

\begin{corollary}
\label{EDAGCOEFF}
(i) For $b\in B_-$, one has
\begin{align}\label{nbce}
\overline{E}_b^\dag\!= X_b +
\!\!\sum_{W(b) \ni c \succq b} q^{n_b(c)} X_c\!\!\mod\! \Si_+(b), 
\!\where n_b(c)\in\Z_{-}.
\end{align}

(ii) The coefficients of $\overline{E}_b^\dag$ belong to $\Z[q^{-1}]$
for any $b\in B$.
\end{corollary}

\proof
($i$) The proof is straightforward using (\ref{ebminusdag})
and (\ref{ebarmodsi}).

($ii$) By (\ref{epolexists}), the denominators of the coefficients
in $E_b$ are of products of factors of the form
$$ (1-q^j\prod_\nu t_\nu^{m_\nu}), \where j,\sum_\nu m_\nu> 0. $$
By Proposition~\ref{EDAG}($iii$), we already know that 
the $\overline{E}_b^\dag$ exist, at least with coefficients
in $\Q_q'$. Hence we may set $t=t_\nu$ for all $\nu$ when
calculating the limits of the coefficients.
As polynomials in $t$, the denominators of $E_b$ then
have leading terms of the form $\pm q^r t^s$ where $r,s>0$.
Since $\overline{E}_b^\dag$ exists, no higher power of $t$ can
appear in the corresponding numerator.
Therefore, we see that the coefficients of $\overline{E}_b^\dag$
belong to $\Z[q^{\pm 1}]$.
Using ($i$) and then (\ref{daginteb}) inductively, now it is easy
to see that the coefficients of $\overline{E}_b^\dag$ lie in
$\Z[q^{-1}]$.
\sq
\medskip

{\em Positivity conjecture.}
We conjecture that the coefficients of 
$\overline{E}_b^\dag$ for all $b\in B$ belong to $\Z_+[q^{-1}]$,
which will hopefully
follow from a more systematic theory of the intertwiners
in the nil-case. 
\smallskip

For $b\in B_-$, we expect
that the polynomials $\overline{E}^\dag_b$ and their untwisted
counterparts (not considered in this work) coincide with  
the corresponding level-one Demazure characters 
for the Kac-Moody algebra $\hat{\mathfrak{g}}$ associated
with $\tR$ upon shifting the $q$\~powers by the 
PBW-degrees from \cite{FFL} (twisted or untwisted). 
We thank Evgeny Feigin for his help with settling the 
conjecture below, which is directly related to our ongoing 
joint research project with him. 

It would give the positivity of the coefficients
of $\overline{E}^\dag_b$ for $b\in B_-$. Moreover, 
then (\ref{daginteb}) can be generally applied to verify that
the coefficients of $\overline{E}_b^\dag$ are  
from $\Z_+[q^{-1}]$ for {\em all\,} $b\in B$,
though their ``geometric" meaning is unclear to us.
We suspect here a connection with the {\em local Weyl modules} 
considered under the PBW\~filtration.

Since the relations of $\overline{E}_b$ 
to the level-one Demazure characters holds only in the twisted
case \cite{Ion1} and simply
to avoid giving the definitions of the untwisted dag-polynomials 
and the twisted PBW-filtration, we state the conjecture in this 
paper only for $\hat{\mathfrak{g}}$ of $AD\!E$\~type. Let  
$\hat{\mathfrak{b}}_+\supset \hat{\mathfrak{n}}_+$ 
be the Borel subalgebra and its radical. 

\begin{conjecture}\label{CONJDOMINB}
For $b\in B_-$,\, let $\,\v_{-b}$ be the Demazure module in 
the so-called basic\, 
representation of $\,\hat{\mathfrak{g}}$, which is 
a $\,\hat{\mathfrak{b}}_+$\~module generated by 
the extremal vector $v_{-b}$ of weight $-b\in B_+$.
In the setting from \cite{San,Ion1},
\begin{align*}
\overline{E}_b\ =\ \sum_{c\succq b,\,g} \hbox{dim}_{-c,g}\, q^g\,X_c
\for c\in B\,,
\end{align*}
where
$\,\hbox{dim}_{-c,g}\,=\,\dim V_{-c,g}\,$ for the subspace 
$V_{-c,g}\,$ of the vectors of degree $g\in \Z_+$ 
in the subspace  $V_{-c}\subset \v_{-b}\,$ of weight $-c\,$ for 
$c\in B$ and the standard Kac-Moody grading. 
Then 
\begin{align*}
&\overline{E}_b^\dag =\!\!
\sum_{B\ni c\succq b,\,f\ge 0}\!\!\!\!\dim 
\bigl(\g_f\,\cap V_{-c,g}\,/\,
\g_{f-1}\,\cap V_{-c,g}\bigr)\, q^{-f-g}\,X_c, \hbox{\ where}\\
&\g_{-1}\,=\,\{0\},\ \g_0\,=\,\C v_{-b},\ \, 
\g_f\ =\ \hat{\mathfrak{n}}_+(\g_{f-1})+\g_{f-1}\for f>0.
\end{align*}
In particular, $-n_b(w(b))$ for $w\in W$ defined in (\ref{nbce})
equals the minimal number of $\be\in R_+$ such that
$w(b)-b=\sum \be$ in the case of $A_n$. 

\vskip -0.3cm
\sq
\end{conjecture}
\medskip

{\em Connection maps.}
Let us introduce $\dHH^{\flat,\dag}$ as the subalgebra of $\HH^\flat$
generated over $\ddot{\Q}_{q,t}^\dag$ by the elements
\begin{align}
\label{dothhdag-gens}
X_b \ (b\in B),\ \ddot{T}_i^\dag \ (i\ge 0),\ \Pi^\flat.
\end{align}
Cf. Definition \ref{DEFHHDOT}.
The defining relations of $\dHH^{\flat,\dag}$ in terms of these
generators are:
\smallskip

($o$)\ \  $ (\ddot{T}_i^\dag-1)(\ddot{T}_i^\dag+t_i^{-1})\ =\
0,\ 0\ \le\ i\ \le\ n$;

($i$)\ \ \ $ \ddot{T}_i^\dag\ddot{T}_j^\dag\ddot{T}_i^\dag\cdots\ =\
\ddot{T}_j^\dag\ddot{T}_i^\dag\ddot{T}_j^\dag\cdots,\ m_{ij}$
factors on each side;

($ii$)\ \   $ \pi_r\ddot{T}_i^\dag\pi_r^{-1}\ =\ \ddot{T}_j^\dag \iif
\pi_r(\al_i)=\al_j$,\, $\pi_r\in \Pi^\flat$;

($iii$)\  $\ddot{T}_i^\dag X_b \ =\ X_b X_{\al_i}^{-1} 
(\ddot{T}_i^\dag)' \iif
(b,\al^\vee_i)=1,\ 0 \le i\le  n$;

($iv$)\ $\ddot{T}_i^\dag X_b\ =\ X_b \ddot{T}_i^\dag$ 
if $(b,\al^\vee_i)=0
\for 0 \le i\le  n$;

($v$)\ \ $\pi_rX_b \pi_r^{-1}\ =\ X_{\pi_r(b)}\ =\
X_{ u^{-1}_r(b)} q^{(\om_{r^*},b)}$.
\smallskip

Let $\bHH^{\flat,\dag}$ be the specialization of $\dHH^{\flat,\dag}$
as all $t_\nu=\infty$ (i.e., $t_\nu^{-1}=0$), and $\overline{\v}^\dag$
the image of the polynomial representation under this specialization.
The involution $\eta$ from (\ref{eta}) gives the 
following isomorphisms:
\begin{align}
\eta : \ &\bHH^\flat \to \bHH^{\flat,\dag}\label{etabar}\\
&\overline{T}_i\mapsto (\overline{T}_i^\dag)',\
X_b \mapsto X_b^{-1},\
\pi_r \mapsto \pi_r,\
q \mapsto q^{-1}\notag\\
\eta^\dag\equal\eta^{-1} : \ 
&\bHH^{\flat,\dag}\mapsto\bHH^\flat\label{etadag}\\
&\overline{T}_i^\dag \mapsto \overline{T}_i',\
X_b \mapsto X_b^{-1},\
\pi_r \mapsto \pi_r,\
q \mapsto q^{-1}.\notag
\end{align}
Due to (\ref{etastar}), one has
\begin{align}
&(H(f))^\ast=\eta(H)(f^\ast),
\for f\in\overline{\v},\ H\in\bHH^\flat\\
&(H(f))^\ast=\eta^\dag(H)(f^\ast),
\for f\in\overline{\v}^\dag,\ H\in\bHH^{\flat,\dag}.
\end{align}
\medskip

\setcounter{equation}{0}
\section{\sc Nonsymmetric Whittaker function}
Let us recall the definition of the Ruijsenaars-Etingof 
limiting procedure from \cite{Ru,Et1}
employed and developed in
\cite{C10} (the symmetric theory). 
Given a difference operator $\cL$
and a function $F$, it is defined by
\begin{align}\label{kapprho}
&\kapp(\cL) = (X_{\rho_k}\Ga_{-\rho_k})\,
\cL\,(X_{\rho_k}\Ga_{-\rho_k})^{-1},\
\kapp(F) =
X_{\rho_k}\Ga_{-\rho_k}(F)\\
&\RE(A) = \lim_{k\rightarrow\infty} \kapp(A),\ \
\RE(F) = \lim_{k\rightarrow\infty} \kapp(F).\notag
\end{align}
This procedure was applied in \cite{C10} to obtain
global Whittaker functions from global {\em symmetric}
$q,t$\~spherical functions in the symmetric case.
The existence a nonsymmetric analogue of this procedure
remained an entirely open question until \cite{ChM}, where it was shown
for the root system $A_1$ as an application of $W$\~spinors.
In \cite{ChO2}, a systematic algebraic study
the corresponding nonsymmetric (spinor) Whittaker functions
(still for $A_1$ only) was carried out; this involved a detailed
analysis of certain subalgebras of the nil-DAHA, most importantly
the {\em core subalgebra}.

In this section, we develop a nonsymmetric (spinor) variant of the
Ruijsenaars-Etingof procedure and apply it to global
{\em nonsymmetric} $q,t$\~spherical function.
Algebraically, these constructions are closely related with
the theory of {\em pseudo-polynomial representation} of nil-DAHA,
which is introduced at the end of this section
and will be the subject of our future works.
\medskip

\subsection{\bf Global spherical functions}
\comment{
One of the  main advantages of the technique of Gaussians is
a possibility to introduce the spherical function as
a reproducing kernel of the Fourier transform from
$\v\ga^{-1}$, the polynomial representation multiplied
by the Gaussian $\ga^{-1}$,
to the $\HH^\flat\,$\~module  $\v\ga$.
In this setting, the construction below is directly
related to the identities (\ref{pbgauss})
(correspondingly, (\ref{pbgaussl}) in the limit).}

By the {\em Gaussians} we mean
\begin{align}
&\tga^{\oplus}=\sum_{b\in B} q^{-(b,b)/2}X_b,\
\tga^{\ominus}=\sum_{b\in B} q^{(b,b)/2}X_b.
\label{gauser}
\end{align}
We need mainly the Gaussian $\tga^\ominus$ in this paper. 
The multiplication by $\tga^{\ominus\,}$ preserves
the space of Laurent series
with coefficients in $\Q[t][[q^{\frac{1}{2\tilde{m}}}]]$, where
$\tilde{m}$ is from the definition of $\Q_{q,t}'.$ 

The Gaussians satisfies the fundamental difference equations
\begin{align}\label{gauseqn}
\Ga_a(\tga^{\oplus})
&=q^{(a,a)/2}X_a\tga^{\oplus},\
\Ga_a(\tga^{\ominus})
=q^{-(a,a)/2}X_a^{-1}\tga^{\ominus}\for a\in B.
\end{align}
Later we will need the following special case of (\ref{gauseqn}):
\begin{align}\label{gausrho}
\Ga_{\rho_k}(\tga^{\ominus})
&=q^{-(\rho_k,\rho_k)/2}X_{-\rho_k}\tga^{\ominus},
\end{align}
provided $\rho_k\in B$ (when $B=P$, the condition 
$k_\nu\in\Z$ is sufficient).

If one uses here the {\em real Gaussians} defined as
\begin{align}\label{gaussreal}
&\ga^{\pm 1}=q^{\pm x^2/2}, \where
X_b\equal q^{x_b}, x_b=(x,b),
x^2=\sum_i\, x_{\al_i}\,x_{\om_i^\vee},
\end{align}
then (\ref{gauseqn}) is satisfied for any complex $a$.
Note that if the series $\tga^{\ominus}$ is considered
as a holomorphic function for $|q|<1$, then 
the function $\tga^{\ominus\,}\ga$ is
$B$\~periodic in terms of $x$.

\comment{
The {\em $q$\~Mehta\~Macdonald identity} from \cite{C5}
\begin{align}
&\langle \tga^{\ominus\,}\mu_\circ\rangle\ =\
\prod_{\al\in R_+}\prod_{ j=1}^{\infty}\Bigl(\frac{
1-t_\al^{-1} q_\al^{(\rho_k,\al^\vee)+j}}{
1-           q_\al^{(\rho_k,\al^\vee)+j} }\Bigr)
\label{mehtamu}
\end{align}
provides the normalization constant
for the $q$\~Gauss integrals
\begin{align}\label{pbgauss}
\lan \ep_b(X) \ep_c^\star(X) \tga^{\ominus\,}\mu_\circ \ran
&=q^{\frac{(b_\#,b_\#)+(c_\#,c_\#)}{2}-(\rho_k,\rho_k)}
\ep_c^\star(q^{b_\#})\lan \tga^{\ominus\,}\mu_\circ\ran,
\end{align}
where $b,c\in B$. It implies the duality formula  (\ref{pebdual}).
There are counterparts of (\ref{pbgauss}) for $\tga^{\,\oplus}$
(treated as an analytic function for $|q|>1$),
and for the Jackson summation taken instead of the constant term
functional. See \cite{C5,C101}.
The results from this paper
can be readily extended to these cases, but we will
stick to the case $|q|<1$.
{\bf Taking the limit.}
Let us tend $t\to 0$ in (\ref{pbgauss}). The
definition of the $P$\~polynomials implies that
\begin{align}\label{limpct}
&\lim_{t\to 0} q^{-(c,\,\rho_k)}
P_c(q^{z-\rho_k})\ =\ q^{(c_+\,,\,z)} \for c\in B_-\,, c_+=w_o(c).
\end{align}
Note that $-(c,\rho_k)=(c_+,\rho_k)$.
For instance, it matches the evaluation
formula in (\ref{pebebs}):
$\lim_{t\to 0} q^{-(c,\rho_k)}
P_c(q^{-\rho_k})=1$.
We come to the following formulas ($c\in B_-\,$):
\begin{align}
&\langle \tga^{\ominus\,}\overline{\mu}_\circ\rangle\ =\
\prod_{i=1}^n\prod_{ j=1}^{\infty}(1-q_i^j),
\label{mehtamul}
\end{align}
\begin{align}\label{pbgaussl}
&\langle \overline{P}_b(X) \overline{P}_c(X)
\tga^{\ominus\,}\overline{\mu}_\circ \rangle\  =\
 q^{\frac{(b,b)+(c,c)}{2}}\, X_{c_+}(q^{b})\,
\langle \tga^{\ominus\,}\overline{\mu}_\circ\rangle.
\end{align}
Here $X_{c_+}(q^{b})=q^{(c_+,b)}=q^{(c,b_+)}=X_{b_+}(q^{c})$.
\medskip
}

We assume that $|q|<1$ and
use the notation $\tga_\la^{\ominus}$ 
for the Gaussians defined for the variable
$\La=q^\la$. Let $\tga^{\ominus}_x=\tga^\ominus$ for
the sake of uniformity.
Thus, $\tga_\la^{\ominus}=\tga^\ominus(q^\la)$.
Accordingly, we use superscripts when applying operators
from the polynomial representation of $\HH$ to functions
of $X$ or $\La$. For instance, we write $T_i^\la$ for the
action of $T_i$ from (\ref{Demazx}) on functions of $\La=q^\la$,
where we replace $X_{\al_i}$ by $\La_{\al_i}$.
When no superscript is used, the action is understood in
terms of $X$.
\smallskip

We will also use the normalization constant
\begin{align}
&\tga^{\ominus}(q^{\rho_k})\ =\ 
\sum_{a\in B} q^{\frac{(a,a)}{2}+(\rho_k,a)}.
\label{lgarhok}
\end{align}

The following theorem results from Theorem 5.4 and 
Corollary 7.3 of \cite{C5}. The function $G(X,\La)$
introduced in (\ref{gexla}) is called 
{\em global nonsymmetric $q,t$\~spherical function}.
  
\begin{theorem}\label{GLOBNSSPH}
(i) The Laurent series
\begin{align}
\Xi(X,\La;q,t)&\equal
\sum_{b\in B} q^{(b_\#,b_\#)/2-(\rho_k,\rho_k)/2} \
\frac{E_b^\star(X) \, E_b(\La)}
{\lan E_b, E_b \ran_\circ}
\label{psiexla}
\end{align}
is well defined with coefficients in
$\Q[t][[q^{\frac{1}{2\tilde{m}}}]]$.
For $|q|<1$, $\Xi$ converges to an entire function
of $X,\La$, provided $t_\nu$ are chosen so that all $E$\~polynomials
exist (by (\ref{epolexists}),
the conditions $|t_\nu|<1$ are sufficient).
Accordingly, $G(X,\La)$ defined via
\begin{align}
\frac{\tga_x^{\ominus}\tga_\la^{\ominus}}
{\tga^{\ominus}(q^{\rho_k})}
G(X,\La) \equal \Xi(X,\La;q,t)
\label{gexla}
\end{align}
is a meromorphic function of
$X,\La$, which is analytic
apart from the zeros of $\tga_x^{\ominus\,}\tga_\la^{\ominus\,}$.


(ii) The function $G(X,\La)$ satisfies
\begin{align}\label{GsymT}
G(X,\La)=G(\La,X),\ T_i^x(G(X,\La)=T_i^\la(G(X,\La),\
1\le i\le n,
\end{align}
and the following extension of (\ref{Yone}):
\begin{align}
Y_a(G(X,\La))=\La_{a}^{-1}G(X,\La) \for a\in B.
\label{LfLa}
\end{align}
For an arbitrary $b\in B$, one has
\begin{align}\label{Shinqt}
& G(X, q^{b_\#})\ =\ 
\frac{E_b(X)}{E_b(q^{-\rho_k})} 
\prod_{\al\in R_+}\prod_{ j=1}^{\infty}\Bigl(\frac{ 
1- q^{(\rho_k,\al)+\nu_\al j}}{
1-t_\al^{-1}q^{(\rho_k,\al)+\nu_\al j} }\Bigr). 
\end{align}
\sq
\end{theorem}

Note that relations (\ref{GsymT}) and (\ref{LfLa})
can be uniformly presented using the anti-involution $\vph$
from (\ref{phianti}) as follows:
\begin{align}
H^x(G(X,\La))\ =\ (\vph(H))^\la (G(X,\La)) \for H\in \HH^\flat.
\label{LfLavph}
\end{align}
It reflects the fundamental fact that $G(X,\La)$ represents the
Fourier transform of DAHA corresponding to
the automorphism $\si$ from (\ref{si}),
which satisfies the relation 
$\vph\si\vph=\si^{-1}$; see \cite{C5}.
\smallskip

\subsection{\bf Action of intertwiners}
Relation (\ref{LfLavph}) results in the following formulas
for the action of the intertwining operators on $E_b^\star$.
Let us first modify the intertwiners $\Psi_i$ and $\tau_+(\pi_r)$
from (\ref{Phijb}) and Theorem \ref{PHIEB}
by applying the automorphism $\tau_-^{-1}$, 
which preserves the elements $Y_b,T_i$ for any $b\in B,i\ge 0$. 
Namely, we set $Y_0=Y_{\al_0}=q^{-1}Y_\vth^{-1}$,
\begin{align}
&\Psi'_{s_i}=\Psi'_i\ =\ 
\tau_-^{-1}\tau_+(T_i) + \frac{t_i^{1/2}-t_i^{-1/2}}
{Y_{\al_i}^{-1}-1},\ \ 0\le i\ge n,\label{Psiiy}\\
&\Psi_{\pi_r}=\,\tau_-^{-1}\tau_+(\pi_r) \for  r\in O,
\ \Psi'_{\hw}=\Psi'_{\pi_r}\Psi'_{i_l}\cdots\Psi'_{i_1},\ 
\notag
\end{align}
where the decomposition 
$\hW^\flat\ni \hw=\pi_r s_{i_l}\cdots s_1$ is
reduced and $\Psi'_{\hw}$ does not depend on its choice.
Then $\Psi'_{\hw}(E_b)$ is
proportional to $E_{c}$ for $c=\hw(\!(b)\!)$
and any $b\in B$. More precisely,
using the action $\tau_-^{-1}(E_b)=q^{b_\#^2/2-\rho_k^2/2}E_b$
\,($b\in B$)\, in the polynomial representation from
\cite{C101},Proposition 3.3.4, 
\begin{align}
t_i^{1/2}\Psi'_{i}\bigl(&q^{(b_+,\rho_k)}E_b\bigr)=  
q^{(c_+,\rho_k)} E_c, \ \,
t_i^{-1/2}\Psi'_{i}\bigl(\frac{q^{(c_+,\rho_k)}
E_c}{\lan E_c,E_c\ran_\circ}\bigr)=  
\frac{q^{(b_+,\rho_k)}E_b}{\lan E_b,E_b\ran_\circ},\notag\\
\label{PsiPrime}
&\hbox{for\ }\, 0\le i\le n,\ c=s_i(\!(b)\!)
\hbox{\ \ \,provided that \ } (\al_i,b+d)>0,\\
&\Psi'_{\pi_r}\bigl(q^{(b_+,\rho_k)}E_b\bigr) =  
q^{(c_+,\rho_k)} E_c \for c=\pi_r(\!(b)\!) \and
r\in O'.\notag
\end{align}

Second,  we set
$\tilde{\Psi}'_{\hw}=\tilde{\Psi}'_{\pi_r}\tilde{\Psi}'_{i_l}
\cdots\tilde{\Psi}'_{i_1}$ for induced decompositions
$\hW^\flat\ni \hw=\pi_r s_{i_l}\cdots s_{i_1}$, where $0\le i\le n,\
r\in O\,$ and
\begin{align}\label{Psiiyy}
&\tilde{\Psi}'_{s_i}=\tilde{\Psi}'_i= 
\,\tau_+^{-1}(T_i) + \frac{t_i^{1/2}-t_i^{-1/2}}
{\si^{-1}(X_{\al_i})-1},\ 
\tilde{\Psi}'_{\pi_r}=\,\tau_+^{-1}(\pi_r).
\end{align}
Note that
\begin{align}\label{Psitilde}
&\tau_-(\si^{-1}(X_b))=\si^{-1}(\tau_+^{-1}((X_b))=
\si^{-1}(X_b), \hbox{\ which results in }\\
&\tau_-(\tilde{\Psi}'_{\hw})\ =\ \tilde{\Psi}''_{\hw}
\equal\tilde{\Psi}''_{\pi_r}\tilde{\Psi}''_{i_l}\cdots
\tilde{\Psi}''_{i_1}
\hbox{\ \,for\,\ } 
\tilde{\Psi}''_{\pi_r}=\si^{-1}(\pi_r),\, r\in O,\notag\\
\tilde{\Psi}''_i\,&=\,\tau_-\tau_+^{-1}\tau_-(T_i) + 
\frac{t_i^{1/2}-t_i^{-1/2}}{\si^{-1}(X_{\al_i})-1}\,=\,
\si^{-1}(T_i) + 
\frac{t_i^{1/2}-t_i^{-1/2}}{\si^{-1}(X_{\al_i})-1}.\notag
\end{align} 
Therefore both families here do not depend on particular
choices of the reduced decompositions of $\hw\in \hW^\flat$ and
intertwine $\{\si^{-1}(X_b)\}$\,:
\begin{align}\label{Psiintert}
\tilde{\Psi}''_{\hw}\ \si^{-1}(X_b)\,
(\tilde{\Psi}''_{\hw})^{-1}=
 \si^{-1}(X_{\hw(b)})=
\tilde{\Psi}'_{\hw}\ \si^{-1}(X_b)\,
(\tilde{\Psi}'_{\hw})^{-1},
\end{align}
where $b\in B,\, \hw\in \hW^\flat$.

\begin{proposition}\label{PSIESTAR}  
For $\,\Xi(X,\La;q,t)\,$ from (\ref{psiexla}) and $\,\hw\in \hW^\flat$,
\begin{align}\label{Psivph}
&(\tau_+^{-1}\vph\tau_+)\bigl(\Psi'_{\hw}\bigr)= 
\tilde{\Psi}'_{\hw^{-1}},\ 
(\Psi'_{\hw})^{\la\,}(\Xi(X,\La))\!=\! 
(\tilde{\Psi}'_{\hw^{-1}})^{x\,}(\Xi(X,\La)).
\end{align}
Accordingly, $\tilde{\Psi}'_{\hw}(E_b^\star)$ is 
proportional to $E_{c}$ for $c=\hw(\!(b)\!)$
and any $b\in B$. More precisely,
\begin{align}
&t_i^{\frac{1}{2}}\tilde{\Psi}'_{i}\bigl(\frac{q^{\frac{(b,b)}{2}}
E_b^\star}{\lan E_b,E_b\ran_\circ}\bigr)\!=\!  
\frac{q^{\frac{(c,c)}{2}}E_c^\star}{\lan E_c,E_c\ran_\circ},\ 
t_i^{-\frac{1}{2}}\tilde{\Psi}'_{i}\bigl(q^{\frac{(c,c)}{2}}
E_c^\star\bigr)=q^{\frac{(b,b)}{2}}E_b^\star\label{PsiEstar}\\ 
&\hbox{for}\ 0\le i\le n,\
b\in B,\ \ c=s_i(\!(b)\!)
\hbox{\ when \ } (\al_i,b+d)<0,\notag\\
&\hbox{and\ \,}\tau_+^{-1}(T_i)(E_b^\star) = 
t_i^{1/2}E_b^\star \hbox{\  when \ }
(\al_i,b+d)=0, \,
i\ge 0\notag,\\  
&\tilde\Psi'_{\pi_r}\bigl(\frac{q^{\frac{(b,b)}{2}}
E_b^\star}
{\lan E_b,E_b\ran_\circ}\bigr)\ \,=\ \, 
\frac{q^{\frac{(c,c)}{2}}E_c^\star}{\lan E_c,E_c\ran_\circ}
\for c=\pi_r(\!(b)\!),\, r\in O'.
\notag
\end{align}
\end{proposition}
{\it Proof.}
To check (\ref{Psivph}), we use that $\tau_+(T_i)=T_i$
for $i>0$ and that
\begin{align}\label{phifixed}
&\tau_+(\pi_r)=\,q^{-\om_r^2/2}X_{\om_r}\pi_r=\,
q^{-\om_r^2/2}X_{\om_r}T_{u_{r^*}}Y_{\om_{r^*}}^{-1}\ (r\in O),
\notag\\
&\tau_+(T_0)=X_0^{-1}T_0^{-1}=\,q^{-1}X_{\vth}
T_{s_{\vth}}Y_{\vth}^{-1},\hbox{\ which result in\ }\notag\\
&\vph\bigr(\tau_+(\pi_r)\bigl)=
\tau_+(\pi_{r^*}) \for r\in O,\ \,
\vph\bigr(\tau_+(T_0)\bigl)=\tau_+(T_0).  
\end{align}
Use (\ref{taux}) and the relation $\pi_r=\pi_{r^*}^{-1}$.
Accordingly,
\begin{align*}
&\tau_+^{-1}\vph\tau_+\bigr(\tau_-^{-1}\tau_+(\pi_r)\bigl)\ =\ 
\tau_+^{-1}\vph\tau_+\bigr(\tau_-^{-1}\vph\tau_+(\pi_{r^*})
\bigl)\\
=&(\tau_+^{-1}\tau_-\tau_+^{-1}\tau_+)(\pi_{r^*})=
(\tau_+^{-1}\tau_-)(\pi_{r^*})=\tau_+^{-1}(\pi_{r^*}),\\
&\hbox{and\ \ \,}\tau_+^{-1}\vph\tau_+\bigr(\tau_-^{-1}
\tau_+(T_0)\bigl)\ =\,\cdots\,=\ 
\tau_+^{-1}(T_0).
\end{align*}
Also, 
\begin{align*}
\tau_+^{-1}\vph\tau_+(Y_b)=
\tau_+^{-1}\tau_-(X_b^{-1})=\tau_+^{-1}\tau_-\tau_+^{-1}(X_b^{-1})=
\si^{-1}(X_b^{-1}).
\end{align*}

This gives the first relation from (\ref{Psivph});
the second follow directly from (\ref{LfLavph}) and
the definition of the series $\,\Xi(X,\La;q,t).$
The remaining formulas result from (\ref{PsiPrime})
and the structure of this series.\sq
\smallskip

We note that one can prove 
this proposition {\em directly} using the
relation
\begin{align*}
\si^{-1}(X_b)=\si^{-2}\bigl(\si(X_b)\bigr)=
\si^{-2}\bigl(Y_b^{-1}\bigr)=
T_{w_0}Y_{w_0(b)}T_{w_0}^{-1}
\end{align*}
from Proposition 3.2.2 from \cite{C101}.
Combined with  (\ref{ebast}), which states that
\begin{align*}
&E_b^\star =
C_b\, T_{w_0}(E_{\varsigma(b)}) \hbox{\ \,for\ \, }
C_b=\!\!\prod_{\nu\in \nu_R} t_\nu ^{l_\nu(u_b)- 
\frac{l_\nu(w_0)}{2}}\!=
C_{\varsigma(b)}\,,
\varsigma(b)=-w_0(b),
\end{align*}
we obtain that
\begin{align}\label{siXEb}
&\si^{-1}(X_a)\bigl(E_b^\star\bigr)\ =\ 
\,q^{(a_{\#},b)}
\bigl(E_b^\star\bigr) \for a,b\in B.
\end{align}
Indeed,
\begin{align*}
C_b^{-1}\,\si^{-1}(X_a)\bigl(E_b^\star\bigr)=
&T_{w_0}Y_{w_0(a)}T_{w_0}^{-1}\bigr(T_{w_0}(E_{\varsigma(b)})\bigl)=
T_{w_0}\bigr(Y_{w_0(a)}(E_{\varsigma(b)})\bigl)\\
=&C_{\varsigma(b)}^{-1}\,q^{-(w_0(a),\varsigma(b)_{\#})}
\bigl(E_b^\star\bigr)=
C_{b}^{-1}\,q^{(a,b_{\#})}
\bigl(E_b^\star\bigr).
\end{align*}

This provides a direct approach to the justification
the symmetries (\ref{LfLavph}), though related
to that based on the interpretation of $G(X,\La)$ as 
the reproducing kernel of the DAHA-Fourier transform.
\smallskip

The symmetry (\ref{LfLavph}) coupled with
(\ref{phifixed}) becomes even simpler for
the standard $\Psi$\~intertwiners (without
the conjugation by $\tau_+^{-1}$). Following
(\ref{normPsiinter}), let 
$\mathfrak{G}'_i=\Psi'_i\psi_i^{-1}\ (i\ge 0)$ for 
$\psi_i =t_i^{1/2}+
(t_i^{1/2}-t_i^{-1/2})/(Y_{\al_i}^{-1}-1)$.
We also set $\,\mathfrak{G}'_{\pi_r}=\Psi'_{\pi_r}\ (r\in O)$
and define \,$\mathfrak{G}'_{\hw}$ \,accordingly. 

Then for
$\,\hw=aw\in \hW^\flat$, where $a\in B,\, w\in W$, 
\begin{align}\label{Psi-straight}
&(\mathfrak{G}'_{\hw^{-1}})^{\la\,}
\bigr(G(X,\La)\bigl)\ =\ 
\hw^{x\,} \bigl(G(X,\La)\bigr)\\
=&\ q^{\frac{a^2}{2}}\,X_{-a}\,
G(E_b^\star(X)\!\mapsto\! \hw\bigl(E_b^\star(X)\bigr),\La).\notag
\end{align}
Indeed, $\vph\bigl(\tau_-^{-1}\tau_+(H)\bigr)=
\vph\bigl(\tau_-^{-1}\vph\tau_+(H)\bigr)=H\ $ for
$\ H=T_i \,(i\ge 0)$\, or for\, $H=\pi_r\, (r\in O)$. Therefore
\begin{align*}
&\vph(\Psi'_i)=T_i+\frac{t_i^{1/2}-t_i^{-1/2}}{X_{\al_i}-1}\!=\!
\frac{t_i^{1/2}X_{\al_i}-t^{-1/2}}{X_{\al_i}-1}s_i, \ 
\vph(\mathfrak{G}'_i)=s_i\\
&\hbox{for\ \ } 0\le i\le n \hbox{\ \ \ and\ \ }
\vph(\Psi'_{\pi_r})\ =\ \vph(\mathfrak{G}'_{\pi_r})\ =\ \pi_{r^*} 
\ \,(r\in O).
\end{align*} 
Formula (\ref{Psi-straight}) is compatible with the Shintani-type
formula (\ref{Shinqt}).
As a matter of fact, this provides a direct
way for establishing (\ref{Shinqt}). Thus the formulas for
the action of the $\Psi$\~intertwiners on the $E$\~polynomials
\underline{}are essentially sufficient for a direct verification of 
Theorem \ref{GLOBNSSPH}.   
\smallskip

\subsection{\bf \texorpdfstring{$W$}{\em W}\~spinors}
By {\em $W$\~spinors}, we simply mean maps $W\to\v$; we denote the
space of $W$\~spinors by $\spin(\v)$.
Thus $\spin(\v)$ is naturally a $\Q_{q,t}'$\~algebra
under pointwise multiplication and addition.
For any $w\in W$,
denote by $\ze_w$ the characteristic function 
$\ze_w(u)=\de_{wu}$. These are pairwise orthogonal
idempotents in $\spin(\v)$. Any element in $\spin(\v)$
can be written uniquely as
\begin{align}
f=\sum_{w\in W}f_w\,\ze_w, \where f_w\equal f(w)\in\v.
\end{align}
We refer to $f_w$ as the $w$\~component of $f$.

We equip $\spin(\v)$ with an action of $W$ via
\begin{align}
(\de(w)f)(u)\equal f(w^{-1}u).
\end{align}
Note that $\de(w)(\ze_v)=\ze_{wv}$
and hence for any $f\in\spin(\v)$,
\begin{align}
(\de(v)f)_w = f_{v^{-1}w}.
\label{vwf}
\end{align}

One has a natural embedding
\begin{align}\label{de-of-f}
\de(f)&\equal\sum_{w\in W} f\,\ze_w, \ f\in\v;
\end{align}
its image is the space of $W$\~invariants of $\spin(\v)$, which
will be denoted by $\spin^\de(\v)$. These $W$\~invariants will be
called {\em $\delta$\~spinors}. 

These definitions are completely independent of the natural action
of $W$ in $\v$ and can be applied to any spaces of functions
instead of $\v$.

Using the action of $W$ in $\v$, there is another embedding
$\varrho:\v\to\spin(\v)$ defined by
\begin{align}\label{rho-spin}
\varrho(f)&\equal\sum_{w\in W} w^{-1}(f)\,\ze_w.
\end{align}
The spinors in its image will be called {\em principal
spinors} or {\em $\varrho$\~spinors}, sometimes, 
simply {\em functions}.

For $W$\~invariant $f \in \v$, the spinors
$\varrho(f)$ and $\de(f)$ coincide;
thus we will simply write $f$ for $W$\~invariant functions.
For arbitrary $f\in\v$ we may also write $f^\varrho\equal\varrho(f)$.
\smallskip

{\em Spinor difference operators.}
Generally, any endomorphism of $\v$ acts pointwise
in $\spin(\v)$. For instance, for a translation
$\Ga_b$, we set
\begin{align}
\Ga_b(f)(u)\equal\Ga_b(f(u)).
\end{align}

We define
\begin{align}
\label{deGab}
&\de(\Ga_b)=\Ga_b^\de\equal\sum_{w\in W}\Ga_b\,\ze_w,\\
&\varrho(\Ga_b)=\Ga_b^\varrho\equal
\sum_{w\in W}\Ga_{w^{-1}(b)}\,\ze_w,
\label{rhoGab}
\end{align}
where $\ze_w$ acts by multiplication in $\spin(\v)$.

Similarly, we may view $X_b$ as a pointwise multiplication
operator:
\begin{align}
X_b(f)(u)\equal X_b(f(u)).
\end{align}

By a {\em spinor difference operator} we mean
any linear combination of operators of the form
$g(X)\,\Ga_b\,\ze_w \ (g(X)\in\Q_{q,t}'(X), b\in B, w\in W)$.
{\em Spinor difference-reflection operators} are defined
as linear combinations of operators of the form
$g(X)\,\Ga_b\,\ze_w\,\de(v)$, where $v\in W$.

Recall that $\cA$ denotes the algebra of difference-reflection
operators (see the beginning of Section \ref{SEC:POLREP}).
Replacing polynomials by rational functions in the definition
of $\spin(\v)$, we define an action of $\cA$ in 
$\spin(\Q_{q,t}'(X))$ by sending 
$$
\phi:\ g\,\Ga_b\,w\ \mapsto\
\varrho(g)\,\varrho(\Ga_b)\,\de(w),
\where
g\in\Q_{q,t}'(X),\ b\in B,\ w\in W.
$$
It is a homomorphism of algebras:
$\phi:\cA\rightarrow\mathrm{End}_{\Q_{q,t}'}(\spin(\Q_{q,t}'(X))).$

We obtain an action of $\HH$ in $\spin(\Q_{q,t}'(X))$
by composing $\phi$ with the polynomial representation
viewed as a homomorphism $\HH\to\cA$;
here $\HH$ acts by spinor difference-reflection operators.
Note that  $\HH$ does not preserve $\spin(\v)$.
\smallskip

\subsection{\bf Spinor \texorpdfstring{$R\!E$}{\em RE}-procedure}

Following \cite{ChM}, we set
\begin{align}\label{kappde}
\kapp^\de(\cL)&\equal \de(X_{\rho_k}\Ga_{-\rho_k})\,\notag
\phi(\cL)\,\de(X_{\rho_k}\Ga_{-\rho_k})^{-1},\\
\kapp^\de(F)&\equal
\de(X_{\rho_k}\Ga_{-\rho_k})\bigl(\varrho(F)\bigr),
\end{align}
for any function $F$ and any difference-reflection
operator $\cL$.
The nonsymmetric variant of (\ref{kapprho}),
{\em $\delta$\~Ruijsenaars-Etingof procedure},
is then defined as
\begin{align}\label{nonsymwhitdef}
\RE^\de(\cL) &\equal \lim_{k\rightarrow\infty} \kapp^\de(\cL),\ \
\RE^\de(F) \equal \lim_{k\rightarrow\infty} \kapp^\de(F),
\end{align}
where $k\to\infty$ means $\Re k_\nu\to\infty$ for all $\nu$
(equivalently, $t_\nu\to 0$).

We record here the following formulas for later use:
\begin{align}
\kapp^\de(v\cL w)&=v\,\kapp^\de(\cL)\,w
\fora v,w\in W,\label{RE-W-eq}\\
\kapp^\de(X_b)&=\sum_{w\in W}\prod_\nu 
t_\nu^{-(\rho_\nu^\vee,w^{-1}(b))}
X_{w^{-1}(b)}\,\ze_w,\label{kapp-de-X}\\
\kapp^\de(\Ga_b)&=\sum_{w\in W}\prod_\nu 
t_\nu^{-(\rho_\nu^\vee,w^{-1}(b))}
\Ga_{w^{-1}(b)}\,\ze_w.\label{kapp-de-Ga}
\end{align}
\medskip

{\em Spinor Whittaker function.}
Let $\Om(X,\La) \equal \RE_x^\de(G(X,\La))/\tga^{\ominus}(1)$, 
where the $\RE_x$ 
is applied to $X$ and we assume that $\rho_k\in B$;
for instance, $\Z\ni k_\nu\to\infty$ is sufficient for $B=P$.  
Equivalently, we can introduce
$$
G'(X,\La)=G(X,\La)\frac
{\tga^{\ominus}(X)\,q^{\frac{x^2}{2}}}
{\tga^{\ominus}(q^{\rho_k})q^{\frac{\rho_k^2}{2}}},\ \
\Om'(X,\La)=
\RE_x^\de(G'(X,\La)),
$$
where $\Re k_\nu\to\infty$ for {\em arbitrary complex} $k$.
Then 
$
\tga^{\ominus}(X)
q^{\frac{x^2}{2}}\Om(X,\La)=\Om'(X,\La).\ 
$ 
Using $G'$ instead of $G$ will not influence the corresponding
operators (studied below) acting on this function since
$\tga^{\ominus}(X)q^{\frac{x^2}{2}}$ is $\hW$\~invariant.
Recall that $|q|<1$, so $t_\nu\to 0$ if $\Re k_\nu\to \infty$.
Equivalently,
$$
\tga^\ominus(X)\,\Om(X,\La)\ =\ 
\lim_{\Re k_\nu\to \infty}
\Ga_{-\rho_k}\bigl(\tga^\ominus(X)
G(X,\La)/\tga^\ominus(q^{\rho_k})
\bigr).
$$

\begin{proposition}
\label{OMLIMIT}
The limit, as \,$\Re k_\nu\to\infty$ for all\, $\nu$,\,
of the series $\Ga_{-\rho_k}^\de(\Xi(X,\La;q,t))$ exists;
here $\Xi(X,\La;q,t)$ is the series from (\ref{psiexla}). 
Accordingly, 
one has
\begin{align}
\Om(X,\La)=(\tga_x^{\ominus}\tga_\la^{\ominus})^{-1}\,
\sum_{b\in B} q^{(b,b)/2}\,
\frac{\overline{E}_b(\La)}
{\lan \overline{E}_b,\overline{E}_b \ran_\circ}\,
\sum_{w\in W} a_{b,w}\,X_{-b_-}\,\ze_w,
\label{omform}
\end{align}
where $a_{b,w}$ is the limit, as all $\,t_\nu\to 0$,
of the coefficient of $X_{-w(b_-)}$ in $E_b^\star$.
In particular, $a_{b,w}\in\Z[q]$ and
$a_{b,u_b^{-1}}=1$, $a_{b,\id}=\de_{b,b_-}$.
\end{proposition}

\proof
First, one has
$$\de(X_{\rho_k}\Ga_{-\rho_k})\,
(\tga_x^{\ominus})^{-1}
=q^{(\rho_k,\rho_k)/2}\,(\tga_x^{\ominus})^{-1}\,
\de(\Ga_{-\rho_k}),$$
as operators acting on spinor-functions of $X$
(due to the $W$\~invariance of $\tga^\ominus$, we omit
$\varrho$ here).
Hence it suffices to consider the limit of
$$
q^{-(b_-,\rho_k)}\,\de(\Ga_{-\rho_k})\cdot
\varrho(E_b^\star(X)),
$$
or equivalently, for each $w\in W$, the limit of
\begin{align}
\label{limcomp}
&q^{-(b_-,\rho_k)}\,\Ga_{-\rho_k}(w^{-1}(E_b^\star))\,\ze_w.
\end{align}
Using (\ref{macd}), the limit of (\ref{limcomp})
as $t_\nu\to 0$ clearly exists and is given as 
$a_{b,w}\,X_{-b_-}\,\ze_w$ for $a_{b,w}$ as claimed.
By (\ref{ebastbar}),
$\overline{E_b^\star}=(\overline{E}_b^\dag)^\ast$.
Hence Corollary \ref{EDAGCOEFF} implies that $a_{b,w}\in\Z[q]$.
\sq
\smallskip

This proposition can be used to justify the existence of
the $\RE$\~limits of the Dunkl operators, which we will call
{\em Toda-Dunkl operators}. Moreover,
We arrive at the following Whittaker counterpart of
Parts $(ii)$ and $(iii)$ of Theorem \ref{GLOBNSSPH}
and formula (\ref{LfLavph}).

\subsection{\bf Main theorem}
\begin{theorem}\label{GLOBNSWH}
(i) The operators
$\RE^\de(H^\vph)$ are well defined for $H\in \dHH^\flat$,
where $H^\vph=\vph(H)$; their coefficients are from
$\Q_q'[X_b,b\in B]$, so they preserve 
$\spin(\overline{\v})\equal \spin(\Q_q'[X_b,b\in B])$.
For instance, the following operators
are well defined:
\begin{align}\label{hatoprs}
\hat{Y}_b\equal\RE^\de(Y_b),\  
\hat{X}_b\equal\RE^\de(\tilde{X}_b) \hbox{ for } \tilde{X}_b\equal 
\ddot{Y}_{-b}^\vph=t^{(b_+,\rho^\vee)}X_b,&\\
\hat{T}_i\equal\RE^\de(\ddot{T}_i)
\hbox{ for } i>0,\ 
\hat{T}_0\equal \RE^\de(\ddot{T}_0^\vph),\
\hat{\pi}_r\equal\RE^\de(\pi_r^\vph) \hbox{ for } r\in O'.&  \notag
\end{align}
  
(ii) The function $\Om(X,\La)$ has the following symmetries: 
\begin{align}\label{GsymTW}
&\hat{T}_i(\Om(X,\La))=T_i^\la(\Om(X,\La)),\
\hat{\pi}_r(\Om(X,\La))=\pi_r^\la(\Om(X,\La))\\ 
&\hbox{for \ \,}
0\le i\le n, \,r\in O',\ 
T_i^\la=T_i\mid_{X\mapsto \La},\ 
\pi_r^\la=\pi_r\mid_{X\mapsto \La}. \notag
\end{align}
It satisfies the following limiting version of the relations from 
(\ref{LfLa}):
\begin{align}
\hat{Y}_{b}(\Om(X,\La))=\La_{b}^{-1}\Om(X,\La),\ 
\overline{Y}^\la_{b}(\Om(X,\La))=\hat{X}_{-b}\,\Om(X,\La),
\label{LfLaW}
\end{align}
where $\overline{Y}_b^\la=\overline{Y}_b\mid_{X\mapsto \La}$ 
for $b\in B$.

(iii) For an arbitrary $c\in B$, let 
$f(q^{c_\dag})\equal f_{u^{-1}}(q^{c_-})$
for a spinor $f=\sum_{w\in W} f_w\,\ze_w$
and for $u=u_c\in W$ from Proposition \ref{PIOM}
such that $u(c)=c_-\in B_-$ and \,$l(u)$ is minimal
possible.
Then
\begin{align}\label{ShinqtW}
& \tga^{\ominus}(1)\Om(q^{c_\dag},\La)\, =\, 
\overline{E}_c(\La)
\prod_{i=1}^n\prod_{ j=1}^{\infty}\frac{1}{1-q_i^{j}},\ \,
\tga^{\ominus}(1)=\sum_{b\in B}q^{b^2/2}. 
\end{align}
Equivalently, see (\ref{omform}) and (\ref{ebarnorms}),
\begin{align}\label{ShinqtWW}
&\sum_{b\in B} q^{(b_--c_-)^2/2}\,
\frac{\overline{E}_b(\La)}
{\lan \overline{E}_b,\overline{E}_b \ran_\circ}\,
a_{b,u_c^{-1}}
\,=\, 
\tga_\la^{\ominus}\,\overline{E}_c(\La)
\prod_{i=1}^n\prod_{ j=1}^{\infty}\frac{1}{1-q_i^{j}}.
\end{align}
\end{theorem}
{\em Proof.} The key claim here is $(ii)$. It follows
from Theorem \ref{GLOBNSSPH} and provides the
existence of the operators in $(i)$.
Let us demonstrate this in the case of $\hat{Y}_b$.
One uses (\ref{LfLa}) as follows:
\begin{align}
Y_b(q^{\frac{(c_\sharp,c_\sharp)}{2}-(\rho_k,\rho_k)}
E_c^\star(X)\,(\tga_x^\ominus)^{-1})
&=\lan Y_b(G(X,\La))\, E_c^\star(\Lambda)\,\tga_\la^\ominus\,
\mu_\circ(\La) \ran,\notag\\
&=\lan \La_{b}^{-1}G(X,\La)\, E_c^\star(\Lambda)\,\tga_\la^\ominus\,
\mu_\circ(\La) \ran.
\label{YG}
\end{align}
Applying $\kapp^\de$ and taking $t_\nu\to 0$, the
right-hand side of (\ref{YG}) is well defined.
Hence it follows that the action of $\RE^\de(Y_b)$ is well
defined on the spinor
$$
\RE^\de(q^{-(c_-,\rho_k)-\frac{(\rho_k,\rho_k)}{2}}
E_c^\star(X)\,(\tga_x^\ominus)^{-1})=
(\sum_w a_{c,w}\,X_{-c_-}\,\ze_w)\,(\tga_x^\ominus)^{-1}.
$$

In general,
one obtains that the action of the operators from $(i)$ is
well defined when they are applied to linear combinations
of spinors of the form $X_b\,\ze_w\,(\tga_x^\ominus)^{-1}$ for 
dominant regular $b$ and $w\in W$.
The operators $\kapp^\de(H)$ have rational coefficients;
nevertheless, this property is sufficient to see that their
coefficients are well defined in the limit $t_\nu\to 0$.
Moreover, this gives that the coefficients of 
$\RE^\de(H)$ actually belong to $\Q_q'[X_b,b\in B]$, i.e.,
do not have nontrivial denominators. As a matter of fact,
this can be formally deduced from (\ref{Gtal})
(see also (\ref{prodforgb})) and the fact that $\RE^\de(H)$ exist.

We will give below a direct and constructive proof of the 
existence of these operators 
and the absence of the denominators.
Actually, we will prove a stronger result based on direct
calculation of the $R\!E^\de$\~limits of the Dunkl operators,
which, for instance, allows to obtain the formulas for the
leading terms of the spinor Toda-Dunkl operators
and clarify their structure (including the
analysis of the vanishing coefficients). 

The claims in $(iii)$ follow directly from Theorem \ref{GLOBNSSPH}.
This is a nonsymmetric generalization of the $q$\~Shintani
formulas, which in turn generalize the classical Shintani-type
formulas in the theory of $p$\~adic Whittaker functions. Note
that the identity from (\ref{ShinqtWW}) does not contain/require
spinors; it involves only the $\overline{E}$\~polynomials and the
coefficients $\{a_{b,u}\}$. 
\sq
\smallskip

{\em Symmetrization.}
The symmetric (nonspinor) $q$\~Whittaker function $\w(X,\La)$
constructed in Theorem 3.2 of \cite{C10} is the symmetrization
of $\Om(X,\La)$. More precisely, one has
\begin{align}
\label{omhatsym}
\de(\w(X,\La))=\sum_{w\in W}\hat{T}_w (\Om(X,\La))
=\sum_{w\in W}T_w^\la(\Om(X,\La)).
\end{align}
In particular, the right-hand side is a diagonal spinor
(in the image of $\de$); all its components coincide. 
See Propositions \ref{SymDun} and \ref{OtherGens} below.

We note that the (nonaffine) hat-symmetrizer
$\sum_{w \in W}\hat{T}_w$ preserves the $\id$\~component of
any $W$\~spinor, which can be readily deduced from 
(\ref{kapp-de-X}) and formulas for $\ddot{T}_i$ acting in
the polynomial representation. See (\ref{RE-de-T}) below
for explicit formulas for $\hat{T}_i (i>0)$, which are 
sufficient to check this claim. Therefore it is not
actually necessary to perform the symmetrization in
(\ref{omhatsym}) because the $\,id$\~component of $\Om(X,\La)$
is exactly $\w(X,\La)$. Using (\ref{ShinqtW})
and formulas (\ref{GsymTW})
for $\hat{T}_i$ with $i>0$, we see that 
this coincidence can be also deduced from formula (\ref{eplimbar}), 
which states that $\overline{E}_{b}\ =\ \overline{P}_{b}$ 
for $b\in B_-$.

Explicitly, one has
\begin{align}\label{wsymsum}
&\w(X,\La)=(\tga_x^{\ominus}\tga_\la^{\ominus})^{-1}\,
\sum_{b\in B_-} q^{(b,b)/2}\,
\frac{X_{b_+}\overline{P}_b(\La^{-1})}
{\prod_{i=1}^n\prod_{j=1}^{-(\al_i^\vee,b)}(1-q_i^j)},\\
&\tga^{\ominus}(1)\w(q^{c},\La)\, =\, 
\overline{P}_c(\La)
\prod_{i=1}^n\prod_{ j=1}^{\infty}\frac{1}{1-q_i^{j}}
\for c\in B_-,\notag\\
\label{ShinqtWP}
&\sum_{b\in B_-} 
\frac{q^{(b-c)^2/2}\overline{P}_b(\La)}
{\prod_{i=1}^n\prod_{j=1}^{-(\al_i^\vee,b)}(1-q_i^j)}
= 
\tga_\la^{\ominus}\,\overline{P}_c(\La)
\prod_{i=1}^n\prod_{ j=1}^{\infty}\frac{1}{1-q_i^{j}}.
\end{align}
The latter formula results from (\ref{ShinqtWW})and
the equality $\overline{E}_b=\overline{P}_b$ for $b\in B_-$.
Indeed, when $c\in B_-$ and hence $u_c=$ id, the coefficient
$a_{b,u_c^{-1}}$ is nonzero only for $b\in B_-$; in this case,
one has $a_{b,id}=1$ and the summation in (\ref{ShinqtWP})
ranges over $b\in B_-$.

It is worth mentioning that  the left-hand side of (\ref{ShinqtWP}) 
becomes zero in the ($p$\~adic) limit $q\to 0$ when $c\not\in B_-$, 
which gives that
\begin{align}
\label{ShinqtWPAA}
\lim_{q\rightarrow 0} \w(q^{c},\La)\, =\,0
\hbox{\ \, unless\ \,} c\in B_-.
\end{align}

\comment{
For $c\in B_-$
and $i>0$, 
\begin{align}
\label{ShinqtWPA}
\hat{T}_i(\Om(q^{c},\La))\, =\, 
\overline{T}_i^\la(\w(q^{c},\La))\,=\, 
\tga^{\ominus}(1)^{-1}\overline{T}_i^\la
(\overline{P}_c(\La))
\prod_{i=1}^n\prod_{ j=1}^{\infty}\frac{1}{1-q_i^{j}}=0
\end{align}
Using formula (\ref{RE-de-T}) below, which reads
\begin{align}
\RE^\de(\ddot{T}_i)=
\sum_{\substack{w\in W\,\mathrm{s.t.}\\w^{-1}(\al_i)<0}}
\ze_w\,(s_i-1) \for i>0,
\label{RE-de-TT}
\end{align}
}

\subsection{\bf Fourier transform}\label{sec:duallimit}
Let us interpret the function 
$\Omega(X,\La)$ as the reproducing kernel of the nil-DAHA
Fourier transform.

The $X\leftrightarrow Y$ symmetry of $\HH^\flat$ is not present in
$\bHH^\flat$. More precisely, the anti-involution $\vph$
from (\ref{phianti}) does not act in $\bHH^\flat$.
To recover this symmetry, we define
$\dHH^{\flat,\vph}\equal\vph(\dHH^\flat)$,
a $\ddot{\Q}_{q,t}'$\~subalgebra of $\HH^\flat$;
see (\ref{dQprime}) for the definition of $\ddot{\Q}_{q,t}'$.

As an abstract algebra,
$\dHH^{\flat,\vph}$ can be described as
follows. Let $\breve{\Pi}^\flat\equal\vph(\Pi^\flat)$,
which is isomorphic to the abelian group $\Pi^\flat$, and let
$\breve{\pi}_r=\vph(\pi_r)$ for $r\in O'$,\,
$\breve{T}_i=\vph(\ddot{T}_i)$ for $0\le i\le n$.
Note that $\breve{T}_i=\ddot{T}_i$ unless $i=0$.
We set $\breve{T}_i'\equal\vph(\ddot{T}_i')$, i.e.,
$\breve{T}_i'=\breve{T}_i-(t_i-1)$.
Then in terms of the generators
$$ \breve{\Pi}^\flat,\ \ \breve{T}_i \ (0\le i\le n),
\ \ Y_b \ (b\in B),$$
the defining relations for $\dHH^{\flat,\vph}$ are:

($o$)\ \ $(\breve{T}_i-t_i)(\breve{T}_i+1)\ =\ 0,
\for 0\ \le\ i\ \le\ n$;

($i$)\ \ \ $\breve{T}_i\breve{T}_j\breve{T}_i\cdots\ =\
\breve{T}_j\breve{T}_i\breve{T}_j\cdots,\ m_{ij}$
factors on each side;

($ii$)\ \ $ \breve{\pi}_r^{-1}\breve{T}_i\breve{\pi}_r\ =\
\breve{T}_j \iif \pi_r(\al_i)=\al_j$,\, $\pi_r\in \Pi^\flat$;

($iii$)\ $\breve{T}_i' Y_b \ =\ Y_b Y_{\al_i}^{-1}
\breve{T}_i \iif (b,\al^\vee_i)=1,\ 0 \le i\le  n$;

($iv$)\ $\breve{T}_i Y_b\ =\ Y_b \breve{T}_i$ if $(b,\al^\vee_i)=0,
0 \le i\le  n$;

($v$)\ \ $\breve{\pi}_r^{-1} Y_b \breve{\pi}_r\ =\
Y_{\pi_r(b)}\ =\ Y_{u^{-1}_r(b)} q^{(\om_{r^*},b)}$.

The elements $\tilde{X}_b \ (b\in B)$ from
(\ref{hatoprs}) belong to $\dHH^{\flat,\vph}$; they
are the images of $\ddot{Y}_{-b}$ under $\vph$.
Namely, $X_b=\vph(Y_{-b})$ and
$\tilde{X}_b\equal\vph(\ddot{Y}_{-b})=q^{(b_+,\rho_k)}X_b$;
in the last equality, we use that $(-b)_+=-w_0(b_+)$ and
$(-w_0(b_+),\rho)=(b_+,\rho)$.
The corresponding relations in $\dHH^{\flat,\vph}$
are obtained by applying $\vph$ to those from
$\dHH^\flat$.

By construction, the anti-involution $\vph$ of $\HH^\flat$
sends $\dHH^{\flat}$ to $\dHH^{\flat,\vph}$.
The automorphism $\si$ from (\ref{si}) also has this
property. Explicitly, $\si$ preserves $q$ and
\begin{align}
\si : \ \label{sidot}
&\ddot{T}_i \mapsto \breve{T}_i \ (i\ge 0),
X_b \mapsto Y_b^{-1},\\ 
&\pi_r \mapsto X_{\om_r}T_{u_r^{-1}}=
\vph(\pi_r^{-1})=\breve{\pi}_r^{-1}.\notag
\end{align}

Now define $\bHH^{\flat,\vph}$ to be the
specialization of $\dHH^{\flat,\vph}$ for
all $t_\nu=0$. 
The anti-isomorphism $\vph$ and the automorphism $\si$
are compatible with this specialization, and we use
the same symbols to denote the resulting maps from
$\bHH^{\flat}$ to $\bHH^{\flat,\vph}$.

Recall the definition of the algebra $\bHH^{\flat,\dag}$ from
Section \ref{sec:infty}.
The involution $\vep$ of $\HH^\flat$ from (\ref{ve}) conjugates
$q, t_\nu$.
Accordingly, we obtain an isomorphism
$\vep^\dag : \bHH^{\flat,\dag} \to \bHH^{\flat,\vph}$ given by
\begin{align}
\label{epdag}
\vep^\dag : \ 
&\overline{T}_i^\dag\mapsto \overline{T}_i^{\,\prime} \ (i>0),\
X_b \mapsto Y_b,\ q \mapsto q^{-1},\\
&\pi_r \mapsto \breve{\pi}_r^{-1}=X_{\om_r} T_{u_r^{-1}}
=\vph(\pi_r^{-1}).\notag
\end{align}

Now we are ready to state the Fourier transform interpretation
of Theorem \ref{GLOBNSWH}.

\begin{corollary}
Define the transforms
\begin{align}
&\bS(f)(X) \equal 
\lan f(\La)\,\Om(X,\La)\,\overline{\mu}_\circ\ran_\circ\label{SFT},\\
&\bE(f)(X) \equal
\lan f(\La)^\ast\,\Om(X,\La)\,\overline{\mu}_\circ\ran_\circ\label{EFT},
\end{align}
acting from functions of $\La$ to $X$\~spinors.
Then one has:
\begin{align}
\bS((\overline{E}_b^\dag)^\ast\tga_\la^\ominus)=
\bE(\overline{E}_b^\dag\tga_\la^\oplus)=
q^{(b,b)/2}(\tga_x^\ominus)^{-1}
\sum_{w\in W} a_{b,w}\,X_{-b_-}\,\ze_w,\label{FTimg}
\end{align}
and
\begin{align}
&\bS(H(f))=\si(H)\bS(f) \for
H \in \bHH^\flat,\ \si(H) \in \bHH^{\flat,\vph},
\label{SH}\\
&\bE(H(f))=\vep^\dag(H)\bE(f) \for
H \in \bHH^{\flat,\dag},\ \vep^\dag(H) \in \bHH^{\flat,\vph},
\label{EH}
\end{align}
provided the existence of the transforms $\bS(f)$ and $\bE(f)$.
\end{corollary}

\proof
The formula (\ref{FTimg}) is immediate from the explicit
expression (\ref{omform}) for $\Om$ and the orthogonality relations
(\ref{ebarnorms}). The intertwining properties 
(\ref{SH}) and (\ref{EH})
follow from the relations $\vep=\vph\star$ and
$\si=\vep\eta=\vph\!\star\!\eta$.
\sq
\medskip

{\em Pseudo-polynomial representation.}
The transform $\bS$ embeds $\overline{\v}_\la\,\tga^\ominus_\la$,
which is $\overline{\v}\,\tga^\ominus$ under $X\mapsto \La$,
into the space of $X$\~spinors.
The image $\bS(\overline{\v}_\la\,\tga_\la^\ominus)$
can be described explicitly as follows.

We need the following properties of the coefficients $a_{b,w}$
from (\ref{omform}):
\begin{align}
&a_{b,w}=0 \unl w \ge u_b^{-1},\label{abw0}\\
&a_{b,w}=a_{b,wy} \iif y \in W^{b_-}.\label{abwy}
\end{align}
Recall that $W^b$ is the centralizer of $b$ in $W$.
These properties are immediate from the description of the
$a_{b,w}$ given in Proposition \ref{OMLIMIT}. For (\ref{abw0}),
one must also use (\ref{macsucc}).
Formula (\ref{abwy}) asserts that $a_{b,w}$ depends only
on the coset of $w$ in $W/W^{b_-}$.
The elements $u_b^{-1}$ for $b\in W(b_-)$
are exactly the minimum length coset representatives for 
$W/W^{b_-}$. 

For any $b\in B$, we set
\begin{align}
\mathbf{X}_b=\sum_{y\in u_b^{-1} W_{b_-}} X_{-b_-}\,\ze_y
\ =\ X_{-b_-}\!\!\!\!\!\sum_{y\in u_b^{-1} W_{b_-}}\,\ze_y.
\end{align}
Then we can write (\ref{FTimg}) as follows:
\begin{align}
\bS((\overline{E}_b^\dag)^\ast\tga_\la^\ominus)=
q^{(b,b)/2}(\tga_x^\ominus)^{-1}
\sum_{\substack{c\in W(b)\\u_c^{-1}\ge u_b^{-1}}} 
a_{b,u_c^{-1}} \mathbf{X}_c,
\end{align}
We conclude that the image $\bS(\overline{\v}_\la\,\tga_\la^\ominus)$
is precisely the span of the elements
$\{\mathbf{X}_b\tga_x^\ominus , b\in B\}$.

Note that when $b\in B$ is regular, i.e., $W^{b_-}=\{\id\}$,
then one simply has $\mathbf{X}_b = X_{-b_-}\,\ze_{u_b^{-1}}$. 
Hence the image of $\bS$ contains all elements 
$X_b\,\ze_w\,\tga_x^\ominus$
for dominant regular $b$ and $w \in W$.
This observation was used above in the proof of 
Theorem\~\ref{GLOBNSWH}
to establish the existence of the limits
$\RE^\de(H)$ for all $H\in \dHH^{\flat,\vph}$ as spinor 
difference-reflection
operators; in particular, it establishes the existence of
$\hat{Y}_b$ for all $b\in B$, which we call
the {\em Toda-Dunkl operators}.

The remainder of the paper is devoted to a direct proof of
the existence of these operators,
which is quite interesting in its own right and provides
valuable information about the structure of these operators.
\smallskip

\setcounter{equation}{0}
\section{\sc Managing G-products}
We will give in the next two sections a 
constructive justification of the existence of
the Toda-Dunkl operators and the other operators
$\RE^\de(H)$ for $H\in\dHH^{\flat,\vph}$, without
any reference to the spinor Whittaker function
$\Om(X,\La)$ obtained above. We assume that $B=P$
for the remainder of the paper, which is sufficient
for establishing the existence and finding the formulas. 

The direct approach involves 
some nontrivial combinatorics of reduced decompositions
and $\la$\~sequences in $\hW$ but gives more exact information
about the structure and the coefficients of such operators.
We will begin with basic estimates and examples and then we 
will proceed by induction in the next section.

\subsection{\bf Basic \texorpdfstring{$t$}{\em t}-estimates}
\label{sec:t-est}
We will denote $\spin(\Q_{q,t}'(X))$ by $\spin(\v')$ in this
and further sections.
For any $f\in\spin(\v')$, $w\in W$, and $\nu \in \nu_R$, we define 
$\ord_w^\nu(f)$ to be the order of the $w$-component
$f_w$ with respect to $t_\nu$.
Hence for $f,g\in\spin(\v')$,
\begin{align}
\ord_w^\nu(fg)=\ord_w^\nu(f)+\ord_w^\nu(g).
\label{ordw-fg}
\end{align}
Note that
\begin{align}
\ord_w^\nu(\de(v)f)=\ord_{v^{-1}w}^\nu(f),
\label{ordw-vf}
\end{align}
due to (\ref{vwf}).

For $f\in\spin(\v')$ and $v,w\in W$, define
\begin{align}
\label{def-ord}
\ord_w^\nu(f\,v)\equal\ord_w^\nu(f),\ \ 
\dord_w^\nu(v\,f)\equal\ord_w^\nu(f).
\end{align}
More generally, we set $\ord_w^\nu (Aw)=\ord_w^\nu(A)=\dord_w^\nu(wA)$
for $w\in W$ and any spinor difference operators $A$,
the sums of the products of $X$\~spinors and 
$\Ga$\~spinors.
\smallskip

Recall that given a reduced decomposition 
$\hu= \pi_r s_{j_l}\cdots s_{j_1}$ for $l=l(\hu)$
and $r\in O$,
\begin{align}\label{prodforgb}
&T_{\hu}\ =\ \hu\, G_{\tal^{\,l}}\cdots G_{\tal^1}\ =\ 
\ G_{-\tbe^{\,l}}\cdots G_{-\tbe^1}\,\hu\\
&\hbox{for\ \ }  
\tal^1=\al_{j_1},
\tal^2=s_{j_1}(\al_{j_2}),\,\ldots\,,\, 
\tbe^r=-b(\tal^{\,r})\in \tR_+,\notag\\
\label{Gtaldef}
&G_{\tal} \ \equal\ t_\al^{1/2}+\frac{t_{\al}^{1/2}-t_{\al}^{-1/2}}
{X_{\tal}^{-1}-1}\,(1-s_{\tal})\\ 
&\ \ \ \ \ \  =\ \frac{t_\al^{1/2}X_{\tal}^{-1}-t_{\al}^{-1/2}}
{X_{\tal}^{-1}-1}-
\frac{t_\al^{1/2}-t_\al^{-1/2}}{X_{\tal}^{-1}-1}\,s_{\tal},\notag\\ 
&G_{-\tal}\ =\ \frac{t_\al^{1/2}X_{\tal}-t_\al^{-1/2}}{X_{\tal}-1}-
\frac{t_\al^{1/2}-t_\al^{-1/2}}{X_{\tal}-1}\,s_{\tal},\notag
\end{align}
where $\tal\in\tR$. 
Recall that 
$X_{\tal}=X_{\al}q^j$\, for $\tal=[\al,j]\,$
and that $\{\tal^1,\tal^2,\ldots,\tal^l\}=\la(\hu)\subset\tR_+$. 
Note that $\hu\, s_{\tal^l}\,\cdots\, s_{\tal^1}=$ $\pi_r.$

We will restrict ourselves to the estimates of the orders
of operators $Y_b^{-1}$ for $b\in P_+$. This is totally
parallel to the case of  $Y_b$ and complementary to the proof 
of Proposition \ref{prop-ord-bound} provided below, where the 
case of positive powers will be addressed.  

Accordingly, 
\begin{align}\label{prodforg}
&T_{\hu}^{-1}\ =
\tilde{G}_{\tal^1}\cdots \tilde{G}_{\tal^l}\, \hu^{-1},\ \,
\tilde{G}_{\tal}=
\frac{t_{\al}-X_{\tal}^{-1}}{1-X_{\tal}^{-1}}+
\frac{t_\al-1}{1-X_{\tal}^{-1}}s_{\tal}.
\end{align}
where, as above, $\tal^1=\al_{j_1}$, $\tal^2=s_{j_1}(\al_{j_2})$,
$\tal^3=s_{j_1}s_{j_2}(\al_{j_3})$ and so on. We set
\begin{align}\label{fgtilde}
\tilde{f}_{\tal}\ =\ \frac{t_\al-X_{\tal}^{-1}}
{1-X_{\tal}^{-1}},\ \ 
\tilde{g}_{\tal}\ =\ \frac{t_\al-1}{1-X_{\tal}^{-1}},
\ \ \tilde{G}_{\tal}\ =\ \tilde{f}_{\tal}+\tilde{g}_{\tal}s_{\tal}\,.
\end{align}

Our aim is to estimate the orders of operators
$\kapp^\de(T_{\hu})$ with respect to $t_\nu$ for
$\nu\in \nu_R$ aiming at the limit $t_\nu\to 0$ for all $\nu$. 
Recall that for $b\in P$ and $\nu\in \nu_R$, 
\begin{align}\label{qrhokb}
&q^{\,(\rho_k\,,\,b)}\,=\,q^{\,\sum_\nu k_\nu(\rho_\nu\,,\,b)}\,=\,
\prod_{\nu}q_\nu^{\,k_\nu\,(\rho_\nu^\vee\,,\,b)}\,=\,
\prod_{\nu}t_\nu^{\,(\rho_\nu^\vee\,,\,b)},\\
&\hbox{also,\ \,} t^{\,l(\hw)/2}\, =\, \prod_{\nu} 
t_\nu ^{\, l_\nu(\hw)/2} = t_\nu ^{\,|\,\la_\nu(\hw)\,|/2}
\for \hw\in \hW.\notag
\end{align}

The following straightforward formulas are actually the key. 
For any $\nu\in \nu_R$, $\tal=[\al,\nu_\al j]$ and 
$\de_{\al,\nu}\equal\de_{\nu_\al,\nu}$,
\begin{align}
\label{ord-g-taln}
\ord_w^\nu(\kapp^\de(\tg_{\tal}))\,&=\,
\begin{cases}
\ \ \ \ \ \ \ \ 0,&\text{if } w^{-1}(\al)>0,\\
-(w^{-1}(\al),\rho^\vee_\nu),&\text{if } w^{-1}(\al)<0,
\end{cases}
\end{align}
and
\begin{align}
\label{ord-f-taln}
\ord_w^\nu(\kapp^\de(\tf_{\tal}))&\,=\,
\begin{cases}
\de_{\al,\nu},&\text{if } w^{-1}(\al)>0,\\
\ 0,&\text{if } w^{-1}(\al)<0.
\end{cases}
\end{align}
Indeed,
\begin{align}
\label{ord-g-talnn}
\kapp^\de(\tg_{\tal})&=
\sum_{w\in W}\frac{t_\nu-1}
{1-q^{(\rho_k\,,\,w^{-1}(\al))}w^{-1}(X_\al^{-1})}\,\ze_w,\\
\label{ord-f-talnn}
\kapp^\de(\tf_{\tal})&=
\sum_{w\in W}\frac{t_\nu\!-q^{(\rho_k\,,\,
w^{-1}(\al))}w^{-1}(X_\al^{-1})}
{1-q^{(\rho_k\,,\,w^{-1}(\al))}w^{-1}(X_\al^{-1})}\,\ze_w.
\end{align}
We also need the $t_\nu$\~orders of $\Ga_b=(-b)$\,:
\begin{align}
\ord_w^\nu(\kapp^\de(\Ga_b))&=
\ord_w^\nu(\sum_{v\in W}
q^{-(\rho_k\,,\,v^{-1}(b))}\Ga_{v^{-1}(b)}\,\ze_v)
=-(\rho_\nu^\vee,w^{-1}(b)).\notag
\end{align}
Applying this formula to $\,s_{\tal}s_{\al}\ =\ 
\Ga_{\al}^{\,j}\,$ for 
$\tal=[\al,\nu_\al j]\in \tR$, we obtain that
$$
\ord_w^\nu(\kapp^\de(s_{\tal}))=
-j\,(\rho^\vee_\nu,w^{-1}(\al)), \and
$$
\begin{align}
\label{ord-g-talns}
\ord_w^\nu(\kapp^\de(\tg_{\tal}s_{\tal}))=&
\begin{cases}
\ \ \ \ \ -j\,(\rho^\vee_\nu,w^{-1}(\al)),
&\text{if } w^{-1}(\al)>0,\\
-(1+j)\,(\rho^\vee_\nu,w^{-1}(\al)),
&\text{if } w^{-1}(\al)<0.
\end{cases}
\end{align}

Recall that we set $\ord_w^\nu (Aw)=\ord_w^\nu(A)=\dord_w^\nu(wA)$
for $w\in W$ and any spinor difference operators $A$,
the sums of the products of $X$\~spinors and 
$\Ga$\~spinors. One also has  that
\begin{align*}
\tg_{\tal}s_{\tal}\ =\ s_\al \Ga^{-j}\tg_{-\tal}, \and
\end{align*}
\begin{align}
\label{ord-g-talnspr}
\dord_w^\nu(\kapp^\de(\tg_{\tal}s_{\tal}))=&
\begin{cases}
(1+j)\,(\rho^\vee_\nu,w^{-1}(\al)),\,
&\text{if } w^{-1}(\al)>0,\\
\ \ \ j\ (\rho^\vee_\nu,w^{-1}(\al)),\,
&\text{if } w^{-1}(\al)<0.
\end{cases}
\end{align}

\smallskip
\subsection{\bf Leading terms}
In this section, we provide the estimates for {\em some}
of the coefficients of the $Y$\~operators, which will be
main ingredients of the direct justification of the existence 
(and invertibility) of $\RE(Y_b)$ in the next section.

It is straightforward
to calculate the top coefficient in the expansion
of $\ddot{T}_{\hu}^{-1}=t^{l(\hu)/2}T_{\hu}^{-1}
\equal \sum_{\hw\in\hW} \ddot{C}_{\hw}\hw$, 
which is $\ddot{C}_{\hu^{-1}}$. 
It can be obtained only by picking the terms without 
$\,s_{\tal}\,$ from all binomials in (\ref{prodforg}). Thus 
\begin{align}\label{cminbt}
\ddot{C}_{\hu^{-1}}\  =\ \prod_{\tal\in \La(\hu)} \tf_{\tal}\ =\ 
\prod_{[\al,j]\in \La(\hu)} \frac{1-t_\al q^j X_{\al}}{1-q^j X_{\al}}. 
\end{align}

Let us apply this formula to $\hu=b\in P_+$. Then 
$$
\la_\nu(b)=\{\tal=[\al,\nu_\al j]\,;\, 
\al>0,\, \nu_\al=\nu,\, 0\le j<(b,\al^\vee)\}
$$
and $(\rho^\vee_\nu-w(\rho^\vee_\nu),b)$ is exactly the number of
$\tal\in \la_\nu(b)$ such that $w^{-1}(\al)<0$, which coincides
with the number of $\tf_{\tal}$ for $\tal\in \la_\nu(b)$ such
that $\ord_w^\nu(\tf_{\tal})=0$. 

Therefore, 
\begin{align}
\label{ord-greatest}
&\ord_w^\nu(\kapp^\de(\ddot{C}_{-b}\Ga_b))=
\sum_{\tal\in\la_\nu(b)}\ord_w^\nu(\tf_{\tal})-
(\rho^\vee_\nu,w^{-1}(b))\\
&=2(\rho^\vee_\nu,b)-(\rho^\vee_\nu-w(\rho^\vee_\nu),b)-
(\rho^\vee_\nu,w^{-1}(b))
=(\rho^\vee_\nu,b).\notag
\end{align}
We see that the 
term $q^{-(\rho_k\,,\,b)}\kapp^\de(\ddot{C}_{-b}\Ga_b)$
in the $\kapp^\de$\~image of the expansion 
\begin{align}\label{ybhwc}
Y_b^{-1}=q^{-(\rho_k\,,\,b)}\ddot{Y}_b^{-1}=q^{-(\rho_k\,,\,b)}
\sum_{\hw\in \hW}\ddot{C}_{\hw}\hw 
\end{align}
is of zero order, as it is supposed to be due to the
existence and invertibility of $\RE(Y_b^{-1})$.
\smallskip

We will prove below that
the order $\ord_w^\nu(\kapp^\de(\Pi))$ is non-negative
for each $w\in W$ and any individual product $\Pi$  contributing to 
$Y_b^{\pm 1}\, (b\in P_+)$, where we pick
terms with and without $s_i$ from 
$T_i^{\pm 1}$ in the polynomial representation. Equivalently,
we can expand $Y_b^{\pm 1}$ with $\hu$ placed on the left.
In the case of $Y_b^{-1}$, which we mainly consider in this section,
$$
\ord_w^\nu(\kapp^\de(\Pi'_{\hw}))\,
\ord_w^\nu(\kapp^\de(\hw))\,\ge\, (\rho^\vee_\nu,b)\,\le\,
\ord_w^\nu(\kapp^\de(\Pi''_{\hw}))\,
\ord_w^\nu(\kapp^\de(\hw))\,
$$
for any individual products $\Pi'_{\hw}$ of $\tf_{\tal^j}$ and 
$\tg_{\tal^j}$ contributing to $\ddot{C}_{\hw}$ from (\ref{ybhwc})
and their counterparts $\Pi''_{\hw}$ for $\dord$, i.e.
for  $\ddot{Y}_b^{-1}= \sum\, \hw \,\Pi''_{\hw}\,$ with $\hw$ collected
on the left. 
\smallskip

\subsection{\bf The lowest terms}
The smallest possible $\hw$ that can be obtained
from $\hu$ is when we always
pick the terms with $\,s_{\tal}\,$ from the binomials
in the product (\ref{prodforg}) for 
$\ddot{T}_{\hu}^{-1}\,=\,t^{l(\hu)/2}T_{\hu}^{-1}
\,=\,\prod_{\nu}t_\nu^{l_\nu(\hu)/2}\,T_{\hu}^{-1}\,.$

It will contribute to the $\ddot{C}$\~coefficient of $\pi_r^{-1}$ 
(maximally distant from $\hw^{-1}$).  There can be of
course other products that contribute to 
$\ddot{C}_{\pi_r^{-1}};$ their number grows exponentially 
with $l(\hu)$. This particular product is as follows:
\begin{align}\label{prodforpi}
&\prod_{r=1}^l \tg_{\al_{j_r}}\ =\ \prod_{r=1}^l
\frac{t_{j_r}-1}{1-X_{\al_{j_r}}^{-1}}\,; 
\end{align}
recall that $\{\al_j\}$ are simple roots.

Let $\nu\in \nu_R,\, w\in W$ and $N^\nu_j(b;w)$ be the number 
of simple reflections $s_{j}$ for $\nu_j=\nu$ in a given 
decomposition of $b\in P_+$ such that $w^{-1}(\al_j)<0$. 
Accordingly, $N_0^\nu$ is zero when $w^{-1}(\vth)<0$ or 
$\nu=\nu_{\lng}$ and the number of $s_0$ in this decomposition 
of  $b$ otherwise.
We set $N_j^\nu(b)$ if all $s_j$ with $\nu_i=\nu$
are counted, and omit the super-index $\nu$ if all 
$\nu\in\nu_R$  are considered. 
Actually here we need only the decomposition of $\pi_r^{-1}b$, 
where $b-\om_r\in Q$.

In such case, we claim that for any reduced 
decomposition of $b=b_+$ and for an arbitrary $w\in W$,
\begin{align}\label{nonjw}
&\ord_w^\nu\Bigl(q^{-(\rho_k\,,\,b)}\kapp^{\de}\bigl(
(\prod_{r=1}^l \tg_{\al_{j_r}})\pi_r^{-1}\bigr)\Bigr)\, \ge\, 0\,
\hbox{,\ \, equivalently\,,}\\
N_0^\nu(b;w) &(\vth,w(\rho_\nu^\vee))\!-\!\!\sum_{j=1}^n\, 
N_j^\nu(b;w)\,(\al_j,w(\rho_\nu^\vee))
\!+\!(\om_{r^*},w(\rho_\nu^\vee))\ge (b,\rho_\nu^\vee).\notag
\end{align}
Here we use that $\pi_r^{-1}=u_r\Ga_{b_r}$ and  
therefore 
\begin{align*}
\kapp^\de(\pi_r^{-1})\!=\!
u_r\bigl(
q^{-(v^{-1}(\om_r),\rho_k)}
\!\!\sum_{v\in W}\Ga_{v^{-1}(\om_r)}^{-1} \ze_v\bigr)\!=\!
q^{-(v^{-1}(\om_r),\rho_k)}
\!\!\sum_{v\in W}\Ga_{v^{-1}(\om_r)}^{-1} \ze_{u_r v}\,&
\\
=\ q^{-(w^{-1}u_r(\om_r),\rho_k)}\sum_{w\in W}
\Ga_{w^{-1}u_r(\om_r)}^{-1} 
\,\ze_w
\ =\ q^{(w^{-1}(\om_{r^*}),\rho_k)}\sum_{w\in W}\Ga_{w^{-1}
(\om_{r^*})}^{+1} \,\ze_w.&
\end{align*}
Recall that $\om_{r^*}=-u_r(\om_r)$ for the 
image $r^*$ of $r$ under the involution of the nonaffine
Dynkin diagram on the set $O$ induced by $-w_0$. Also,
$(b,\rho^\vee_\nu)=(-w_0(b),\rho^\vee_\nu)$ for $\nu\in\nu_R$.

The existence of a reduced decomposition of $b$ 
such that (\ref{nonjw}) holds for any $w$ can be justified
directly using that the left-hand side of (\ref{nonjw}) is
additive for $b=a+c$, where $a\in Q\cap P_+$ and $c\in P_+$
(so is its right-hand side). We use here that $l(a+c)=l(a)+l(c)$
and therefore $N_j^\nu(a+c;w)=N_j^\nu(a;w)+N_j^\nu(c;w)$ if the reduced 
decomposition of $a$ and $c$ are combined together (without further 
using the homogeneous Coxeter transformations). 

For instance, for
$w=\id$ and $b\in P_+$, the inequality becomes the estimate
$N_0(b)\ge (b-\om_r,\rho^\vee)/(\vth,\rho^\vee)$. Taking
$w=s_j(j>0)$ and assuming that the root system is not $A_1$,
we arrive at 
$$
N_0(b)(\vth,\rho^\vee-\al_j^\vee)+
N_j(b)\ge (b-\om_r,\rho^\vee)+\de_{r^*,j}.
$$
For $b=\om_r$, it means that $N_j(\om_r)\ge \de_{r^*,j}$ and
that $u_r=\pi_r^{-1}\om_r$ contains 
at least one $s_{r^*}$ in any of its reduced decompositions.
This holds since otherwise $u_r(\om_{r^*})=\om_{r^*}$, 
which is impossible because $u_r(\om_{r^*})=-\om_{r}$.

When $w=w_0$, we obtain the inequality, which is of some 
interest and clarifies the combinatorics related to our
estimates:
$$
\sum_{j=1}^n N_j(b)\ =\ 2(b,\rho^\vee)-N_0(b)\ 
\ge\ (b+\om_r,\rho^\vee).
$$
It can be readily transformed to
$N_0(b)\ \le\ (b-\om_r,\rho^\vee)$
and is sharp for $b=\om_r$, $r\in O'$.
\medskip

\subsection{\bf Taking one \texorpdfstring{$g$}{\em g}}
Let us now  
pick only one term $\tg_{\tal^i}s_{\tal^i}\, (1\le i\le l)\,$ 
from the binomials in the product (\ref{prodforg}) and then 
move $\ts^i\equal s_{\tal^i}$ to the left.
We will also use the notation $s^i\equal s_{\al^i}$ for
$\tal^i=[\al^i,\nu_{\al^i} j]$. Recall that
\begin{align}\label{prodforgr}
&T_{\hu}^{-1}\ =
\tilde{G}_{\tal^1}\cdots \tilde{G}_{\tal^l}\, \hu^{-1},\ \ 
\tilde{G}_{\tal}=
\frac{t_{\al}-X_{\tal}^{-1}}{1-X_{\tal}^{-1}}+
\frac{t_\al-1}{1-X_{\tal}^{-1}}s_{\tal}.
\end{align}

The 
contribution of this term to 
$\ddot{Y}_b^{-1}\equal 
q^{(\rho_k\,,\,b)}Y_b^{-1}$ for $b\in P_+$\, is
\begin{align}\label{prodforb1}
\ts^i\,\Bigl(\,\tf_{\ts^i(\tal^1)}\tf_{\ts^i(\tal^2)}
\cdots\tf_{\ts^i(\tal^{i-1})}
\bigl(\tg_{-\tal^i}\bigr)
\tf_{\tal^{i+1}}\cdots\tf_{\tal^{l}}\,\Bigr)\, \Ga_b.
\end{align}
We claim that the order of this expression is no smaller than
$(\rho^\vee,b)$. Using (\ref{ord-greatest}) and
(\ref{ord-g-talnspr}), it suffices to
check the following lemma.

\begin{lemma}\label{lemordoneg}
\begin{align}\label{ordwjump}
&\ord_w^\nu\bigl(\tf_{\ts^i(\tal^1)}\tf_{\ts^i(\tal^2)}
\cdots\tf_{\ts^i(\tal^{i-1})}\bigr)+
\dord_w^\nu(\tg_{\tal^i}(\ts^i))\\
&\ \ge \ord_w^\nu\bigl(\tf_{\tal^1}\tf_{\tal^2}
\cdots\tf_{\tal^{i}}\bigr).\notag
\end{align}
\end{lemma}
Let us begin its proof with the example of simplest $w=\id$.
We set here and below 
$\tal^j=[\al^j,\nu_{\al^j} m_j]$ for $j=1,\cdots, l$, 
for instance
$\tal^i=[\al^i,\nu_{\al^i} m_i]$, and
$\tbe^r\equal \ts^i(\tal^r)=[\be^r,\nu_{\be^r} k_r]$ for 
$r=1,2,\ldots,i$; thus $\tbe^i=\tal^i$. 

If 
$w=\id$, then $w^{-1}(\al^i)=\al^i>0$ and we need to check that
for $\tbe^r$ such that $\nu_{\be^r}=\nu$ 
and for a given $\nu\in \nu_R$,
\begin{align}\label{wisidm}
(m_i+1)(\rho^\vee_\nu,\al^i)\ge 
|\,\{\,\tbe^r\,;\, \be^r=s^i(\al^r)<0,\, 
 1\le r\le i\,\}\,|.
\end{align}
Using (\ref{talinla})
we obtain that if $\al^i\neq \be\in \la_\nu(s^i)$,
then for each $\,k\,$, either $[\be,\nu_\be k]$ or 
$[\be',\nu_\be k]$ for $\be'=-s^i(\be)$ can belong
to the set $\{\tbe^r,\, 1\le r \le i\}$, 
but not both. 

More exactly, Claim $(ii)$ from
Lemma \ref{LASTAL} gives that there will be exactly
one such occurrence for each $0\le k\le (\al^i,\be^\vee)m_i$.
The coefficient $(\al^i,\be^\vee)$ is $1$ unless
$\be$ is short and $\al^i$ is long, when it
equals $\nu_{\lng}$, i.e. it is $\eta_{\be,\al^i}$ in the
notation from this Lemma. 

Therefore the total number $\mathfrak{N}$ 
of $\tbe=[\be,\nu_\be k]\in \{\tbe^r\}$ such
that $\al^i\neq \be\in \la(s^i)$ is as follows:
\begin{align*}
\mathfrak{N}-(\rho,\al^\vee)=
&\ m_i\bigl((\rho_{\lng},\al^\vee)-1+\nu_{\lng}
(\rho_{\sht},\al^\vee)\bigr)\\
=&\ m_i\bigl((\rho_{\lng}^\vee,\al)-1+(\rho_{\sht},\al)\bigr)\\
=&\ m_i\bigl((\rho^\vee,\al)-1\bigr)
\hbox{\ \, for long\ }\al=\al^i,\\
\mathfrak{N}-(\rho^\vee,\al)=&\ m_i\bigl((\rho_{\lng}^\vee,\al)+
(\rho_{\sht},\al)-1\bigr)\\
=&\ m_i\bigl((\rho^\vee,\al)-1\bigr)
\hbox{\ \, for short\ }\al=\al^i.\\
\end{align*}
Accordingly, for $\mathfrak{N}_\nu$ defined when 
$\al^i\neq \be\in \la_\nu(s^i)$,
\begin{align*}
\mathfrak{N}_\nu-(\rho_\nu,\al^\vee)
=&\ m_i\bigl((\rho^\vee_\nu\,,\al)-\de_{\al,\nu})\bigr)
\hbox{\ \, for long\ }\al=\al^i,\\
\mathfrak{N}_\nu\,-\,(\rho_\nu^\vee,\al)
=&\ m_i\bigl((\rho^\vee_\nu\,,\al)-\de_{\al,\nu})\bigr)
\hbox{\ \, for short\ }\al=\al^i.\\
\end{align*}

Allowing the roots 
$\{\,[\al,\nu_\al k]\,,\, 0\le k\le m_i\,\},\,$
i.e. omitting the restriction $\be\neq\al^i$,
we obtain that the right-hand
side of (\ref{wisidm}) 
equals $(1+m_i)((\rho^\vee_\nu,\al^i)-
\de_{\al^i,\nu})+
\de_{\al^i,\nu}(1+m_i)$
for short $\al^i$, which
does coincide with the left-hand side, 
and is strictly smaller than the left-hand side
if $\al^i$ is long, $(\al^i,\rho_{\sht})\neq 0$
and $\nu=\nu_{\sht}$. Recall that the latter constraint
means that we calculate the order with respect
to $t_{\sht}$.
\smallskip

\comment{
\begin{align}
\la_\diamond(s_{\al})= 
\{\,[\al,\nu_\al k]\,&\mid\, 0\le k\le 2j\,\}\,\cup\,
\{\,\tbe>0\,\mid\, 0\le \nu_\be k<(\al,\be)j\,\}
\notag\\
\cup\,&
\{\,[\be,(\al,\be)j]\,\mid\, \be<0 \hbox{\ or\ }
\be\in \la(s_{\al})\,\}.\label{lashtlngax}
\end{align} 
}

As another example, let
us consider $w=w_0$. One needs to verify that
\begin{align*}
m_i(\rho^\vee_\nu,\al^i) \le 
|\,\{\,\tbe^r\,;\, \be^r=s^i(\al^r)<0,\,
\nu_{\al^r}=\nu=\nu_{\be^r},\, 1\le r<i\,\}\, |
\end{align*}
in this case.
The right-hand side for $1\le r\le i$, i.e. with
$i$ added to the range, has been already calculated above. 
So the required cardinality is
$m_i(\rho^\vee_\nu,\al^i)-\de_{\al^i,\nu}$ plus 
$(\rho^\vee_\nu,\al^i)$
for short $\al^i$ and  $(\rho_\nu,(\al^i)^\vee)$ for long
$\al^i$. Thus the inequality holds and
is strict unless $\al^i$ is a simple root and $\nu=\nu_i$.
\smallskip

\subsection{\bf One \texorpdfstring{$g$}{\em g} and any
\texorpdfstring{$w$}{\em w}}\label{sec:ONEGANYW}
This case is actually the key for the general consideration,
which will be managed by induction with respect to the number
of terms $\tg_{\tal^i}\ts^i$ taken in the products.

We will use that $\tbe^r=\ts^i(\tal^i)\not\in 
\{\tal^i,\ldots,\tal^1\}$ for $r<i\,$
if and only if $\tbe^r<0$. Indeed,
$$
\tbe^r=\tal^r-2\frac{(\al^i,\al^r)}{(\al^i,\al^i)}\tal^i
$$
and if $\tbe^r>0$, then either $\tbe^r$ is a sum of two roots from
$\la_\nu(b)$ with positive coefficients and therefore belongs to
$\{\tal^i,\ldots,\tal^1\}$ or $(\al^i,\al^r)>0$ and 
$\tal^r$ is such a sum of $\tbe^r$ and $\tal^i$. In the latter
case, $\tbe^r$ must belong to $\la_\nu(b)$ because otherwise
$\tal^i$ would occur before $\tal^r$ in this sequence; see 
Claim $(d')$ Theorem 2.1 from \cite{C103}. Moreover, $\tbe^r$
must appear before $\tal^i$ (actually before $\tal^r$) in this case.
This proves that  $\tbe^r<0$. 
\smallskip

For arbitrary $w\in W$ and $\nu\in\nu_R$, 
let us first consider $\tbe^r$ such that
\begin{align}\label{ordwtbe}
&\ord_w^\nu(\tf_{\tbe^r}) < \ord_w^\nu(\tf_{\tal^r})\for r<i,
\hbox{\  equivalently,}\\
&w^{-1}(\be^r)=w^{-1}(s_{\al^i}(\al^r))<0 \and w^{-1}(\al^r)>0.\notag
\end{align}
We assume here and below that $\nu_{\al^r}=\nu=\nu_{\be^r}$;
otherwise we have $0=0$ in (\ref{ordwtbe}).
Setting $\ep'=\sgn (w^{-1}(\al^i))$ and
$\al'\equal \ep'w^{-1}(\al^i)\in R_+$, we have 
\begin{align*}
&w^{-1}(\al^r) \in \la_\nu(s_{\al'}) 
\hbox{\ \, and therefore\ \,}
 (w^{-1}(\al^r),\al')=\ep'(\al^r,\al^i)>0.
\end{align*}
We can stick here only to negative $\tbe^r$ due to the remark above;
only such $\tbe^r$ may influence the difference 
$\sum_{r=1}^i\ord_w^\nu(\tf_{\tbe^r})-$
$\sum_{r=1}^i\ord_w^\nu(\tf_{\tal^r})$.
\medskip

$(a)$\, {\em The case of\,
$w^{-1}(\al^i)<0$.} 
Accordingly, $\ep'=-$ and $(\al^r,\al^i)<0$ for $\tal^r$ 
under consideration. However then the corresponding 
$\tbe^r=\tal^r-2\frac{(\al^i,\al^r)}{(\al^i,\al^i)}\tal^i$
belong to $\la_\nu(b)$ and, moreover, coincides with certain 
$\tal^p$ for $p<i$, as it was noted above. We see that 
\begin{align}\label{sumorwtf}
\sum_{r=1}^i\ord_w^\nu(\tf_{\tbe^r})-
\sum_{r=1}^i\ord_w^\nu(\tf_{\tal^r})\ge 0.
\end{align}
If $m_i=0$, i.e. $\al^i$
is nonaffine, this concludes the verification of 
(\ref{ordwjump}) in the considered case. Otherwise, the difference
from (\ref{sumorwtf}) must be large enough to compensate the 
negative $t_\nu$\~order due to $\tg_{\tal^i}\ts^i$. 
Let us address this. 
\smallskip

We can now assume that $m_i\ge 1$. 
It suffices to consider $\tal^r=[\al^r,\nu_{\al^r} m_r]$ with 
$\nu_{\al^r} m_r>0$ and $\tal^r=[\al^r,0]$ with 
$\al^r\not\in \la(s^i)$. The roots $\tal^r$ from
$\la(s^i)$ have been already used; recall that $s^i=s_{\al_i}$.
Let us evaluate the corresponding number of pairs
$\{\tal^r,\tbe^r<0\}$ such that 
$\nu_{\al^r}=\nu=\nu_{\be^r}$ and 
\begin{align}\label{ordwtbe+}
&\ord_w^\nu(\tf_{\tbe^r}) > \ord_w^\nu(\tf_{\tal^r})\for r<i,
\hbox{\  equivalently,}\\
&w^{-1}(\be^r)=w^{-1}(s_{\al^i}(\al^r))>0 \and w^{-1}(\al^r)<0.\notag
\end{align}
Thus $\al^i,\al^r\in \la_\nu(w^{-1})$. We continue to assume
that $\ep'=-$; so $\al'=-w^{-1}(\al^i)>0$.

Let us take an arbitrary $\ga\in \la_\nu(s_{\al'})$; then
$\ga'=-s_{\al'}(\ga)$ belongs to $\la_\nu(s_{\al'})$ as well.
We will first consider the $ADE$\~case. For an arbitrary $j$ such
that  $m_i>j>0$, we claim that either $[-w(\ga),j]$ or 
$[-w(\ga'),m_i-j]$ belongs to the sequence 
$\{\tal^r,\,1\le r \le i\}$. 
Indeed,
\begin{align}\label{vminusga}
&w(-\ga)+w(-\ga')=w(-\al')=\al^i \and\\
&[w(-\ga),j]+[w(-\ga'),m_i-j]=[w(-\al'),m_i]=\tal^i.\notag
\end{align}
Therefore exactly one of these two roots must occur in the sequence
$\la(b)$ before $\tal^i$. Generally speaking this one can be with a 
negative nonaffine component, but all nonaffine components
are positive in $\la(b)$ since $b\in P_+$.

Moreover, we can assume that $-w(\ga)>0$, since 
at least one of $w(-\ga)$ and $w(-\ga')$ must be positive.
Then the claim is that either $w(-\ga)=\tal^r\not\in\la(s^i)$
for certain $r<i$ or $[w(-\ga'),m_i]=\tal^r$ for $r\le i$,
which readily follows from (\ref{vminusga}) with $j=0$.


We obtain that the number of $\tal^r\not\in\la(s^i)$ in the 
form  $[-w(\ga),j]$ for $\ga\in \la(s_{\al'})$ will be 
$m_i((\rho, \al')-1)+m_i=m_i(\rho, \al')$, where 
the second $m_i$ counts the roots 
$[\al^i, j]$ for $1\le j\le m_i$. 

In the $BCFG$\~case, the
calculation is very similar to that for (\ref{wisidm});
the general answer for $\tal^r\not\in \la(s^i)$ satisfying
(\ref{ordwtbe+}) is $m_i(\rho^\vee_\nu, \al')$. 
Indeed, for any root system $R,$
we need to count the number $\mathfrak{N}_\nu$ 
of \,$\tga=[\ga,\nu_\ga k]$ \, such
that 
$$
\ga\in \la_\nu(s_{\al'})\setminus\{\al'\},\ 
0< \nu_\ga k< \nu_{\al^i} m_i,\  [w(\ga),\nu_\ga k]
\in \{\tal^r,\,1\le r \le i\};
$$ 
the consideration of the roots $\tal^r\not\in \la(s^i)$,
with $k=0$ and those for $\nu_\ga k= \nu_{\al^i} m_i$, is left
to the readers. We obtain that
\begin{align*}
\mathfrak{N}_\nu=
&(m_i-1)\bigl((\rho^\vee_\nu,\al)-\de_{\al,\nu}\bigr)
\hbox{\ \ \,for long or short\ }\al=\al'.
\end{align*}
The remaining $\tal^r$ with $k=0$ 
and those for $\nu_\ga k= \nu_{\al^i} m_i$ will change
$(m_i-1)$ in $\mathfrak{N}_\nu$
by $m_i$. Then we add the roots with $\ga=\al'$ if $\nu_\al=\nu$
and obtain the required $m_i(\rho^\vee_\nu, \al')$.
\smallskip

To finalize this calculation we use (\ref{ord-g-talnspr})\,:
\begin{align*}
\dord_w^\nu(\kapp^\de(\tg_{\tal^i}s_{\tal^i}))=
\dord_w^\nu(\kapp^\de(s_{\al^i}\Ga_{\al^i}^{-m_i}
\tg_{-\tal^i}))=
m_i\,(\rho^\vee_\nu,w^{-1}(\al^i));
\end{align*}
recall that $w^{-1}(\al^i)<0$. Therefore (\ref{ordwjump}) 
holds in this case.
\medskip

$(b)$\, {\em The case of\,
$w^{-1}(\al^i)>0$.} Now $\ep'=+ \and (\al^r,\al^i)>0$.
We fix $\nu\in \nu_R$.
We can essentially follow the same verification as for
(\ref{ordwtbe+}) with $\breve{w}=w_0w$ instead of $w$. 
Indeed, $\breve{w}^{-1}(\al^i)<0$ and we need to count the number 
of pairs $\{\tal^r,\tbe^r\equal\ts^i(\tal^r)\}$ satisfying the same
positivity conditions as in  (\ref{ordwtbe+}), namely,
$\breve{w}^{-1}(\al^r)<0$ and $\breve{w}^{-1}(\be^r)>0$. 

However now the switch from 
$\,\ord_w^\nu(\tf_{\tbe^r})\,$ to 
$\,\ord_w^\nu(\tf_{\tal^r})$
in the corresponding pair $\{\tal^r,\tbe^r\}$ will 
increase 
\begin{align}\label{ordtbesum}
\ord_w^\nu\bigl(\tf_{\tbe^1}\tf_{\tbe^2}
\cdots\tf_{\tbe^{i-1}}\bigr)=\ord_w^\nu(\tf_{\tbe^1})
+\ldots+\ord_w^\nu(\tf_{\tbe^{i-1}})
\end{align}
from (\ref{ordwjump}) by one. Accordingly, it suffices to know 
the upper bound for the number of such pairs, not the lower bound 
(actually, the exact number) needed in (\ref{ordwtbe+}). 
Relation (\ref{ordwjump}), which compares the sum in
(\ref{ordtbesum}) plus  $\dord(\tg_{\tal^i}\ts^{i})$ 
with that for $\tal^r$ $(1\le r\le i)$, 
holds (only) due to the positive $t_\nu$\~orders of 
$\tg_{\tal^i}$ and $\ts^{i}=s_{\tal^i}$.
\smallskip

In contrast to (\ref{ordwtbe+}), we now have to include  
nonaffine $\tal^r$ from $\la(s^i)$. This is straight and
will be considered below. We obtain that
the number of pairs $\{\tal^r,\tbe^r\}$ such that the substitution
$\tbe^r\mapsto \tal^r$ increases $(\ref{ordtbesum})$
can be no greater than
$$
(m_i+1)((\rho^\vee_\nu, \al')-\de_{\al',\nu})+
\de_{\al',\nu}m_i=
(m_i+1)(\rho^\vee_\nu, \al')-\de_{\al',\nu};
$$
recall that $\al'=w^{-1}(\al^i)$ and $\nu\in \nu_R$
is fixed.

On the other hand, (\ref{ord-g-talnspr}) results in
\begin{align*}
\dord_w^\nu(\kapp^\de(\tg_{\tal^i}s_{\tal^i}))=
(m_i+1)\,(\rho^\vee_\nu,w^{-1}(\al^i)).
\end{align*}
Thus the change of this order versus 
$\ord_w^\nu(\kapp^\de(\tf_{\tal^i}))=\de_{\al',\nu}$ 
is exactly $(m_i+1)(\rho^\vee_\nu, \al')-
\de_{\al',\nu}$ 
and therefore it
``compensates" the total sum of negative differences 
$\ord_w^\nu(\tf_{\tbe^r}-\tf_{\tal^r})$ for $r<i$.  
We establish that the inequality from (\ref{ordwjump}) 
holds when $w^{-1}(\al^i)>0$ and conclude the justification
of Lemma \ref{ordwtbe}.
\smallskip

\setcounter{equation}{0}
\section{\sc Toda-Dunkl Operators}
The general case will be managed by induction with respect to
the number of $s^i$ taken in the products. We will begin
with some notations and basic estimates.

Recall that for $\tal=[\al,\nu_\al j]\in \tR$, $t_{\tal}=t_{\al}=
t_{\nu_\al}$,
\begin{align*}
&G_{\tal}^+\equal t_{\al}^{1/2}+\frac{t_{\al}^{1/2}-t_{\al}^{-1/2}}
{X_{\tal}^{-1}-1}(1-s_{\tal})=
t_{\al}^{-1/2}(f_{\tal}+g_{\tal}\,s_{\tal}),\\
&\notag \ \ \ \ \, 
f_{\tal}=\frac{t_{\al} X_{\tal}^{-1}-1}{X_{\tal}^{-1}-1}, \ \
g_{\tal}=\frac{t_{\al}-1}{1-X_{\tal}^{-1}};\\
&G_{\tal}^-\equal t_{\al}^{-1/2}+
\frac{t_{\al}^{1/2}-t_{\al}^{-1/2}}{1-X_{\tal}}(1-s_{\tal})
=t_{\al}^{-1/2}(f_{\tal}-s_{\tal}\,g_{\tal}).
\comment{
&\notag \ \ \ \ \, f_{\tal}'=
\frac{t_{\tal}-X_{\tal}^{-1}}{1-X_{\tal}^{-1}}, \ \
g_{\tal}'=\frac{t_{\tal}-1}{X_{\tal}-1}.
}
\end{align*}
\comment{
Recall that
$$
G_{\al_i}^+=s_i\,T_i,\ \
G_{-\al_i}^+=T_i\,s_i,\ \
G_{\al_i}^-=s_i\,T_i^{-1},\and
G_{-\al_i}^-=T_i^{-1}\,s_i.
$$
}
Also, $\ddot{G}_{\tal}^\pm\equal t_{\tal}^{1/2}G_{\tal}^\pm$
and $Y_b=q^{(b_-,\rho_k)}\,b\,
\ddot{G}_{\tal^l}^{\hbox{\tiny\sgn}(\ep_l)}\cdots
\ddot{G}_{\tal^1}^{\hbox{\tiny\sgn}(\ep_1)}$ for $b\in P$;
see~(\ref{ybsgn}).
Note that $\ddot{G}^{-}_{-\tal}=\tilde{G}_{\tal}$ was used
in (\ref{prodforg}). 

See the definition and basic properties of
$\ord_w^\nu$ and $\dord_w^\nu$ in the beginning
of Section~\ref{sec:t-est}.

The following proposition provides the basic tool needed for
a direct proof of the existence of the Toda-Dunkl operators.
The proposition consists of two parts.
In part A, given $u\in W$ and a reduced decomposition
$u=s_{j_l}\cdots s_{j_1}$, we will consider
arbitrary products of the form
\begin{align}
\label{u-prod}
\ddot{G}_{\al^p}^{\pm}\cdots \ddot{G}_{\al^r}^{\pm},\ \  
\ddot{G}_{-\al^r}^{\pm}\cdots \ddot{G}_{-\al^p}^{\pm} 
\for 1\leq r\leq p\leq l,
\end{align}
and expand such products by choosing from each 
$\ddot{G}_{\al}^{\pm}$ either $f_\al$ or $g_\al \, s_\al$. 
The statements in part B are more restrictive, though they
directly result in the existence of the Toda-Dunkl operators;
see Corollary~\ref{YBORDS}.
We will see in Proposition~\ref{SpinDunkl} that it is 
sufficient to manage only the nonaffine products from 
(\ref{u-prod}) in order to prove the existence of
 the Toda-Dunkl operators. 

\subsection{\bf The key step}
The following proposition is the key in our approach.
\begin{proposition}
\label{prop-ord-bound}
A. For $u\in W$ and its reduced decomposition
$u=s_{j_l}\cdots s_{j_1}$, let $\,1\leq r\leq p\leq l$ and $w\in W$.

(i) The $\ord_w^\nu$ of any product in the expansion of
$\kapp^\de(\ddot{G}_{\al^p}^+\cdots\ddot{G}_{\al^r}^+)$
is bounded below by
$\ord_w^\nu(\kapp^\de(f_{\al^p}\cdots f_{\al^r}))$.

(ii) The $\ord_w^\nu$ of any product in the expansion of
$\kapp^\de(\ddot{G}_{-\al^r}^+\cdots\ddot{G}_{-\al^p}^+)$
is bounded below by
$\ord_w^\nu(\kapp^\de(f_{-\al^r}\cdots f_{-\al^p}))$.

(iii) The $\dord_w^\nu$ of any product in the expansion of
$\kapp^\de(\ddot{G}_{\al^p}^-\cdots \ddot{G}_{\al^r}^-)$
is bounded below by
$\ord_w^\nu(\kapp^\de(f_{\al^p}\cdots f_{\al^r}))$.

(iv) The $\dord_w^\nu$ of any product in the expansion of
$\kapp^\de(\ddot{G}_{-\al^r}^-\cdots\ddot{G}_{-\al^p}^-)$
is bounded below by
$\ord_w^\nu(\kapp^\de(f_{-\al^r}\cdots f_{-\al^p}))$.
\smallskip

B. Let $b=\pi_r s_{j_l}\cdots s_{j_1}$ be a reduced
decomposition of $b\in P_+$ and let $1\leq p\leq l$, $w\in W$.

(i) The $\ord_w^\nu$ of any product in the expansion of
$\kapp^\de(\ddot{G}_{\tal^p}^+\cdots\ddot{G}_{\tal^1}^+)$
is bounded below by
$\ord_w^\nu(\kapp^\de(f_{\tal^p}\cdots f_{\tal^1}))$.

(ii) The $\dord_w^\nu$ of any product in the expansion of
$\kapp^\de(\ddot{G}_{-\tal^1}^-\cdots\ddot{G}_{-\tal^p}^-)$
is bounded below by
$\ord_w^\nu(\kapp^\de(f_{-\al^1}
\cdots f_{-\al^p}))$.
\end{proposition}
\begin{corollary}\label{YBORDS}
Let $b\in P_+$. The $\ord_w^\nu$ of
any particular product in the
expansions of  
$$
\kapp^\de(b^{-1}Y_b)\, =\, G_{\tal^l}^+\cdots G_{\tal^1}^+ 
\hbox{\ \ or\ \ } 
\kapp^\de(Y_b^{-1} b)\, =\, G_{\tal^1}^-\cdots G_{\tal^l}^- 
$$ 
is no smaller than $\ord_w^\nu(b)$ and $\ord_w^\nu(b^{-1})$
correspondingly. Therefore, the operators $\hat{Y}_b=\RE^\de(Y_b)$ and
$\hat{Y}'_b=\RE^\de(Y_b^{-1})$ are well defined and invertible;\,
$Y_b Y_b^{-1}=1$ results in $\hat{Y}_b \hat{Y}'_b=1$.
Moreover, the $\RE^\de$\~limits of all products in 
their $G$\~expansions
are well defined.
\end{corollary}

The corollary readily follows from Part B of the proposition.
For instance, use $(ii)$ and calculations
performed in (\ref{ord-greatest}) and (\ref{ybhwc}) 
in the case of $Y_b^{-1}$. 

We provide below a complete proof only for Claim ($i$) from Part A
of the proposition.
Statements ($ii,iii,iv$) and Part B can be proved by similar 
arguments. The justification is based on the following lemma.

\begin{lemma}
\label{lem-ord-est}
A. Let $u\in W$, choose a reduced decomposition
$u=s_{j_l}\cdots s_{j_1}$, and let $1\leq r\leq p\leq l$.
Then for any $p\geq i\geq r$ and $w\in W$, one has
\begin{align}
&\ord_w^\nu(\kapp^\de(g_{\al^i}f_{s^i(\al^{i-1})}
\cdots f_{s^i(\al^r)}))
\geq\ord_w^\nu(\kapp^\de(f_{\al^i}\cdots f_{\al^r})).
\label{u-ord-est}
\end{align}

B. We use the notation $\ts^k=s_{\tal^k}$, where $1\le k\le l$.
For a reduced decomposition
$b=\pi_rs_{j_l}\cdots s_{j_1}\in P_+$ and any $1\leq i\leq l$,
$w\in W$,
\begin{align}
&\ord_w^\nu(\kapp^\de(g_{\tal^i}f_{\ts^i(\tal^{i-1})}
\cdots f_{\ts^i(\tal^1)}\,\ts^i))
\geq\ord_w^\nu(\kapp^\de(f_{\tal^i}\cdots f_{\tal^1})).
\label{u-ord-esta}
\end{align}
\end{lemma}
\smallskip

\subsection{\bf The justifications}$\,$\\
{\em Proof of Lemma \ref{lem-ord-est}.\ }
We prove Part A only; the second part
is quite parallel to Section \ref{sec:ONEGANYW}
instead of the arguments below. Actually,
formula (\ref{ordwtbe+}) above is the only really 
special feature of the affine case. This formula
and related ones were given in Section \ref{sec:ONEGANYW} 
in the case of $Y_b^{-1}$; the adjustments needed 
for (\ref{u-ord-esta}), the case of $Y_b$, are straightforward.
This formula is exactly
the reason why we need to make $r=1$ when extending 
(\ref{u-ord-est}) to the affine case.
\smallskip

Let us begin with some basic orders.
We note first that for $\tal=[\al,\nu_\al j]$
\begin{align}
\label{ord-g-tal}
\ord_w^\nu(\kapp^\de(g_{\tal}))&=
\begin{cases}
0,&\text{if } w^{-1}(\al)>0,\\
-(\rho_\nu^\vee,w^{-1}(\al)),&\text{if } w^{-1}(\al)<0,
\end{cases}
\end{align}
and
\begin{align}
\label{ord-f-tal}
&\ord_w^\nu(\kapp^\de(f_{\tal}))=
\begin{cases}
0,&\text{if } w^{-1}(\al)>0,\\
\de_{\nu,\nu_\al},&\text{if } w^{-1}(\al)<0.
\end{cases}
\end{align}
The second line in (\ref{ord-f-tal}) follows from
the fact that $(\rho_\nu^\vee,w^{-1}(\al))\neq 0$ for all $w\in W$
provided $\nu=\nu_\al$.
\smallskip

Let $u=s_{j_l}\cdots s_{j_1}$ be the reduced decomposition
from Part A and $r\le i\le p\le l$.
Write $\al=\al^i$ (so $s_\al=s^i$) and take $\be=\al^k$ for any
$i> k \geq r$.
Using (\ref{ord-f-tal}), one has
$$ \ord_w^\nu(\kapp^\de(f_{\be}))\leq
\ord_w^\nu(\kapp^\de(f_{s_\al(\be)})) $$
unless
\begin{align}
\nu=\nu_\be, \ \ \ w^{-1}(\be)<0, \and w^{-1}(s_\al(\be))>0.
\label{bad-be}
\end{align}
An equivalent description of (\ref{bad-be}) is
\begin{align}
\ord_w^\nu(\kapp^\de(f_\be))=1
\and \ord_w^\nu(\kapp^\de(f_{s_\al(\be)}))=0,
\label{be-drop}
\end{align}
Assuming (\ref{bad-be}) holds, there are
two cases to consider:
either $w^{-1}(\al)>0$ or $w^{-1}(\al)<0$.

Suppose $w^{-1}(\al)>0$.
Then $\ord_w^\nu(\kapp^\de(g_\al))=\ord_w^\nu(\kapp^\de(f_\al))=0$.
If (\ref{be-drop}) occurs,
then one must have $(\be,\al)<0$.
Hence $s_\al(\be)$ belongs to $\la(u)$ and
by Lemma \ref{la-lem},
one has $s_\al(\be)=\al^j$ where $i>j>k$.
Therefore, the application of $s_\al$ to the product
$f_{\al^{i-1}}\cdots f_{\al^r}$ reverses the positions of the
factors $f_{s_\al(\be)}$ and $f_\be$ for all pairs
$\{\be,s_\al(\be)\}$, where $\be$ satisfies (\ref{be-drop});
the $\ord_w^\nu$ of any other factors in this product
can only increase upon the application of $s_\al$.
This proves (\ref{u-ord-est}) when $w^{-1}(\al)>0$.

It remains to consider 
the case when $w^{-1}(\al)<0$.
We note that
$w^{-1}(s_\al(\be))=s_{w^{-1}(\al)}(w^{-1}(\be))$.
By (\ref{lashtlng1}), one has
\begin{align}
\label{l-ineq}
& l_\nu(s_{w^{-1}(\al)})\leq-2(\rho_\nu^\vee,w^{-1}(\al))-
\de_{\nu,\nu_\al}.
\end{align}
(The only case when (\ref{l-ineq}) is not an equality is
$\nu_\al = \nu_{\lng}$ and $\nu=\nu_{\sht}$.)
Combining this with (\ref{ord-g-tal}) and (\ref{ord-f-tal}) yields
\begin{align}
&\ord_w^\nu(\kapp^\de(g_\al))\geq\ord_w^\nu(\kapp^\de(f_\al))
+\frac{l_\nu(s_{w^{-1}(\al)})-\de_{\nu,\nu_\al}}{2}.
\label{g-corr}
\end{align}
Using (\ref{lashtlng1}), one sees that
$$
\frac{l_\nu(s_{w^{-1}(\al)})-\de_{\nu,\nu_\al}}{2}
$$
is the maximum possible number of
$\be$ satisfying (\ref{bad-be}).
In other words, (\ref{g-corr}) compensates for all drops in the order
coming from (\ref{be-drop}) when applying $s_\al$ to the product
$f_{\al^{i-1}}\cdots f_{\al^r}$. This establishes (\ref{u-ord-est}).
\sq
\smallskip

Now we are going to prove Statement ($i$) from Part A of
Proposition \ref{prop-ord-bound}.
We argue by induction on the number of factors of the
form $g_\al\,s_\al$ chosen to form a particular product in
the expansion --- the base case being the product when no such
factors are chosen, i.e., 
$\p^\emptyset\equal\kapp^\de(f_{\al^p}\cdots f_{\al^r})$.

First, let us consider some particular cases.
Suppose that just one factor of the form $g_\al\,s_\al$, say
$g_{\al^i}\,s^i$, is chosen. In other words, take the product
\begin{align}
\p^i\equal\kapp^\de(f_{\al^p}\cdots f_{\al^{i+1}}\,
g_{\al^i}\,s^i\,f_{\al^{i-1}}\cdots f_{\al^r}).\notag
\end{align}
Due to (\ref{def-ord}),
\begin{align}
\ord_w^\nu(\p^i)
&=\ord_w^\nu(\kapp^\de(f_{\al^p}\cdots f_{\al^{i+1}}\,
g_{\al^i}\,f_{s^i(\al^{i-1})}\cdots f_{s^i(\al^r)}))\notag\\
&=\ord_w^\nu(\kapp^\de(f_{\al^p}\cdots f_{\al^{i+1}}))
+\ord_w^\nu(\kapp^\de(g_{\al^i}\,f_{s_{\al^i}(\al^{i-1})}\cdots
f_{s_{\al^i}(\al^r)})),\notag
\end{align}
Then (\ref{u-ord-est}) gives
$\ord_w^\nu(\p^i)\geq\ord_w^\nu(\p^\emptyset)$,
as claimed.

Next, let us consider the case when two factors of $g_\al\,s_\al$
are chosen:
\begin{align*}
\p^{ij}\equal\kapp^\de(f_{\al^p} \cdots f_{\al^{i+1}}\,
(g_{\al^i}\,s^i)\,f_{\al^{i-1}}
\cdots f_{\al^{j+1}}\,(g_{\al^j}\,s^j)\,
f_{\al^{j-1}}\cdots f_{\al^r}).
\end{align*}
Due to (\ref{def-ord}),
\begin{align}
\ord_w^\nu(\p^{ij})=\ord_w^\nu(\kapp^\de(
f_{\al^p} \cdots f_{\al^{i+1}}\,&g_{\al^i}\,f_{s^i(\al^{i-1})}
\cdots f_{s^i(\al^{j+1})}\notag\\
\times g_{s^i(\al^j)}&\,f_{s^i s^j(\al^{j-1})}\cdots
f_{s^i s^j(\al^r)})).
\label{two-g}
\end{align}
Apply (\ref{ordw-vf}) and (\ref{u-ord-est}) as follows:
\begin{align*}
\ord&_w^\nu(\kapp^\de(g_{s^i(\al^j)}\,f_{s^i s^j(\al^{j-1})}\cdots
f_{s^i s^j(\al^r)}))\\
=&\ord_{s^i w}^\nu(\kapp^\de(g_{\al^j}\,f_{s^j(\al^{j-1})}\cdots
f_{s^j(\al^r)}))\\
\geq\ord_{s^i w}^\nu&(\kapp^\de(f_{\al^j}\,f_{\al^{j-1}}\cdots
f_{\al^r}))\!=\!\ord_w^\nu(\kapp^\de(f_{s^i(\al^j)}\,
f_{s^i(\al^{j-1})}\cdots
f_{s^i(\al^r)})).
\end{align*}
Returning to (\ref{two-g}), one then has
\begin{align}
\label{two-g-one-g}
\ord_w^\nu(\p^{ij})&\geq
\ord_w^\nu(\kapp^\de(f_{\al^p}\cdots f_{\al^{i+1}}\,g_{\al^i}\,
f_{s^i(\al^{i-1})}\cdots f_{s^i(\al^r)}))\\
&=\ord_w^\nu(\p^i)\geq\ord_w^\nu(\p^\emptyset).\notag
\end{align}

In general, for any decreasing sequence 
$p\geq i_1 > i_2 > \cdots > i_m \geq r,$
we set 
$$ \p^{i_1\ldots i_m}\equal
\kapp^\de(h_p\cdots h_r), $$
where $h_i = g_{\al^i}\,s^i$ whenever $i\in\{i_1,\ldots,i_m\}$
and $h_i=f_{\al^i}$ otherwise.
The same reasoning used to arrive
at (\ref{two-g-one-g}) shows that
\begin{align}
\ord_w^\nu(\p^{i_1 \ldots i_m})\geq\ord_w^\nu(\p^{i_1 \ldots i_{m-1}}),
\end{align}
which gives the induction step and completes the proof.
\sq
\medskip

{\em Omitting all $f_\al$.}
As an example clarifying the nature of the estimates
in Lemma \ref{lem-ord-est}, let us discuss the extremal case of 
Proposition \ref{prop-ord-bound} when
we choose all $g_\al\,s_\al$ when expanding
$\kapp^\de(\ddot{G}_{\al^l}^+\cdots\ddot{G}_{\al^1}^+)$.
That is, consider the product
\begin{align}
\kapp^\de(g_{\al^l}\,s^l\cdots g_{\al^1}\,s^1)=
\kapp^\de(g_{-u^{-1}(\al_{j_l})}
\cdots g_{-u^{-1}(\al_{j_1})}\,u^{-1}).\notag
\end{align}
Let $M_j^\nu(u; w)$ denote the number of simple reflections
$s_j$ for $\nu_j=\nu$ in the given decomposition
$u=s_{j_l}\cdots s_{j_1}$ such that $(uw)^{-1}(\al_j)>0$
(so it depends on the choice of the reduced decomposition).
Using (\ref{ord-g-tal}) and (\ref{ord-f-tal}),
Proposition \ref{prop-ord-bound}($i$) then translates to
\begin{align}
\sum_{j=1}^n M_j^\nu(u;w)\,(\al_j,uw(\rho_\nu^\vee))\,
\geq \sum_{\al\in\la(u)\cap\la(w^{-1})}\de_{\nu,\al}.
\notag
\end{align}

Similarly, Proposition \ref{prop-ord-bound}($iii$) 
in such a case leads to the product
\begin{align}
\kapp^\de(s^l\,g_{\al^l}\cdots s^1\,g_{\al^1})=
\kapp^\de(u^{-1} g_{\al_{j_l}}
\cdots g_{\al_{j_1}})\notag
\end{align}
and gives that
\begin{align}\label{nonjwx}
-\sum_{j=1}^n \de_{\nu,\al_j}\,
N_j(u;w)\,(\al_j,w(\rho_\nu^\vee))\,
\geq\, \sum_{\al\in\la(u)\cap\la(w^{-1})}\de_{\nu,\al},
\end{align}
where $N_j(u; w)$ is the number of simple reflections
$s_j$ in $u=s_{j_l}\cdots s_{j_1}$ such that $w^{-1}(\al_j)<0$.
It is a counterpart of formula (\ref{nonjw}) with
$u\in W$ instead of $b\in P_+.$   

When $w=w_0$, 
relation (\ref{nonjwx}) becomes the identity
$l_\nu(u)=l_\nu(u)$.  Taking $w=u^{-1}$ and assuming  
that there is only one simple $\al_j$ in $\la_\nu(u)$, we obtain
the following upper bound for $l_\nu(u)$\,: 
$$
-N_j(u)(u(\al_j),\rho_\nu^\vee)\ge l_\nu(u),\ 
N_j(u)= |\{j_r=j,\, r=1,2,\ldots,l\}|.
$$

\comment{
$$
-(\al_j,u^{-1}\rho_\nu^\vee)=
(\al,\rho_\nu^\vee-u^{-1}(\rho_\nu))-1=
\sum_{\al\in\la(u^{-1})}(\al_j,\al)-1=
-\sum_{\al\in u(\la(u))}(\al_j,\al)-1=
-\sum_{\al\in \la(u)}(\al_j,u(\al))-1.
$$
}
\smallskip

\subsection{\bf \texorpdfstring{$Y$}{\em Y}-hat operators}
In this section, we use the considerations above
to give a direct and constructive
justification of the existence of the Toda-Dunkl operators, based
on the consideration of minuscule weights, $\vth$  and 
short roots (considered as weights).

\begin{proposition}
\label{SpinDunkl}
The limit $\RE^\de(Y_b)$ exists when

$(1)$ $b=\om_r \ \ (r\in O')$,\,\ \ \ \ \ 
$(2)$ $b$ is a short positive root,

$(3)$ $b=-\om_r \ (r\in O')$,${}^{\,}$\ \ \ \
$(4)$ $b$ is a short  negative root,\\
and therefore it exists for any $b\in P$.
\label{DEREY}
\end{proposition}
\proof
The final claim holds because $P$ is generated by the minuscule 
weights together with $Q$, which in turn is generated by 
the short roots.

For (1) and (2), we consider first $b\in P_+$.
Write $b=\pi_r\tilde{w}=\pi_rs_{j_l}
\cdots s_{j_1} \ ( l= l(b))$
in $\hat{W}$ and form $\tilde{\al}^p \ 
(1\leq p\leq  l)$ from (\ref{tal}).
Since $b\in P_+$, one has $l_\nu(b)=2(b,\rho_\nu^\vee)$ 
and $Y_b=T_b=\pi_r T_{j_l}\cdots T_{j_1}$.
Hence
$Y_b=q^{-(b,\,\rho_k)}\Ga_{-b}\,\ddot{G}_{\tal^l}^
+\cdots\ddot{G}_{\tal^1}^+$.
Using (\ref{kapp-de-Ga}), we can write
$$ \kapp^\de(Y_b)=\sum_{w\in W}
q^{-(b,\,\rho_k-w(\rho_k))}\,\Ga_{-w^{-1}(b)}\,\ze_w\,
\kapp^\de(\ddot{G}_{\tal^l}^+\cdots\ddot{G}_{\tal^1}^+).$$

We claim that for any $b\in P_+$,
\begin{align}
\label{xi-de-Y}
\xi^\de(Y_b)\equal \sum_{w\in W}
q^{-(b,\,\rho_k-w(\rho_k))}\,\Ga_{-w^{-1}(b)}\,\ze_w\,
\kapp^\de(f_{\tal^l}\cdots f_{\tal^1})
\end{align}
is regular at $t_\nu=0$.
Indeed, one has
$$ q^{-(b,\,\rho_k-w(\rho_k))}=
\prod_\nu t_\nu^{-(b,\,\rho_\nu^\vee-w(\rho_\nu^\vee))} $$
and the exponents $(b,\rho_\nu^\vee-w(\rho_\nu^\vee))$
count the number of $\tal=[\al,\nu_\al j]\in\la_\nu(b)$
such that $w^{-1}(\al)<0$.
This follows from (\ref{xlambi})
and the following:
\begin{align}
\label{rho-wrho}
\rho_\nu^\vee-w(\rho_\nu^\vee)=
\sum_{\al\in\la_\nu(w^{-1})}\al^\vee.
\end{align}
Now the regularity of $\xi^\de(Y_b)$ is immediate from
(\ref{ord-f-tal}).
\smallskip

(1) Let $b=\om_r$ for $r\in O'$; recall that $\om_r=\pi_r u_r$.
Using Proposition \ref{prop-ord-bound}, where we take $u=u_r$,
the regularity of $\kapp^\de(Y_{\om_r})$ follows from that of
$\xi^\de(Y_{\om_r})$.
\smallskip

\comment{
(2) Let $b=\vth$. We have $\vth=s_0 s_\vth$ and
$ l(\vth)=1+ l(s_\vth)$.
Let $ l= l(s_\vth)$, and write $\al^1,\cdots,\al^l$ for
$\tal^1,\cdots,\tal^l$.
Clearly, $\tal^{l+1}=[\vth,1]$. Hence
\begin{align}
\kapp^\de(Y_\vth)=\Bigl(\sum_{w\in W}
q^{-(\vth,\rho_k-w(\rho_k))}
\Ga_{-w^{-1}(\vth)}\otimes w\Bigr)\kapp^\de
(\ddot{G}_{\tal^{l+1}}^+)\cdots\ &,\\
\kapp^\de(\ddot{G}_{\tal^{l+1}}^+)=
\Bigl(\kapp^\de(f_{[\vth,1]})+\kapp^\de(g_{[\vth,1]})
(\sum_{w\in W} t^{-(w^{-1}(\vth),\rho^\vee)}
\Ga_{w^{-1}(\vth)}\otimes w)\,
s_\vth\Bigr)&.\notag
\end{align}

Choosing $\kapp^\de(f_{[\vth,1]})$ from the binomial
$\kapp^\de(\ddot{G}_{[\vth,1]}^+)$, we apply Proposition
\ref{prop-ord-bound}
(with $u=s_\vth$) to the remaining factors in $\kapp^\de(Y_\vth)$ 
and then use the regularity of $\xi^\de(Y_\vth)$.

On the other hand, choosing 
$\kapp^\de(g_{[\vth,1]}\,s_{[\vth,1]})$, we are led to
consider$s_\vth\,\tilde{t}\,\kapp^\de(\ddot{G}_{\al^l}^+)
\cdots\kapp^\de(\ddot{G}_{\al^1}^+)$,
where $\tilde{t}$ is given by
\begin{align}
\label{th-tw}
\sum_{w\in W} t^{-e_w}\otimes w, \for e_w =
\begin{cases}
(\vth,\rho^\vee),&\text{if } w^{-1}(\vth)<0,\\
(\vth,\rho^\vee-w(\rho^\vee)),&\text{if } w^{-1}(\vth)>0.
\end{cases}
\end{align}
The order of (\ref{th-tw}) is bounded below by that of
\begin{align}
\sum_{w\in W}\prod_{\al\in\la(s_\vth)\cap\la(w^{-1})}t^{-1}
\otimes w.
\end{align}
Thus we can again apply Proposition \ref{prop-ord-bound} with
$u=s_\vth$.
\smallskip
}

(2) Suppose $b=\al$ is any short positive root.
Using Lemma \ref{lem-sht-pos}, find a reduced expression
$s_\vth=s_{j_1}\cdots s_{j_p}s_ms_{j_p}\cdots s_{j_1}$ such that
$\al=s_{j_r}\cdots s_{j_1}(\vth)$ where $0\leq r\leq p$.
Let $l=l(s_\vth)=2p+1$ and construct
$\la(s_\al)=\{\al^1,\ldots,\al^l\}$ using the chosen
reduced decomposition.

Recall that $\vth=s_0 s_\vth$ and $l(\vth)=l(s_\vth)+1$.
Accordingly, one has $Y_\vth=T_0T_{s_\vth}$ and
$\la(\vth)=\la(s_\vth)\cup\{[\vth,1]\}$.

Due to (\ref{inn-prod-min1}) and (\ref{TYTL}), one has
\begin{align*}
Y_\al&=(T_{j_r}^{-1}\cdots T_{j_1}^{-1})\,T_0\,(T_{j_1}
\cdots T_{j_p}
T_m T_{j_p}\cdots T_{j_{r+1}}).
\end{align*}
Hence, for $v=s_{j_r}\cdots s_{j_1}$, we can write
\begin{align}
\label{v-Y-al}
&v^{-1}\,Y_\al\,v=q^{-(\vth,\,\rho_k)}\,\ddot{G}_{\al^r}^
-\cdots\ddot{G}_{\al^1}^-\,
\Ga_{-\vth}\,\ddot{G}_{[\vth,1]}^+\,\ddot{G}_{\al^l}^
+\cdots\ddot{G}_{\al^{r+1}}^+.
\end{align}
We note that by (\ref{RE-W-eq})
$$ \kapp^\de(v^{-1}\,Y_\al\,v)=v^{-1}\,\kapp^\de(Y_\al)\,v. $$
Hence it suffices to prove that $\kapp^\de(v^{-1}\,Y_\al\,v)$
is regular at $t_\nu=0$.

By Proposition \ref{prop-ord-bound}($i,iii$), it is enough to consider
\begin{align*}
q^{-(\vth,\,\rho_k)}
\kapp^\de(f_{\al^r}\cdots f_{\al^1}\,\Ga_{-\vth}\,
\ddot{G}^+_{[\vth,1]}\,f_{\al^{l}}\cdots f_{\al^{r+1}}).
\end{align*}
We expand this product by choosing either $f_{[\vth,1]}$ or
$g_{[\vth,1]}\,s_{[\vth,1]}$ from $\ddot{G}_{[\vth,1]}^+$.

Choosing $f_{[\vth,1]}$ from $\ddot{G}_{[\vth,1]}^+$, we arrive at
$\xi^\de(Y_\vth)$, which is known to be regular at $t_\nu=0$
(cf. (\ref{xi-de-Y})).

Thus it remains to choose $g_{[\vth,1]}\,s_{[\vth,1]}$. This yields
$$ q^{-(\vth,\,\rho_k)}
\kapp^\de(f_{\al^r}\cdots f_{\al^1}\,
g_{[\vth,-1]}\,s_\vth\,f_{\al^{l}}\cdots f_{\al^{r+1}}), $$
where we have used that 
$\Ga_{-\vth}\,g_{[\vth,1]}\,s_{[\vth,1]}=g_{[\vth,-1]}\,s_\vth$.
According to (\ref{def-ord}), when calculating $\ord_w^\nu$,
one must move $s_\vth$ to the right:
$$ s_\vth\,(f_{\al^l} \cdots f_{\al^{r+1}})
= (f_{-\al^1}\cdots f_{-\al^{l-r}})\,s_\vth,
$$
where we have used (\ref{tbe-opp}).
By (\ref{ord-g-tal}), we need to show that
\begin{align}
\label{gth-ord}
\ord_w^\nu(&\kapp^\de(f_{\al^r}\cdots f_{\al^1} f_{-\al^1}\notag
\cdots f_{-\al^{l-r}}))\\
&\geq
\begin{cases}
(\vth,\rho_\nu^\vee),&\text{if } w^{-1}(\vth)>0,\\
(\vth,\rho_\nu^\vee+w(\rho_\nu^\vee)),&\text{if } w^{-1}(\vth)<0,
\end{cases}
\end{align}
for any $0\leq r \leq p$.

To this end, assume first that $w^{-1}(\vth)>0$.
Clearly we have
$$ \ord_w^\nu(\kapp^\de(f_{\al^r}\cdots f_{\al^1}
f_{-\al^1}\cdots f_{-\al^r}))
=\sum_{i=1}^r \de_{\nu,\al^i}. $$
For the remaining factors in the left-hand side
of (\ref{gth-ord}), one has
$$ \ord_w^\nu(\kapp^\de(f_{-\al^{r+1}}\cdots f_{-\al^{l-r}}))
\geq\de_{\nu,\vth}+\sum_{i=r+1}^p \de_{\nu,\al^i}. $$
This can be seen as follows.
First, $\al^{p+1}=\vth$ and hence by (\ref{ord-f-tal})
$\ord_w^\nu(\kapp^\de(f_{-\al^{p+1}}))=\de_{\nu,\vth}$.
Second, for each $r+1\leq i\leq p$, at least one of $w^{-1}(\al^i)$ or
$w^{-1}(\al^{l-i+1})$ must be positive. 
This follows from Lemma~\ref{la-lem}($ii$), (\ref{tbe-opp}),
and the assumption $w^{-1}(\vth)>0$.
Therefore, altogether one has
\begin{align}
\ord_w^\nu(\kapp^\de(f_{\al^r}\cdots f_{\al^1} f_{-\al^1}
\cdots f_{-\al^{l-r}}))
\geq \de_{\nu,\vth}+\sum_{i=1}^p \de_{\nu,\al^i}.
\label{gth-ord+}
\end{align}
Finally, using Lemma \ref{LASTAL}($i$), one finds that
the right-hand side of (\ref{gth-ord+}) is exactly
$(\vth,\rho_\nu^\vee)$.

Now assume $w^{-1}(\vth)<0$.
One has
\begin{align}
\ord_w^\nu(\kapp^\de(f_{-\al^1}\cdots f_{-\al^{l-r}}))
=\sum_{\substack{1\leq i\leq l-r\\w^{-1}(\al^i)>0}}\de_{\nu,\al^i}
\label{ordf-1lr}
\end{align}
and
\begin{align}
\ord_w^\nu(\kapp^\de(f_{\al^r}\cdots f_{\al^1}))
&=\sum_{\substack{1\leq i\leq r\\w^{-1}(\al^i)<0}}\de_{\nu,\al^i}
\geq \sum_{\substack{l-r+1 \leq i\leq l\\w^{-1}(\al^i)>0}}
\de_{\nu,\al^i}.
\label{ordfr1}
\end{align}
The inequality in (\ref{ordfr1}) follows from 
Lemma~\ref{la-lem}($ii$), (\ref{tbe-opp}),
and the assumption that $w^{-1}(\vth)<0$.
In particular, if 
$w^{-1}(\al^{l-i+1})>0$ for $1\leq i\leq p$, then necessarily
$w^{-1}(\al^i)<0$.
Putting (\ref{ordf-1lr}) and (\ref{ordfr1}) together, one has
\begin{align}
\label{gth-ord-}
\ord_w^\nu(\kapp^\de(f_{\al^r}\cdots 
f_{\al^1} f_{-\al^1}\cdots f_{-\al^{l-r}}))
\geq
\sum_{\substack{1\leq i\leq l\\w^{-1}(\al^i)>0}}\de_{\nu,\al^i}.
\end{align}
Finally, to get (\ref{gth-ord}), we observe that
\begin{align}
\label{rho+wrho}
\rho_\nu^\vee+w(\rho_\nu^\vee)=
\sum_{\substack{\al>0,\,\nu_\al=\nu\\w^{-1}(\al)>0}}\al^\vee
\end{align}
and consequently $(\vth,\rho_\nu^\vee+w(\rho_\nu^\vee))$
is exactly the number of $\tal\in\la_\nu(\vth)$ with
$w^{-1}(\al)>0$. Note that since $w^{-1}(\vth)<0$, such
$\tal$ must belong to $\la_\nu(s_\vth)\setminus\{\vth\}$.

This completes the proof of (\ref{gth-ord}) and hence
the proof of (2) as well.
\smallskip

\rmk
The relation $T_i^{-1} Y_b T_i^{-1}=Y_{s_i(b)}$ from (\ref{TYTL}),
which was used at the beginning of (2),
is valid only when $(b,\al_i^\vee)=1$. In particular, 
it does not hold for $b=\al_m$ and $i=m$.
In this case,
\begin{align}
\label{T-Y-al}
T_m^{-1} Y_{\al_m} T_m^{-1}=Y_{\al_m}^{-1}+
(t_m^{1/2}-t_m^{-1/2})T_m^{-1}.
\end{align} 
One cannot use (\ref{T-Y-al}) to pass from $Y_{\al_m}$ to 
$Y_{\al_m}^{-1}$ in a way that is compatible with the limit
$t_\nu\to 0$. Nevertheless, we can reach $Y_{\al_m}^{-1}$, along
with all the operators corresponding to negative short roots,
by starting from $Y_\vth^{-1}$.
This is carried out in (4) below.
\sq
\smallskip

Before (3) and (4), let us make some general remarks
about $Y_{-b}$ for arbitrary $b\in P_+$.
Write $b=\pi_r s_{j_l}\cdots s_{j_1} \ (l=l(b))$
and construct $\la(b)=\{\tal^1,\cdots,\tal^l\}$ using
this reduced decomposition.
Then $-b=\pi_r^{-1} s_{\pi_r(j_1)}\cdots s_{\pi_r(j_l)}$,
which is a reduced decomposition.

For $1\leq p\leq l$, let
$$ \tbe^p=-b(\tal^{l-p+1})=
s_{\pi_r(j_l)}\cdots s_{\pi_r(j_{l-p+2})}
(\al_{\pi_r(j_{l-p+1})}), $$
so that $\la(-b)=\{\tbe^1,\cdots,\tbe^l\}$.
We can write
$$ Y_{-b}=q^{-(b,\rho_k)}\Ga_b\,
\ddot{G}_{\tbe^l}^-\cdots\ddot{G}_{\tbe^1}^-
=\ddot{G}_{-\tal^1}^-\cdots\ddot{G}_{-\tal^l}^-\,
q^{-(b,\rho_k)}\,\Ga_b.
$$
Hence
\begin{align}
\kapp^\de(Y_{-b})=
\kapp^\de(\ddot{G}_{-\tal^1}^-\cdots\ddot{G}_{-\tal^l}^-)
\sum_{w\in W}q^{-(b,\,\rho_k+w(\rho_k))}\,\Ga_{w^{-1}(b)}\,\ze_w.
\end{align}
We claim that for any $b\in P_+$,
\begin{align}
\xi^\de(Y_{-b})\equal
\kapp^\de(f_{-\tal^1}\cdots f_{-\tal^l})
\sum_{w\in W}q^{-(b,\,\rho_k+w(\rho_k))}\,
\Ga_{w^{-1}(b)}\,\ze_w.
\end{align}
is regular at $t_\nu=0$. The proof is similar to that for
$\xi^\de(Y_b)$ from (\ref{xi-de-Y}) that was given before step (1).
One uses (\ref{rho+wrho}) instead of (\ref{rho-wrho}).
\smallskip

(3) In the case of $b=\om_r \ (r\in O')$,
the regularity of $\kapp^\de(Y_{-b})$ is immediate
from that of $\xi^\de(Y_{-b})$, due to
Proposition \ref{prop-ord-bound}($ii$).
\smallskip

\comment{
Take $b=\vth$. We have $\tal^{l}=[\vth, 1]$.
The argument now proceeds as in step (2).
In particular, we consider
\begin{align}
\kapp^\de\Bigl(\frac{t-1}{q^{-1}X_\vth^{-1}-1}\,
\Ga_\vth\,s_\vth\Bigr)
\sum_w t^{-(\vth,\rho^\vee-w(\rho^\vee))}\otimes w
\end{align}
After moving $s_\vth$ to the right, this has the same order as
\begin{align}
\label{thmin-tw}
\sum_{w\in W} t^{-e_w}\otimes w,\for e_w=
\begin{cases}
(\vth,\rho^\vee),&\text{if } w^{-1}(\vth)>0,\\
(\vth,\rho^\vee+w(\rho^\vee)),&\text{if } w^{-1}(\vth)<0.
\end{cases}
\end{align}
Note that (\ref{thmin-tw}) has order bounded below by that of
\begin{align}
\sum_{w\in W}\prod_{\al\in\la(s_\vth)
\cap(R_+\setminus\la(w^{-1}))}t^{-1}\otimes w,
\end{align}
and now (\ref{u-prod-pr}) applies with $u=s_\vth$.
}

(4) For $b$ equal to any negative short root $\al$, the proof
is similar to (2). Use Lemma \ref{lem-sht-pos} to choose a 
reduced decomposition
$s_\vth=s_{j_1}\cdots s_{j_p}s_m s_{j_p}\cdots s_{j_1}$
such that
$s_{j_r}\cdots s_{j_1}(-\vth)=-\al$ where $0\leq r\leq p$.
Then, starting from
$Y_\vth^{-1}=T_{s_\vth}^{-1}T_0^{-1}$,
we use (\ref{inn-prod-min1}) and (\ref{TYTL}) to get
\begin{align}
&Y_{s_{j_1}(\vth)}^{-1}=T_{j_1}Y_\vth^{-1}T_{j_1},\ \notag
Y_{s_{j_2}s_{j_1}(\vth)}^{-1}=T_{j_2}T_{j_1}
Y_\vth^{-1}T_{j_1}T_{j_2}, \ \ldots \ ,\\
&Y_\al^{-1}=T_{j_r}\cdots T_{j_1} 
Y_\vth^{-1} T_{j_1}\cdots T_{j_r}.\notag
\end{align}
Then, as in (2), 
the regularity of $\kapp^\de(Y_{\al}^{-1})$ can be shown
using Proposition \ref{prop-ord-bound}($ii,iv$).
\sq
\smallskip


\comment{
\begin{proposition}
The $\de$\~$\RE$ limit of $Y_b$ exists for (i) $b=\om_r \ (r\in O')$
and (ii) $b=\vth$.
In these cases,
\begin{align}
\RE^\de(Y_b)=\xi^\de(Y_b)+
\sum_{v\in W,\,b_\ast}f_{b_\ast,v}\,\Ga_{-b_\ast}^\varrho v, 
\ \ b_\ast =
\begin{cases}
\om_r &\text{for } (i),\\
\vth, 0 &\text{for } (ii),
\end{cases}
\label{REdeYr}
\end{align}
where $f_{b_\ast,v}\in\spin(\overline{\v})$ and
\begin{align}
\xi^\de(Y_{b})\equal\sum_{w\in W} \Ga_{-w^{-1}(b)}
\prod_{\substack{\tal\in\la(b) \
\mathrm{ s.t.}\\(w^{-1}(\tal),\rho^\vee)=-1}}
(1-X_{w^{-1}(\tal)})\otimes w.
\label{xideY}
\end{align}
\label{RE-dunkl}
\end{proposition}
}

\subsection{\bf Other generators}
Recall the definition of $\dHH^{\flat,\vph}$ from
Section~\ref{sec:duallimit} (we continue to take $B=P$).
We will calculate the $\RE^\de$ limits of the
remaining generators $\breve{T}_i\ (i\geq 0)$ and $\breve{\Pi}$.
Recall that $\breve{T}_i=\ddot{T}_i$ for $i>0$ and
$\breve{T}_0=\vph(\ddot{T}_0)=t_0^{1/2}T_{s_\vth}^{-1}T_0^{-1}$.
We also consider $\tilde{X}_b\equal q^{(b_+,\rho_k)}X_b=
\vph(\ddot{Y}_{-b})$.

\begin{proposition}
\label{OtherGens}
(i) The operators $\hat{T}_i=\RE^\de(\breve{T}_i)$ exist for all 
$i=0,\ldots, n$. Moreover,
\begin{align}
\hat{T}_i=\RE^\de(\ddot{T}_i)=
\sum_{\substack{w\in W\,\mathrm{s.t.}\\w^{-1}(\al_i)<0}}
\ze_w\,(s_i-1) \for i>0.
\label{RE-de-T}
\end{align}

(ii) For any $b\in P$,
\begin{align}
\RE^\de(\tilde{X}_b)=
\sum_{\substack{w\in W \ \mathrm{s.t.}\\w^{-1}(b)=b_+}}
X_{b_+}\,\ze_w.
\label{RE-de-dX}
\end{align}

(iii) For any $r \in O'$, 
\begin{align}
\RE^\de(\breve{\pi}_r^{-1})=
&\sum_{w\in W}\Bigl(X_{w^{-1}(\om_r)}\notag
\!\!\prod_{\substack{\al\in\la(u_r) \ \mathrm{s.t.}\\
(w^{-1}(\al),\rho^\vee)=1}}(1-X_{w^{-1}(\al)}^{-1})\,\ze_w\Bigr)\, 
u_r^{-1}\\
&+\sum_{v<u_r^{-1}}\bigl(\sum_{w\in W}f_{v,w}\,\ze_w\bigr)\,v 
\hbox{\ \ for certain\ } f_{v,w}\in \Q_q'[X_b , b\in B].
\label{RE-de-pi}
\end{align}
\end{proposition}
\smallskip

\proof
($i$) Using $\kapp^\de(s_i)=s_i$ and
(\ref{kapp-de-X}), we readily arrive
at (\ref{RE-de-T}) for $i>0$. The case of $i=0$
is significantly more involved. We have
$\breve{T}_0=\vph(\ddot{T}_0)=
t_0^{1/2}T_{s_\vth}^{-1}X_\vth^{-1}$,
where $\vph$ is the duality anti-involution defined in
(\ref{phianti}).
Write $s_\vth=s_{j_l}\cdots s_{j_1}=
s_{j_1}\cdots s_{j_l} \ (l=l(s_\vth))$.
Let $\al^p=s_{j_1}\cdots s_{j_{p-1}}(\al_{j_p})\in\la(s_\vth)$ for
$p=1,\ldots, l$.

Now
\begin{align}
\breve{T}_0=t_0^{1/2}\prod_\nu t_\nu^{-l_\nu(s_\vth)/2}\,
\ddot{G}_{-\al^1}^-\cdots\ddot{G}_{-\al^l}^-\,
s_\vth\,X_\vth^{-1}.
\label{breve-T0}
\end{align}
By Lemma \ref{LASTAL}($i$),
one has $l_\nu(s_\vth)=2(\vth,\rho_\nu^\vee)-\de_{\nu,\vth}$.
Hence
$$ t_0^{1/2}\prod_\nu t_\nu^{-l_\nu(s_\vth)/2}=
\prod_\nu t_\nu^{-(\vth,\,\rho_\nu^\vee)+\de_{\nu,\vth}}. $$
Returning to (\ref{breve-T0}), we have
\begin{align}
\kapp^\de(\breve{T}_0)=
\kapp^\de(\ddot{G}_{-\al^1}^- \cdots\ddot{G}_{-\al^l}^-)
\sum_{w\in W}t_{\sht}\,q^{-(\vth,\,\rho_k+w(\rho_k))}
\,X_{w^{-1}(\vth)}\,\ze_w
\,s_\vth.
\label{kapp-de-brT0}
\end{align}
By Proposition \ref{prop-ord-bound}($iv$),
\begin{align*}
&\dord_w^\nu(\kapp^\de(\ddot{G}_{-\al^1}^- 
\cdots\ddot{G}_{-\al^l}^-)\\
\geq\, \ord_w^\nu&(\kapp^\de(f_{-\al^1}\cdots f_{-\al^l}))
=\!\!\sum_{\al\in\la(s_\vth)\cap(R_+
\setminus\la(w^{-1}))}\de_{\al,\nu}.
\notag
\end{align*}
The claim now follows from (\ref{rho+wrho}) and the description
of the sets $\la_\nu(\vth)$ due to $\vth=s_{0}s_{\vth}$. 
The $t_{\sht}$ factor in 
(\ref{kapp-de-brT0})
accounts for the case when $w^{-1}(\vth)>0$, because
$\la(s_\vth)=\la(\vth)\setminus\{[\vth,1]\}$.
\smallskip

($ii$) By definition, $\tilde{X}_b=q^{(b_+,\,\rho_k)}X_b$; 
hence
$$ \kapp^\de(\tilde{X}_b)=
\sum_{w\in W}q^{(b_+-w^{-1}(b),\,\rho_k)}
X_{w^{-1}(b)}\,\ze_w. $$
Now (\ref{RE-de-dX}) follows, using the fact that
$b_+ \geq w^{-1}(b)$ for all $w\in W$.
\smallskip

($iii$) Recall that $\breve{\pi}_r^{-1}=X_{\om_r}T_{u_r^{-1}}$.
Let $u_r=s_{j_l}\cdots s_{j_1}$ be a reduced decomposition.
Construct $\la(u_r)=\{\al^1,\ldots,\al^l\}$ using this
decomposition. Then
$$
\breve{\pi}_r^{-1}=q^{-(\om_r,\,\rho_k)}X_{\om_r}
\ddot{G}_{-\al^1}^+\cdots\ddot{G}_{-\al^l}^+\,u_r^{-1}.
$$
We have used here that
$l_\nu(u_r)=l_\nu(\om_r)=2(\om_r,\rho_\nu^\vee)$.
Hence
$$
\kapp^\de(\breve{\pi}_r^{-1})=
\bigl(\sum_{w\in W}q^{-(\om_r,\,\rho_k+w(\rho_k))}\,X_{w^{-1}(\om_r)}
\,\ze_w\bigr)\,
\kapp^\de(\ddot{G}_{-\al^1}^+\cdots\ddot{G}_{-\al^l}^+)\,u_r^{-1}.
$$
Now, by (\ref{rho+wrho}) and Proposition \ref{prop-ord-bound}($ii$),
the limit $\RE^\de(\breve{\pi}_r^{-1})$ exists. Then
(\ref{RE-de-pi}) follows readily.
\sq
\smallskip

Propositions \ref{SpinDunkl} and \ref{OtherGens} provide
a direct justification, independent of Theorem \ref{GLOBNSWH},
of one of the key results of this paper:
the action of $\bHH^{\flat,\vph}$, the limit of $\dHH^{\flat,\vph}$ 
as $t_\nu\to 0$, in $\spin(\Q_q'(X))$.
Moreover, we obtain that the generators of $\bHH^{\flat,\vph}$
have no nontrivial denominators and therefore preserve
$\spin(\overline{\v})$, as in Part $(i)$ of Theorem \ref{GLOBNSWH}.
Finally, we also see that it is not necessary to work over the field
$\Q_q'=\Q(q^{1/(2m)})$; the action of the generators of 
$\bHH^{\flat,\vph}$, including that of the Toda-Dunkl operators, 
is defined over the ring $\Z[q^{\pm 1/(2m)}]$.
\smallskip

{\em Symmetrization.}
Let $\spin^\de(\overline{\v})$
denote the space of $W$\~invariants of $\spin(\overline{\v})$
under the $\de$\~action, which is simply
$\de(\overline{\v})$. 
Recall that $\RE$ (without the super index $\de$) is the
non-spinor Ruijsenaars-Etingof procedure defined in
(\ref{kapprho}). We will use the operators $\cL_f$ and 
$L_f$ from (\ref{def-Lf}), where $f\in \Q_q'[X]^W$.

\begin{proposition}\label{SymDun}
Upon the restriction to $\spin^\de(\overline{\v})$,
one has
$$
\RE^\de(\cL_f)=\sum_{w\in W}\RE(L_f)\,\ze_w,
$$
for any $f\in\Q_q'[X]^W$.
\end{proposition}
\proof
First, $f(Y)$ is central in $\dHH_Y^\vph$, the subalgebra
of $\dHH^\vph$ generated by $\breve{T}_i=T_i\ (i>0)$ and 
$Y_b\ (b\in P)$,
and hence
$\RE^\de(\cL_f)$ commutes with $\RE^\de(\ddot{T}_i)$ from
(\ref{RE-de-T}) for $i=1,\ldots,n$. Second,
an element $g\in\spin(\overline{\v})$ belongs to
$\spin^\de(\overline{\v})$ if and only if 
$\hat{T}_i(g)=0$ for $i=1,\ldots,n$.
Indeed, applying (\ref{RE-de-T}) to 
$g=\sum_{w\in W}g_w\,\ze_w$ gives
$$
\hat{T}_i(g)=
\sum_{\substack{w\in W\,\mathrm{s.t}\\
w^{-1}(\al_i)<0}}(g_{s_i w}-g_w)\,\ze_w.
$$
The right-hand side vanishes if and only if
$g_{s_i w}=g_w$ whenever one has $w^{-1}(\al_i)<0$.
The latter condition is always met either by 
$w$ or by $w'=s_iw$. Thus all differences
$g_{s_i w}-g_w$ must vanish.

We conclude that $\RE^\de(\cL_f)$ preserves $\spin^\de(\v)$.
Therefore it has the form $\sum_{w\in W}M\,\ze_w$
upon the restriction to $\spin^\de(\overline{\v})$
for some difference operator $M$.
By considering the $\id$\~component of $\RE^\de(\cL_f)$,
one sees that $M=\RE(L_f)$.
\sq

\subsection{\bf Examples}
\label{examples}
{\em Toda-Dunkl operators.}
For the root system $A_1$,
\begin{align*}
\hat{Y}_{\om}&=\Ga_{-\om}^\varrho\Bigl(
(\ze_{\id}+(1-X_{\al}^{-1})\ze_{s})\id
+(-\ze_{\id}+X_{\al}^{-1}\ze_{s})\,s\Bigr),\\
\hat{Y}_\om^{-1}&=\hat{Y}_{-\om}=
\Bigl((1-X_{\al}^{-1})\ze_{\id}+\ze_s\Bigr)\Ga_{\om}^\varrho
+\Bigl(\ze_{\id}-X_\al^{-1}\ze_s\Bigr)\Ga_{-\om}^\varrho\,s,
\end{align*}
in terms of the fundamental weight $\om$ and simple root $\al$,
where $s=s_\al$ and $\Ga_{-\om}^\varrho$ is from (\ref{rhoGab}).
Upon the restriction to $\spin^\de(\overline{\v})$,
one has
$$\hat{Y}_{\om}+\hat{Y}_{\om}^{-1}\ =\ 
\Bigl((1-X_\al^{-1})\Ga_\om+\Ga_{-\om}\Bigr),
$$
a special case of Proposition~\ref{SymDun}.


For the root system $A_2$, one has
\begin{align*}
\hat{Y}_{\om_1}=
\Ga_{-\om_1}^\varrho\Bigl(
(\ze_{\id}
+(1-X_{\al_1}^{-1})\ze_{s_1}
+\ze_{s_2}
+(1-X_{\al_2}^{-1})\ze_{s_1 s_2}\hspace{.5in}&\\
+(1-X_{\al_1}^{-1})\ze_{s_2 s_1}
+(1-X_{\al_2}^{-1})\ze_{s_1 s_2 s_1})&\id\\
+(-\ze_{\id}
+X_{\al_1}^{-1}\ze_{s_1}
-\ze_{s_2}
-(1-X_{\al_1}^{-1})\ze_{s_2 s_1}
+X_{\al_2}^{-1}\ze_{s_1 s_2 s_1})&\,s_1\\
+(\ze_{s_2}+X_{\al_1+\al_2}^{-1}\ze_{s_1 s_2}
-X_{\al_1}^{-1}\ze_{s_2 s_1}
-X_{\al_1+\al_2}^{-1}\ze_{s_1 s_2 s_1})&\,s_1 s_2\\
+(-\ze_{s_1}
-\ze_{s_2}
+(1-X_{\al_1}^{-1})X_{\al_2}^{-1}\ze_{s_1 s_2}
+X_{\al_1}^{-1}\ze_{s_2 s_1}
+X_{\al_1+\al_2}^{-1}\ze_{s_1 s_2 s_1})&\,s_1 s_2 s_1
\Bigr).
\end{align*}
The operator $\hat{Y}_{\om_2}$ is obtained by interchanging the
indices 1 and 2 of $\om_i$, $s_i$, and $\al_i$ in the above
formula. The operators $\hat{Y}_{\om_i}$ are invertible (their
inverses are $\hat{Y}_{-\om_i}$), and one has
\begin{align*}
\hat{Y}_{\om_1}^{-1}+\hat{Y}_{\om_1}\hat{Y}_{\om_2}^{-1}+\hat{Y}_{\om_2}
=(1-X_{\al_1}^{-1})\Ga_{\om_1}+(1-X_{\al_2}^{-1})\Ga_{-\om_1+\om_2}+
\Ga_{-\om_2}
\end{align*}
upon the restriction to $\spin^\de(\overline{\v})$, where
the right-hand side is the $q$\~Toda operator
$\RE(L_{-\om_1})$; 
cf. Proposition~\ref{SymDun} and (\ref{dLa}).

For the root system $B_2$, with $\al_1$ long and $\al_2$ short,
the fundamental weight $\om_2$ is minuscule, while $\om_1=\vth$ is
not. One has
\begin{align*}
\hat{Y}_{\om_2}=
\Ga_{-\om_2}^\varrho\Bigl(
(\ze_{\id}+\ze_{s_1}
+(1-X_{\al_1}^{-1})(\ze_{s_2 s_1}+\ze_{s_1 s_2 s_1})\hspace{.5in}&\\
+(1-X_{\al_2}^{-1})(\ze_{s_2}+\ze_{s_1 s_2}
+\ze_{s_2 s_1 s_2}+\ze_{s_1 s_2 s_1 s_2}))&\id\\
+(-(\ze_{\id}+\ze_{s_1})
+X_{\al_2}^{-1}(\ze_{s_2}+\ze_{s_1 s_2 s_1 s_2})
-(1-X_{\al_2}^{-1})\ze_{s_1 s_2}
\hspace{.5in}&\\
-(1-X_{\al_1}^{-1})\ze_{s_1 s_2 s_1})&\,s_2\\
+(\ze_{s_1}+(1-X_{\al_2}^{-1})\ze_{s_1 s_2}
+X_{\al_1+\al_2}^{-1}\ze_{s_2 s_1}
-X_{\al_1}^{-1}\ze_{s_1 s_2 s_1}
\hspace{.5in}&\\
+X_{\al_1+\al_2}^{-1}(1-X_{\al_2}^{-1})\ze_{s_2 s_1 s_2}
-X_{\al_1+\al_2}^{-1}\ze_{s_1 s_2 s_1 s_2})&\,s_2 s_1\\
+(-(\ze_{s_1}+\ze_{s_2})-(1-X_{\al_2}^{-1})\ze_{s_1 s_2}
+X_{\al_1}^{-1}(1-X_{\al_2}^{-1})\ze_{s_2 s_1}
\hspace{.5in}&\\
+X_{\al_1+2\al_2}^{-1}\ze_{s_2 s_1 s_2}
+X_{\al_1}^{-1}\ze_{s_1 s_2 s_1})&\,s_2 s_1 s_2\\
+(-(\ze_{s_1}+\ze_{s_2 s_1})+X_{\al_2}^{-1}\ze_{s_1 s_2}
+X_{\al_2}^{-1}(1-X_{\al_1}^{-1})\ze_{s_2 s_1 s_2}
\hspace{.5in}&\\
+X_{\al_1+\al_2}^{-1}\ze_{s_1 s_2 s_1}
-X_{\al_1+\al_2}^{-1}(1-X_{\al_2}^{-1})\ze_{s_1 s_2 s_1 s_2})
&\,s_1 s_2 s_1\\
+(\ze_{s_1}-X_{\al_2}^{-1}\ze_{s_1 s_2}-X_{\al_1+
\al_2}^{-1}\ze_{s_1 s_2 s_1}
+X_{\al_1+2\al_2}^{-1}\ze_{s_1 s_2 s_1 s_2})&\,s_1 s_2 s_1 s_2
\Bigr).
\end{align*}
\smallskip

{\em Application to the nonsymmetric Whittaker function.}
Let us give the values of the coefficients $a_{b,w}$ from
(\ref{nbce}) and 
Proposition~\ref{OMLIMIT} in the case of the root system $A_2$.
These coefficients are the only ingredient of the theory
of $\overline{E}^\dag$\~polynomials necessary for an explicit 
description of the nonsymmetric Whittaker function $\Om$.
See \cite{ChM} or \cite[(2.7)]{ChO2} for the $A_1$\~case.
\smallskip

Note that $a_{0,w}=1$ for all $w\in W$. The following tables
give the values of $a_{b,w}$ for nonzero $b$ with a fixed $b_-$
and all $6$ elements $w\in\S_3$.
\medskip

For $b_-=n\om_2 \ (n<0)$:

\begin{tabular}{|c||c|c|c|c|c|c|}
\hline
$b\setminus w$ & $\id$ & $s_1$ & $s_2$ & $s_2 s_1$ & $s_1 s_2$ & $s_1 s_2 s_1$ \\
\hline\hline
$b\ =\ n\om_2$ & $1$ & $1$ & $q^{n}$ & $q^{n}$ & $q^{n}$ & $q^{n}$
\\ \hline
$b\!=\!s_2(b_-)$ & $0$ & $0$ & $1$ & $1$ & $q^{n}$ & $q^{n}$\\ \hline
$b=-n\om_1$ & $0$ & $0$ & $0$ & $0$ & $1$ & $1$\\ \hline
\end{tabular}    
\medskip

\noindent
We put $a_{b,w}$ in the corresponding row and column above and in 
the following tables.
\medskip

For $b_-=n\om_1 \ (n<0)$:

\begin{tabular}{|c||c|c|c|c|c|c|}
\hline
$b\setminus w$ & $\id$ & $s_1$ & $s_2$ & $s_2 s_1$ & $s_1 s_2$ & $s_1 s_2 s_1$ \\
\hline\hline
$b\ =\ n\om_1$ & $1$ & $q^{n}$ & $1$ & $q^{n}$ & $q^{n}$ & $q^{n}$
\\ \hline
$b\!=\!s_1(b_-)$ & $0$ & $1$ & $0$ & $q^{n}$ & $1$ & $q^{n}$\\ \hline
$b=-n\om_2$ & $0$ & $0$ & $0$ & $1$ & $0$ & $1$\\ \hline
\end{tabular}
\medskip

\noindent
Note that one can pass from either of these two tables to the other
by relabeling the indices $1$ and $2$ in the columns.
\medskip

For $b_-=n_1\om_1+n_2\om_2 \ (n_1,n_2<0)$:

\begin{tabular}{|c||c|c|c|c|c|c|}
\hline
$b\setminus w$ & $\id$ & $s_1$ & $s_2$ & $s_2 s_1$ & $s_1 s_2$ 
& $s_1 s_2 s_1$ \\
\hline\hline
$b_-$ & $1$ & $q^{n_1}$ & $q^{n_2}$ & $q^{n_1+n_2}$ & $q^{n_1+n_2}$ 
& $q^{n_1+n_2}$
\\ \hline
$s_1(b_-)$ & $0$ & $1$ & $0$ & $q^{n_1}$ & $q^{n_1+n_2}$ & $0$ \\ \hline
$s_2(b_-)$ & $0$ & $0$ & $1$ & $q^{n_1+n_2}$ & $q^{n_2}$ & $0$ \\ \hline
$s_2 s_1(b_-)$ & $0$ & $0$ & $0$ & $1$ & $0$ & $q^{n_2}$ \\ \hline
$s_1 s_2(b_-)$ & $0$ & $0$ & $0$ & $0$ & $1$ & $q^{n_1}$ \\ \hline
$b_+$ & $0$ & $0$ & $0$ & $0$ & $0$ & $1$ \\ \hline
\end{tabular}
\medskip
\medskip

Furthermore, let us provide the values of $a_{b,w}=q^{n_b(w(b))}$ 
for $A_3$ in the notation from Conjecture \ref{CONJDOMINB} for 
$b=b_-$, i.e. for antidominant $\,b$. 
Confirming this conjecture,  
$\,-n_b(w(b))$ coincides with the lowest $q$\~degree of the 
coefficient of $X_{w(b)-b}$ in the product 
$\prod_{\al\in R_+}(1-qX_\al)^{-1}$ from Lusztig's definition
of the (nonaffine) Kostant $q$\~partition function. See 
\cite{JLZ} and \cite{FFL} concerning using the Kostant 
$q$\~partition function in the theory of the BK\~filtration and
the PBW\~filtration.
 
Namely,  $\ n_b(w(b))=(b,\ga_w)\,$ for 
$w(b)\succq b\in P_-$ and proper $\ga_w\in P_+$.
In this case, $\ga_w$ is the maximal positive root in the 
set $\la(w)$ ($\ga_w=0$ for $w=\id$), except for the following
permutations:
\begin{align}\label{abwa3}
&\ga_w=\ep_{12}+\ep_{34} \for w=(2143),\ 
\ga_w=\ep_{14} \for w=(3142)\notag,\\
&\ga_w=\ep_{14}+\ep_{23} \for w=(3412),\ \ (4312),
\ \ (3421),\ \ (4321),
\end{align}
where $\ep_{ij}=\ep_i-\ep_j$, $\al_i=\ep_i-\ep_{i+1}$ in the
notation from the tables of \cite{Bo}.

We note that one can compute the $\overline{E}^\dag$\~polynomials 
for $A_n$ using the SAGE software for the $E$\~polynomials based on the 
formula due to Haiman-Haglund-Loehr followed by $t\to\infty$. 
However a direct usage of the intertwining operators 
of nil-DAHA is more efficient (and we need them for all root
systems).

\comment{
We also conjecture that for $b\in B_-$,\ 
$n_b(c)=-\sum_j(\ga_j^\vee,b)\,$ where $\Th^u=\{\ga\}$ is the set of 
all maximal roots in $\la(u)$ for $u=u_c$, i.e., those
satisfying $\ga_j+\al\not\in \la(u)$
for any $\al\in R_+$. 
Following the classical theory of
maximal roots, $(\ga,\al)\ge 0$ for any $\ga\in \Th^u$ 
and all $\al\in \la(u)$. Then  $(\ga,\ga')=0$ if 
$\ga\neq \ga'\in \Th^u$,
because otherwise $\ga=(\ga-\ga')+\ga'$ provided that $\ga-\ga'>0$. 
Therefore $(\ga',\al_j)=0$ for $\ga\neq \ga'$ and all $\al_j$ 
in the presentation of $\ga$ in terms of simple roots. 
Denoting the set of such $\al_j$
by  $\Ga^\ga$, it is connected in the Dynkin diagram $\Ga$. The
subsets $\Ga^\ga$ and $\Ga^{\ga'}$ are not linked in $\Ga$ if
$\ga\neq \ga'$. 
Correspondingly, $u=\prod_\ga u^\ga$ for $\ga\in \Th^u$,
where $u^\ga$ is from the Weyl
group associated with $\Ga^\ga$ and $\ga$ is a unique maximal
root in $\la(u^\ga)$. Also, the elements $u^\ga$ are pairwise 
commutative and $\la(u)$ is a disjoint union of $\la(u^\ga)$.
This conjecture is compatible with Claim $(i)$ from
Corollary \ref{EDAGCOEFF} combined with the relations
from (\ref{daginteb}). We hope it can be deduced from
these statements.
}

\bibliographystyle{unsrt}

\begin{thebibliography} {ABCD}
\vfil




\bibitem [BeF] {BeF} 
{R.~Bezrukavnikov}, and {M.~Finkelberg},
{\em Equivariant Satake category and Kostant-Whittaker 
reduction}, Preprint arXiv:0707.3799v2 (2007).

\bibitem [BC] {BC}
{A.~Borodin}, and {I.~Corwin},
{\em Macdonald processes}, Preprint
arXiv:1111.\-4408v4 [math.PR] (2013). 

\bibitem [B] {Bo}
{N.~Bourbaki},
{\em Groupes et alg\`ebres de Lie}, Ch. {4--6},
Hermann, Paris (1969).

\bibitem [BrF] {BF}
{A.~Braverman}, and {M.~Finkelberg},
{\em Finite-difference quantum Toda lattice via equivariant
$K$-theory}, Transformation Groups { 10} (2005), 363--386.

\bibitem [CS] {CS}
{W.~Casselman}, and {J.~Shalika},
{\em The unramified principal series of $p$-adic groups, II.
The Whittaker function}, Comp. Math. {41} (1980), 207--231.

\bibitem [C1] {C13}
{I.~Cherednik},
{\em Integration of quantum many-body problems
by af\-fine Knizh\-nik--Za\-mo\-lod\-chi\-kov equations},
Pre\-print RIMS {  776} (1991),
Adv. Math. { 106} (1994), 65--95.

\bibitem [C2] {C15}
\bysame
{\em Double affine Hecke algebras,
Knizhnik-Za\-mo\-lod\-chi\-kov equa\-tions, and Mac\-do\-nald's
ope\-ra\-tors},
IMRN {9} (1992), 171--180.

\bibitem [C3] {C2}
\bysame,
{\em Double affine Hecke algebras and Macdonald's
conjectures},
Annals of Mathematics {141} (1995), 191--216.


\bibitem [C4] {C4}
\bysame,
{\em Nonsymmetric Macdonald polynomials},
IMRN {10} (1995), 483--515.

\bibitem [C5] {C5}
\bysame,
{\em Difference Macdonald-Mehta conjecture},
IMRN {10} (1997), 449--467.

\bibitem [C6] {C1}
{\bysame},
{\em Intertwining operators of double affine Hecke algebras},
Se\-lec\-ta Math. New ser. {3} (1997), 459--495

\bibitem [C7] {C12}
\bysame,
{\em Double affine Hecke algebras and difference
Fourier transforms},
Inventiones Math. {152} (2003), 213--303.

\bibitem [C8] {C101}
\bysame,
{\em Double affine Hecke algebras},
London Mathematical Society Lecture
Note Series, {319}, Cambridge University Press, Cambridge, 2006.

\bibitem [C9] {C103}
\bysame,
{\em Non-semisimple Macdonald polynomials},
Selecta Math. (2008).

\bibitem [C10] {C10}
\bysame,
{\em Whittaker limits of difference spherical functions}, 
IMRN {20} (2009), 3793--3842; arXiv:0807.2155 (2008).

\bibitem [CF] {ChF}
\bysame, and {B.~Feigin},
{\em Rogers-Ramanujan type identities and Nil-DAHA},
Preprint arXiv:1209.1978v4 [math.QA] (2012).

\bibitem [CM] {ChM}
\bysame, and {X.~Ma},
{\em Spherical and Whittaker functions via DAHA I,II},
Preprint arXiv: 0904.4324 (2009), Selecta Math.

\bibitem [CO1] {ChO1}
\bysame, and {D.~Orr},
{\em One-dimensional nil-DAHA and Whittaker functions I},
Transformation Groups {17:4} (2012), 953--987. 

\bibitem [CO2] {ChO2}
\bysame, and {D.~Orr},
{\em One-dimensional nil-DAHA and Whittaker functions II},
to appear in Transformation Groups (2013); arXiv:1104.3918v3 (2011).



\bibitem [Et] {Et1}
{P.~Etingof},
{\em  Whittaker functions on quantum groups and $q$-deformed Toda
operators}, AMS Transl. Ser. 2, {194}, 9--25, AMS, Providence,
Rhode Island, 1999.


\bibitem [FFL] {FFL}
{E.~Feigin}, and {G.~Fourier}, and {P.~Littelmann}, 
PBW–-filtration over $\Z$  and compatible bases
for $V_{\Z}(\la)$ in type $A_n$ and $C_n$,
Preprint arxiv: 1204.1854v1.

\bibitem [GLO1] {GLO1}
{A.~Gerasimov}, and {D.~Lebedev}, and {S.~Oblezin},
{\em On $q$-deformed $\mathfrak{gl}_{l+1}$-Whittaker functions, I},
Preprint arXiv: 0803.0145. 


\bibitem [GiL] {GiL}
{A.~Givental}, and {Y.-P.~Lee},
{\em Quantum $K$-theory on flag
manifolds, finite-difference Toda lattices and quantum groups},
Inventiones Math. {151} (2003), 193--219.

\bibitem [GW] {GW}
{R.~Goodman}, and {N.~R.~Wallach},
{\em Conical vectors and Whittaker vectors},
J. Functional Analysis, {39} (1980), 199--279.

\bibitem [HC] {HC}
{Harish-Chandra},
{\em Discrete series for
semisimple Lie groups, II.},
Acta Math. {116} (1963), 1--111.



\bibitem [Hu] {Hu}
{J.~Humphreys},
{\em Reflection groups and Coxeter Groups},
Cambridge University Press (1990).

\bibitem [Ion] {Ion1}
{B.~Ion},
{\em Nonsymmetric Macdonald polynomials and Demazure characters},
Duke Mathematical Journal {116}:2 (2003), 299--318.






\bibitem [JLZ] {JLZ}
{A.~Joseph}, and {G.~Letzter}, and {S.~Zelikson},
{\em On the Brylinski-Kostant filtration},
JAMS {13}:4 (2000), 945--970.

\bibitem [KS] {KnS}
{F.~Knop}, and {S.~Sahi},
{\em A recursion and a combinatorial formula for Jack
polynomials}, Inventiones Math. {128}:1 (1997), 9--22.

\bibitem [L] {L}
{G.~Lusztig},
{Affine Hecke algebras and their graded version},
J. of the AMS {2}:3 (1989), 599--635.


\bibitem [M1] {M2}
\bysame,
{\em A new class of symmetric functions},
Publ. I.R.M.A., Strasbourg, Actes 20-e Seminaire Lotharingen,
(1988), 131--171 .

\bibitem [M2] {M4}
\bysame,
{\em Affine Hecke algebras and orthogonal polynomials},
S\'e\-mi\-naire Bour\-baki {47}:797 (1995), 01--18.

\bibitem [M3] {M1}
\bysame,
{\em Affine Hecke algebras and orthogonal polynomials},
Cambridge Tracts in Mathematics {157},
Cambridge University Press, 2003.


\bibitem [Op] {O2}
{E.~Opdam},
{\em Harmonic analysis for certain
representations of graded Hecke algebras},
Acta Math. {175} (1995), 75--121.

\bibitem [Ru] {Ru}
{S.N.M.~Ruijsenaars},
{\em Factorized weight functions vs. factorized scattering},
Commun. Math. Phys. {228} (2002), 467-–494.

\bibitem [San] {San}
{Y.~Sanderson},
{\em On the Connection Between Macdonald Polynomials and
Demazure Characters},
J. of Algebraic Combinatorics, {11} (2000), 269--275.

\bibitem [Sev] {Sev}
{A. Sevostyanov}, 
{\em Quantum deformation of Whittaker modules and the Toda lattice},
Duke Math. J., { 105}:2 (2000), 211--238.



\bibitem [Sto] {Sto2}
{J.~Stokman}
{\em The c-function expansion of a basic hypergeometric function
associated to root systems},
Preprint arXiv:1109.0613 (2011).


\bibitem [Wa] {Wa}
{N.~R.~Wallach},
{\em Real Reductive Groups II},
Academic Press, Boston, 1992.

\end{thebibliography}

\end{document}